\documentclass[reqno,12pt]{amsart}

\usepackage{latexsym}
\usepackage{amsthm,amsmath, amssymb,amsfonts}
\usepackage{mathrsfs}

\usepackage{comment}
\usepackage{hyperref}
\theoremstyle{plain}
\usepackage[dvipsnames]{xcolor}
\usepackage{enumitem}
\usepackage{graphicx}

\theoremstyle{plain}
\newtheorem{thm}{Theorem}[subsection]
\newtheorem{cor}[thm]{Corollary}
\newtheorem{lem}[thm]{Lemma}
\newtheorem{prop}[thm]{Proposition}

\newtheorem*{question}{Question}
\newtheorem*{problem}{Problem}

\makeatletter
\@addtoreset{thm}{section}
\makeatother

\renewcommand{\thethm}{%
\ifnum\value{subsection}=0
\thesection.\arabic{thm}%
\else
\thesubsection.\arabic{thm}%
\fi
}

\usepackage{hyperref}
\usepackage[dvipsnames]{xcolor}
\hypersetup{colorlinks=true, linktocpage=true,citecolor=Blue, linkcolor=RedViolet, urlcolor=RedViolet}

\theoremstyle{definition}

\newtheorem{rmk}[thm]{Remark}
\newtheorem{example}[thm]{Example} 

\newcommand{\eps}{\varepsilon}
\def\Xint#1{\mathchoice
{\XXint\displaystyle\textstyle{#1}}%
{\XXint\textstyle\scriptstyle{#1}}%
{\XXint\scriptstyle\scriptscriptstyle{#1}}%
{\XXint\scriptscriptstyle\scriptscriptstyle{#1}}%
\!\int}
\def\XXint#1#2#3{{\setbox0=\hbox{$#1{#2#3}{\int}$}
	\vcenter{\hbox{$#2#3$}}\kern-.5\wd0}}

\def\dashint{\Xint-}

\newcommand\restr[2]{{
	\left.\kern-\nulldelimiterspace 
	#1 
	\right|_{#2} 
}}

\makeatletter

\def\l@section{%
\@tocline{1}{5pt}{0pt}{2.3em}{\bfseries}%
}

\def\l@subsection{%
\@tocline{2}{0pt}{1.8em}{3.2em}{\normalfont}%
}

\makeatletter
\newcommand{\subsectionnotoc}[1]{%
\subsection{#1}%
\let\subsectionnotoc@heading\@svsechd
\def\@svsechd{%
	\begingroup
	\let\addcontentsline\@gobblethree
	\subsectionnotoc@heading
	\endgroup
}%
}
\makeatother

\makeatletter
\newcommand{\subsectionstarnotoc}[1]{%
\subsection*{#1}%
\let\subsectionstarnotoc@heading\@svsechd
\def\@svsechd{%
	\begingroup
	\let\addcontentsline\@gobblethree
	\subsectionstarnotoc@heading
	\endgroup
}%
}
\makeatother
\begin{document}

\title[]{SPHERES WITH TRANSITIVE FAMILIES OF CMC
	HYPERSURFACES}

\author{Luciano L. Junior}
\address{Università di Trento, Dipartimento di Matematica\\
	Via Sommarive 14, 38123 Trento, Italy}
\email{luciano.luzzijunior@unitn.it}

\author{Ivan Miranda}
\address{IMPA -- Instituto de Matemática Pura e Aplicada\\
	Rio de Janeiro, RJ, Brasil, 22460-320}
\email{ivan.miranda@impa.br}

\begin{abstract}
	The goal of this paper is to investigate the following
	question: what geometric restrictions does the existence
	of a transitive family of constant mean curvature
	hyperspheres impose on a Riemannian metric on $\mathbb S^n$?
	The notion of a transitive family can be made precise
	through the theory of Zoll families developed by Ambrozio,
	Marques, and Neves. Within this framework, we partially
	extend their local rigidity theorem for even conformal
	factors on $\mathbb S^3$ to every dimension $n\ge4$,
	assuming each hypersphere in the minimal Zoll family has
	the area of a round equator.
	Moreover, we give a partial answer to the following
	generalized version of a question posed by Bryant: must a metric on $\mathbb S^n$, conformal to the
	round metric and admitting a Zoll family of hyperspheres
	with a common nonzero constant mean curvature, have
	constant sectional curvature?
	Finally, we investigate inverse problems following ideas of Hangan and prove a rigidity result about metrics for which the canonical geodesic spheres have constant mean curvature.
\end{abstract}

\maketitle

{
	\tableofcontents
	\thispagestyle{empty}
}

\section{Introduction}

A Riemannian metric $g$ on the sphere $\mathbb S^2$ is called Zoll if all its geodesics are simple, closed, and have the same length. Hilbert asked if every such metric had constant curvature. In~\cite{Zoll1903}, Zoll provided a negative answer by explicitly constructing examples of non-round Zoll metrics. Later, Funk \cite{Funk13} proved the following infinitesimal rigidity, for a conformal deformation of the round metric through Zoll metrics with fixed common geodesic length, the infinitesimal conformal factor must be odd. The non obstructed case was treated by Guillemin~\cite{Guillemin}, who proved that every smooth odd function can be realized by the velocity of such a deformation. The complementary rigidity statement is due to Green~\cite{Green63}, and later was reproved by LeBrun and Mason~\cite{LeBrun_Mason02}, which state that, if $g=e^{2\rho}\mathrm{can}$ is Zoll and $\rho$ is even, then $\rho$ is constant.  Here, $\mathrm{can}$ denotes the canonical round metric of the sphere.

Ambrozio, Marques, and Neves extended this picture to higher dimensions.
Following their definition, a Riemannian metric $g$ on $\mathbb S^n$
belongs to $\mathcal Z$ if there exists a smooth family
$\{\Sigma_\sigma\}_{\sigma\in\mathbb{RP}^n}$ of smoothly embedded
$(n-1)$-dimensional minimal spheres in $(\mathbb S^n,g)$ such that,
for each $x\in\mathbb S^n$ and each hyperplane
$\pi\subset T_x\mathbb S^n$, there exists a unique
$\sigma\in\mathbb{RP}^n$ with $x\in\Sigma_\sigma$ and
$T_x\Sigma_\sigma=\pi$. Such a family is called a Zoll family.
We write $g\in\mathcal Z'$ if, moreover,
$\mathrm{Area}_g(\Sigma_\sigma)
=\mathrm{Area}_{\mathrm{can}}(\mathbb S^{n-1})$
for every $\sigma\in\mathbb{RP}^n$.

For every $n\geq3$, they proved that each smooth odd function $f$
on $\mathbb S^n$ is realized as the infinitesimal conformal factor
of a smooth path $g_t=e^{2\rho_t}\mathrm{can}\in\mathcal Z'$,
with $\rho_0=0$ and $\dot\rho_0=f$%
~\cite[Theorem~A]{lucas_coda_andre}.
They also established a local analogue on $\mathbb S^3$ of
Green's rigidity theorem. More precisely, if $g=e^{2\rho}\mathrm{can}\in\mathcal Z$,
with $\rho$ even and sufficiently close to zero in the smooth
topology, then $\rho$ is constant%
~\cite[Theorem~B]{lucas_coda_andre}.

The objective of this work is to use the perspective of Ambrozio, Marques and Neves to approach the following extension of Hilbert's question, posed by Bryant
and recorded by Herreros~\cite{Herreros12}, which concerns curves
of prescribed nonzero geodesic curvature.
\begin{question}
	If all curves of a fixed nonzero geodesic curvature on a Riemannian surface are closed, must the surface have constant Gaussian curvature?
\end{question}

In higher dimensions we can reformulate Bryant's question as a question about curves of prescribed constant mean curvature (CMC). More generally, we can use the language introduced by Ambrozio, Marques, and Neves to extend this question to spheres of arbitrary dimension. To make this precise, consider the definition.

For a fixed $h\in\mathbb R$, we say that a smooth metric $g$
on $\mathbb S^n$ belongs to $\mathcal Z_h$ if there exists a
smooth family $\{\Sigma_\sigma\}_{\sigma\in\mathbb S^n}$ of embedded
oriented hyperspheres, each with constant mean curvature $h$
with respect to its chosen orientation, such that, for every $x\in\mathbb S^n$ and every
oriented hyperplane $\pi$ in $T_x\mathbb S^n$, there is a unique
$\sigma\in\mathbb S^n$ for which $x\in\Sigma_\sigma$ and
$T_x\Sigma_\sigma=\pi$ as oriented hyperplanes.
Writing $\sigma=\sigma(x,\pi)$ for this unique parameter, we further
require that the assignment $ (x,\pi)\longmapsto\Sigma_{\sigma(x,\pi)} $ be of class $C^3$. The set $\{\Sigma_\sigma\}_{\sigma \in \mathbb{S}^n}$ is called $h$-Zoll family.

With this machinery we propose the following higher-dimensional analogue of Bryant's question in the conformal class of the round metric.

\begin{question}
	Let $n\geq 2$, $h>0$, and let $\mathrm{can}$ denote the canonical
	Riemannian metric on $\mathbb S^n$.
	If $g=e^{2\rho}\mathrm{can}\in\mathcal Z_h$, then $g$ is isometric to the round metric, up to homothety.
\end{question}

Regarding partial results. In dimension two, Herreros~\cite{Herreros12} proved
rigidity under perturbations supported in a sufficiently small
region of a surface of constant curvature (see also Asselle and
Lange~\cite{Asselle_Lange20} and Asselle and
Benedetti~\cite{Asselle_Benedetti22}\footnote{These results, among others cited here, belong
to the framework of magnetic systems. At the end of the
	introduction, we provide a brief dictionary between the theory there within and the perspective of prescribed mean curvature theory
	to clarify how our results fit into the literature.}). In higher dimensions, the axiom of spheres, see Leung and Nomizu~\cite{LeungNomizu71}, cover the case where the $h$-Zoll family consists of totally umbilical submanifolds. 

Beyond these rigidity results, we prove that the proposed question holds locally in dimensions two and three, for almost every prescribed mean curvature. More specifically, we prove the following.

{
	\renewcommand{\thethm}{A}
	\begin{thm}\label{maintheoremA}
		Let $n= 2$ or $3$. There exists a subset
		$\mathcal K^n\subset(0,\infty)$ of full Lebesgue measure with the
		following property. For every $h\in\mathcal K^n$, there exists an open neighborhood
		$\mathcal V_h$ of $\mathrm{can}$ in its conformal class,
		endowed with the smooth topology, such that, if $g\in\mathcal V_h\cap\mathcal Z_h$ and $\mathrm{Vol}_g(\mathbb S^n)
		=\mathrm{Vol}_{\mathrm{can}}(\mathbb S^n)$, then
		$$
		g=\varphi^*\mathrm{can}
		$$
		for some conformal diffeomorphism
		$\varphi:\mathbb S^n\to\mathbb S^n$.
	\end{thm}
}

The proof of Theorem~\ref{maintheoremA} follows in two steps. First,
we study the infinitesimal obstructions to conformal
deformations of the canonical metric within $\mathcal Z_h$.
This leads to the study of the kernel and inverse of
a Funk--Radon type integral operator. Sobolev
estimates for its inverse are available for almost every $h>0$ with respect
to Lebesgue measure. Second, we turn this infinitesimal
information into a local rigidity statement, generalizing
the arguments of Kiyohara~\cite{Kiyohara87} and Ambrozio,
Marques, and Neves~\cite{lucas_coda_andre}. The estimates used in this step hold for every $h>0$, so the restriction to a set of full measure reduces to
the estimates in the first step. A heuristic explanation of these theorems can be found in Section~\ref{subsection:Motivation}.

Adapting ideas from the proof of Theorem~\ref{maintheoremA} we also conclude the following generalized version of the local rigidity theorem proved by Ambrozio, Marques, and Neves on $\mathbb S^3$.

{
	\renewcommand{\thethm}{B}
	\begin{thm}\label{maintheoremB}
		Fix $n \geq 4$. The canonical metric in $\mathbb{S}^n$  has a open neighborhood $\mathcal U$ in the smooth topology such that if $\rho$ is even and $e^{2\rho} \mathrm{can} \in \mathcal U \cap \mathcal Z'$, then $\rho=0$.
	\end{thm}
}

In contrast with the minimal case, the uniformization theorem implies that, in dimension two, Theorem~\ref{maintheoremA} prevents volume-normalized deformations of the round metric within $\mathcal{Z}_h$ for almost every $h>0$. In higher dimensions, however, conformal directions no longer exhaust the possible deformation directions at the canonical metric. To better understand deformations beyond the conformal class, we also study the following inverse problem: given a fixed family of hyperspheres, characterize all Riemannian metrics for which every member of the family has constant mean curvature $h>0$.

In the minimal case, Hangan~\cite{Hangan96} classified
the Riemannian metrics on $\mathbb R^n$ for which all
hyperplanes are minimal. The corresponding problem
on the sphere was studied by Ambrozio, Marques, and
Neves~\cite{lucas_coda_andre}. See~\cite{lucas_survey26} for a survey on this and related problems. Here, we consider the following CMC counterpart of this problem. For $r\in(0,\pi/2)$ and $\sigma\in\mathbb S^n$, let
$$
\Sigma_\sigma^r
\doteq
\{p\in\mathbb S^n:d_{\mathrm{can}}(p,\sigma)=r\},
$$
be oriented by the normal pointing into the canonical
geodesic ball of center $\sigma$ and radius $r$.

\begin{problem}
	Given $r\in(0,\pi/2)$ and $h>0$, classify
	the Riemannian metrics $g$ on $\mathbb S^n$ for which
	every $\Sigma_\sigma^r$, $\sigma\in\mathbb S^n$, has constant
	mean curvature $h$.
\end{problem}

As highlighted before, from the axiom of spheres, the totally umbilical case is rigid in dimensions $n\geq3$. In dimension two, the rigidity theorem for magnetic systems of Burns and Matveev~\cite{BurnsMatveev06} implies that any such metric
is a constant multiple of the canonical one. Moreover, their argument contemplate the stronger case where $h$ is allowed to be a positive smooth function. 

In these notes, we extend the surface rigidity result of Burns-Matveev to the hypersurface case in arbitrary dimensions. Moreover, since our
argument is essentially local and relies only on
properties that all space forms enjoy, we obtain
the following stronger rigidity theorem.

{
	\renewcommand{\thethm}{C}
	\begin{thm}\label{maintheoremC}
		Let $(M^n,\bar g)$, $n\geq 2$, be a connected Riemannian manifold of
		constant sectional curvature $\kappa\in\mathbb R$, such that
		$$
		r^*
		:=
		\inf_{p\in M}\mathrm{inj}_{\bar g}(p)
		>0.
		$$
		Fix $0<r<r^*$ (if $\kappa>0$, assume also $r<{\pi}/{2\sqrt{\kappa}}$). Let $g$ be a Riemannian metric on $M$. Suppose that there exists a smooth function 
		$f: M \to (0,\infty)$ such that, for every $\sigma\in M$, the $\bar g$-geodesic sphere
		of radius $r$ centered at $\sigma$,
		$$
		\Sigma_\sigma^r
		\doteq
		\left\{
		p\in M:
		d_{\bar g}(p,\sigma)=r
		\right\},
		$$
		has prescribed mean curvature $\restr{f}{\Sigma_{v}^r}$ in $(M,g)$. Then $g$ is homothetic to $\bar g$ and $f$ is constant. 
	\end{thm}
}

To conclude the introduction, we complement Bryant's question with further motivation to study Zoll and generalized Zoll metrics. We first discuss how, in dimensions two and three, these metrics arise as rigidity cases in certain variational inequalities on the sphere. In dimension two, we also describe their dynamical interpretation within the theory of magnetic systems.

The Three Geodesics Theorem of Lusternik and Schnirelmann states that every Riemannian two-sphere contains at least three distinct simple closed
geodesics~\cite{LusternikSchnirelmann47}. The variational principle underlying the existence of three distinct geodesics lies in the topology of the space of unparametrized embedded circles in $\mathbb S^2$ deformation retracts onto $\mathbb{R}P^2$. Applying min--max methods
to the length functional gives three positive critical values
$$
\ell_1(g)\leq\ell_2(g)\leq\ell_3(g),
$$
each realized by a simple closed geodesic. Mazzucchelli
and Suhr~\cite{MazzucchelliSuhr18} proved that these three
values coincide if and only if $g$ is Zoll.

Hatcher's proof of the Smale conjecture gives a parallel
picture in dimension three: the space of unparametrized
embedded two-spheres in $\mathbb S^3$ deformation retracts
onto the family of great spheres, parametrized by
$\mathbb{RP}^3$~\cite{Hatcher83}. Again, this topology suggests
that the area functional should have at least four
critical points, motivating the classical Yau's conjecture that every
Riemannian metric on the three-sphere contains at least four distinct
embedded minimal two-spheres~\cite{Yau82}.
The corresponding min--max theory in the Simon--Smith
setting gives four spherical area widths
$$
\sigma_1(g)\leq\sigma_2(g)\leq\sigma_3(g)\leq\sigma_4(g).
$$
Ambrozio, Marques, and Neves~\cite{AmbrozioMarquesNeves24}
proved that these four values coincide if and only if
$g$ admits a Zoll family of embedded minimal two-spheres.
Thus, generalizing the two-dimensional scenario. 

The Lusternik--Schnirelmann theorem and Yau's conjecture have
natural counterparts in the CMC setting. For a given
Riemannian metric on $\mathbb S^n$, with $n=2$ or $3$,
Arnold's conjecture and its three-dimensional analogue
predict the existence of at least two distinct immersed
hyperspheres with constant mean curvature $h$, for each
$h>0$ (see
\cite[Conjectures~1.3 and~1.6]{MazurowskiZhou25}).
As before, a topological motivation comes from the fact
that the space of oriented, unparametrized, smoothly
embedded $(n-1)$-spheres in $\mathbb S^n$ is homotopy
equivalent to $\mathbb S^n$.
Recent progress (see for instance  \cite{ChengZhou21,ChengZhou23,SarnataroStrykerCurves,MazurowskiZhou25}) suggests that variational methods provide
a promising approach to these questions. It is therefore
natural to conjecture that generalized Zoll metrics
again arise as the rigidity case of a suitable variational
inequality. In view of Theorem~A, we further expect
that this rigidity within the conformal class occurs only for metrics of constant curvature, in contrast to the minimal case.

In dimension two, generalized Zoll metrics also have a natural
dynamical interpretation through the theory of magnetic systems.
We briefly recall the correspondence between prescribed curvature
and magnetic systems.

A magnetic system on the oriented sphere $\mathbb S^2$ is a pair
$(g,f)$, where $g$ is a Riemannian metric and
$f\colon\mathbb S^2\to\mathbb R$ is a smooth function.
A $(g,f)$-geodesic is a unit-speed curve
$\gamma\colon I\to\mathbb S^2$ whose geodesic curvature is
$f\circ\gamma.$ Existence and uniqueness for this equation define a smooth flow
on the unit tangent bundle.
This flow describes the motion of a charged particle on
$\mathbb S^2$ subject to the magnetic field $f\,dV_g$ (see~\cite{Asselle_Lange20}).

A magnetic system is called Zoll if its flow has the same
orbits as a free $S^1$-action on the unit tangent
bundle. When $f\equiv0$, this recovers the usual notion of a Zoll metric. Moreover, every metric $g \in \mathcal Z_h$ determines a Zoll magnetic system with constant magnetic function $f \equiv h$. Conversely, a Zoll magnetic system
$(g,h)$ provides a metric in $\mathcal Z_h$ given that its prime
orbits project to embedded curves. In particular, for each fixed $h>0$, these notions coincide for metrics in a neighborhood of the round metric.

With this small dictionary, the aforementioned works regarding magnetic systems can be translated to the notion of prescribed mean curvature. In particular, we record the rigidity theorem of Asselle and Benedetti~\cite[Theorem~1.9]{Asselle_Benedetti22} which implies that, if a metric $g$ on $\mathbb S^2$ satisfies
$
g\in\bigcap_j\mathcal Z_{h_j}
\;\text{for some sequence }h_j\rightarrow+\infty,
$
then $g$ has constant Gaussian curvature.

Finally, we record the work of Bohr and Paternain~\cite{BohrPaternain26} on the existence of families of Zoll magnetic systems on $\mathbb S^2$ near $(\mathrm{can},0)$. Their construction allows the magnetic function to vary, in particular does not guarantee the constancy of the magnetic field. These examples are therefore compatible with the rigidity results proved here.

\subsectionstarnotoc{Organization of the Paper}
In Section~\ref{prelim}, we introduce the notation and establish
the basic formulae for normal graphical perturbations of geodesic
spheres in $\mathbb S^n$. We conclude by showing that, under the
hypotheses of Theorems~\ref{maintheoremA} and~\ref{maintheoremB},
the corresponding Zoll families are unique up to reparametrization
and, for metrics sufficiently close to the round metric, are
represented by the normal graphs provided by White's perturbation
theory. The representation of Zoll families by normal graphs,
after reparametrization, is the only step in our argument
that requires the dimensional assumption of
Theorem~\ref{maintheoremA}. This step can be interpreted
as a local uniqueness and graphical representation for
Zoll families associated with a fixed metric.
In dimension two, the required uniqueness follows
directly from uniqueness for the ODE's governing curves of constant geodesic curvature.
In dimension three, it follows from the uniqueness
theorem of Gálvez and Mira for transitive families
\cite[Theorem~2.4]{GalvezMira20}.
Consequently, a higher-dimensional extension of this
local uniqueness and representation result would suffice
to extend the rigidity argument.

Section~\ref{section:maintheorem} contains the proofs of
Theorems~\ref{maintheoremA} and~\ref{maintheoremB}.
In Section~\ref{subsection:Motivation}, we give a heuristic
explanation of the two main steps: establishing infinitesimal
rigidity and converting into local rigidity. We then state the
required estimates and prove both theorems assuming these estimates.
In Section~\ref{subsection:funk-radon}, we study the Funk--Radon
type integral operator underlining the infinitesimal rigidity, more specifically we study bounds for its inverse in various Sobolev spaces. In
Section~\ref{subsection:nonlinearbound}, we establish the nonlinear
estimates needed to pass from infinitesimal to local rigidity. 

Finally, in Section~4, we address the inverse problem and prove
Theorem~\ref{maintheoremC}.

\subsectionstarnotoc{Acknowledgements}  The authors are grateful to Lucas Ambrozio for his interest in this work and his valuable suggestions on earlier drafts of the paper.

\subsectionstarnotoc{Funding} The first author was supported by the European Research Council (ERC) under the European Union’s Horizon 2020 research and innovation programme (grant agreement No. 947923). The second author was supported in part by the Coordenação de Aperfeiçoamento de Pessoal de Nível Superior – Brasil (CAPES) – Finance Code 001.

\subsectionstarnotoc{AI Disclaimer} 
We have used Large Language Models (LLMs) to enrich the bibliography and revise the final manuscript. The LLMs did not produce any text or mathematical content present in this paper.
\section{Graphical perturbations}\label{prelim}

The purpose of this section is to fix notation and recall a few fundamental formulae for graphical perturbations of hypersurfaces in $\mathbb S^n$. The theory derived in this section is inspired by~\cite{lucas_coda_andre}.  

\subsectionnotoc{Notation}

Let $\mathrm{can}$ denote the canonical metric on $\mathbb S^n$. Fix $r\in(0,\pi/2]$. For each $\sigma\in\mathbb S^n$, define
$$
\Sigma_\sigma^r
=
\{x\in\mathbb S^n:d_{\mathrm{can}}(x,\sigma)=r\}
=
\{x\in\mathbb S^n:\langle x,\sigma\rangle=\cos r\}
$$
and
$$
\Omega_\sigma^r
=
\{x\in\mathbb S^n:d_{\mathrm{can}}(x,\sigma)\leq r\}.
$$
The outward unit normal to $\Sigma_\sigma^r$ can be expressed as
$$
N_\sigma^r(x)
=
\frac{\cos r\,x-\sigma}{\sin r}.
$$
We use this normal to parametrize graphical perturbations. The orientation of $\Sigma_\sigma^r$ is determined by the unit normal pointing into $\Omega_\sigma^r$, and its mean curvature is always computed with respect to this orientation.

\begin{lem}
	The family $\{\Sigma_\sigma^r\}_{\sigma\in\mathbb S^n}$ consists of embedded hyperspheres with constant mean curvature
	$$
	h(r)=(n-1)\cot r.
	$$
	In particular, it is an $h(r)$-Zoll family for the canonical metric.
\end{lem}

The incidence set of the canonical family is
\begin{align*}
	\mathcal I_r
	&=
	\{(x,\sigma)\in\mathbb S^n\times\mathbb S^n:x\in\Sigma_\sigma^r\}\\
	&=
	\{(x,\sigma)\in\mathbb S^n\times\mathbb S^n:d_{\mathrm{can}}(x,\sigma)=r\}.
\end{align*}
This set allows us to describe graphical perturbations of all the canonical spheres simultaneously. Let $\Phi\in C^\infty(\mathcal I_r)$ take values in $(-r,\pi-r)$, and write $\Phi_\sigma(x)=\Phi(x,\sigma)$. Define the perturbed hypersphere
\begin{equation}\label{graph}
	\Sigma_\sigma^r(\Phi)
	=
	\left\{
	\cos(\Phi_\sigma(x))x
	+
	\sin(\Phi_\sigma(x))N_\sigma^r(x):
	x\in\Sigma_\sigma^r
	\right\}.
\end{equation}
Clearly, $\Sigma_\sigma^r(0)=\Sigma_\sigma^r$. It will be also convenient to introduce the map $E_r:\mathcal I_r\times[-r,\pi-r]\longrightarrow\mathbb S^n$, defined by
$$
E_r(x,\sigma,s)
=
\cos(s)x+\sin(s)N_\sigma^r(x).
$$
We will also use $\Sigma_\sigma^r(\Phi)$ to denote the parametrizing map $\Sigma_\sigma^r(\Phi):\Sigma_\sigma^r\longrightarrow\mathbb S^n$, given by
$$
\Sigma_\sigma^r(\Phi)(x)
=
E_r(x,\sigma,\Phi_\sigma(x)).
$$
We think of $\Sigma_\sigma^r(\Phi)$ as the normal graph over $\Sigma_\sigma^r$ of the function $\Phi_\sigma$.

If $\Phi$ is sufficiently small in $C^1(\mathcal I_r)$, these maps are embeddings for every $\sigma\in\mathbb S^n$. We denote by $\Omega_\sigma^r(\Phi)$ the connected component of $\mathbb S^n\setminus\Sigma_\sigma^r(\Phi)$ containing $\sigma$, so that $\Sigma_\sigma^r(\Phi)=\partial\Omega_\sigma^r(\Phi)$. Let $N_\sigma^r(\Phi)$ denote the outward pointing normal of $\Omega^r_\sigma(\Phi)$ with respect to the canonical metric. By our convention, we orient $\Sigma_\sigma^r(\Phi)$ by the inward pointing normal, and its mean curvature will always be computed with respect to this orientation.

We next collect some elementary formulae describing the geometry of these graphs. Set
$$
j_r(s)
=
\frac{\sin(r+s)}{\sin r}
$$
and define
$$
(\Sigma_\sigma^r)^\perp(\Phi)(x)
=
-\sin(\Phi_\sigma(x))x
+
\cos(\Phi_\sigma(x))N_\sigma^r(x).
$$

\begin{lem}\label{lembasicformulae}
	For every $x\in\Sigma_\sigma^r$ and $X,Y\in T_x\Sigma_\sigma^r$, the following identities hold:
	\begin{itemize}
		\item[$i)$] The differential of the graph map is
		$$
		d\Sigma_\sigma^r(\Phi)_x[X]
		=
		j_r(\Phi_\sigma)X
		+
		d\Phi_\sigma(X)(\Sigma_\sigma^r)^\perp(\Phi).
		$$
		\item[$ii)$] The induced metric is
		$$
		\left\langle
		d\Sigma_\sigma^r(\Phi)[X],
		d\Sigma_\sigma^r(\Phi)[Y]
		\right\rangle
		=
		j_r(\Phi_\sigma)^2\langle X,Y\rangle
		+
		d\Phi_\sigma(X)d\Phi_\sigma(Y).
		$$
		\item[$iii)$] The outward unit normal to $\Sigma_\sigma^r(\Phi)$ is
		$$
		N_\sigma^r(\Phi)
		=
		\frac{
			j_r(\Phi_\sigma)(\Sigma_\sigma^r)^\perp(\Phi)
			-
			\nabla^{\Sigma_\sigma^r}\Phi_\sigma
		}{
			\sqrt{
				j_r(\Phi_\sigma)^2
				+
				\left|\nabla^{\Sigma_\sigma^r}\Phi_\sigma\right|^2
			}
		}.
		$$
		\item[$iv)$] The tangential Jacobian is
		$$
		\mathrm{Jac}\Sigma_\sigma^r(\Phi)
		=
		j_r(\Phi_\sigma)^{n-2}
		\sqrt{
			j_r(\Phi_\sigma)^2
			+
			\left|\nabla^{\Sigma_\sigma^r}\Phi_\sigma\right|^2
		}.
		$$
	\end{itemize}
	Here, $\nabla^{\Sigma_\sigma^r}$, the corresponding norm, and all inner products on $\Sigma_\sigma^r$ are computed with respect to the metric induced by $\mathrm{can}$.
\end{lem}

\begin{proof}
	Let $X$ be a tangent vector to $\Sigma_\sigma^r$ at $x$. Differentiating the graph map gives
	$$
	d\Sigma_\sigma^r(\Phi)_x[X]
	=
	j_r(\Phi_\sigma(x))X
	+
	d\Phi_\sigma(X)(\Sigma_\sigma^r)^\perp(\Phi)(x),
	$$
	which proves $i)$.
	
	The vector $(\Sigma_\sigma^r)^\perp(\Phi)(x)$ has unit length and is orthogonal to the tangent space of $\Sigma_\sigma^r$ at $x$. Thus, for tangent vectors $X$ and $Y$,
	$$
	\left\langle
	d\Sigma_\sigma^r(\Phi)_x[X],
	d\Sigma_\sigma^r(\Phi)_x[Y]
	\right\rangle
	=
	j_r(\Phi_\sigma)^2\langle X,Y\rangle
	+
	d\Phi_\sigma(X)d\Phi_\sigma(Y),
	$$
	proving $ii)$.
	
	For $iii)$, set
	$$
	W
	=
	j_r(\Phi_\sigma)(\Sigma_\sigma^r)^\perp(\Phi)
	-
	\nabla^{\Sigma_\sigma^r}\Phi_\sigma.
	$$
	
	Since $(\Sigma_\sigma^r)^\perp(\Phi)$ is orthogonal to the gradient of $\Phi_\sigma$, we have
	$$
	|W|^2
	=
	j_r(\Phi_\sigma)^2
	+
	\left|\nabla^{\Sigma_\sigma^r}\Phi_\sigma\right|^2.
	$$
	
	Moreover, both terms defining $W$ are orthogonal to the position vector $\Sigma_\sigma^r(\Phi)(x)$. Therefore, $W$ is tangent to the ambient sphere at that point. For every tangent vector $X$,
	\begin{align*}
		\left\langle
		W,d\Sigma_\sigma^r(\Phi)_x[X]
		\right\rangle
		&=
		j_r(\Phi_\sigma)
		\left[
		d\Phi_\sigma(X)
		-
		\left\langle
		\nabla^{\Sigma_\sigma^r}\Phi_\sigma,X
		\right\rangle
		\right] \\
		&=0.
	\end{align*}
	
	Thus, $W$ is normal to the graph. Since
	$$
	\left\langle
	W,(\Sigma_\sigma^r)^\perp(\Phi)
	\right\rangle
	=
	j_r(\Phi_\sigma)>0,
	$$
	it points outward. Dividing by its norm proves $iii)$.
	
	Finally, let $e_1,\ldots,e_{n-1}$ be an orthonormal basis of the tangent space of $\Sigma_\sigma^r$ at $x$, and set
	$$
	u_i
	=
	j_r(\Phi_\sigma)e_i
	+
	d\Phi_\sigma(e_i)(\Sigma_\sigma^r)^\perp(\Phi).
	$$
	Then
	$$
	g_{ij}
	=
	\langle u_i,u_j\rangle
	=
	j_r(\Phi_\sigma)^2\delta_{ij}
	+
	d\Phi_\sigma(e_i)d\Phi_\sigma(e_j).
	$$
	
	The determinant of this rank-one perturbation is
	$$
	\det(g_{ij})
	=
	j_r(\Phi_\sigma)^{2(n-2)}
	\left(
	j_r(\Phi_\sigma)^2
	+
	\left|\nabla^{\Sigma_\sigma^r}\Phi_\sigma\right|^2
	\right).
	$$
	Therefore proving $iv)$.
\end{proof}

The generalized Gauss map of the perturbed family is the map $\mathcal G_r(\Phi):\mathcal I_r\longrightarrow T_1\mathbb S^n$ defined by
$$
\mathcal G_r(\Phi)(x,\sigma)
=
\bigl(
\Sigma_\sigma^r(\Phi)(x),
N_\sigma^r(\Phi)(x)
\bigr).
$$

\begin{prop}
	If $\Phi$ is sufficiently close to zero in $C^2(\mathcal I_r)$, then $\mathcal G_r(\Phi)$ is a diffeomorphism.
\end{prop}

\begin{proof}
	Clearly the map $\mathcal G_r(0)$ is a diffeomorphism. Moreover, the map $\Phi\mapsto\mathcal G_r(\Phi)$ is continuous from the $C^2$ topology to the $C^1$ topology. Hence, for $\Phi$ sufficiently close to zero in $C^2(\mathcal I_r)$, the map $\mathcal G_r(\Phi)$ is $C^1$-close to $\mathcal G_r(0)$, and therefore the conclusion follows.
\end{proof}

\subsectionnotoc{Weighted area functionals}

Given a Riemannian metric  $g$ on $\mathbb S^n$ and $h\in\mathbb R$ the $h$-area
functional acting on smooth domains $\Omega\subset\mathbb S^n$ is defined by
$$
\mathcal A_h(g,\Omega)
=
\mathrm{area}_g(\partial\Omega)
-
h\,\mathrm{vol}_g(\Omega).
$$
We notice that $\Omega$ is a critical point of this functional if and only if its boundary is $h$-CMC. Indeed, if $\Omega_t$ is a variation of $\Omega$ with outward normal velocity $f$, then
$$
\left.\frac{d}{dt}\right|_{t=0}
\mathcal A_h(g,\Omega_t)
=
\int_{\partial\Omega}(H_g-h)f\,dA_g.
$$

Fix $r\in(0,\pi/2]$ and let $g=e^{2\rho}\mathrm{can}$. For every $\Phi$
sufficiently close to zero, the weighted area functional defines a function
on $\mathbb S^n$ by
$$
\mathcal A(\rho,\Phi)(\sigma)
=
\mathcal A_{h(r)}
\bigl(g,\Omega_\sigma^r(\Phi)\bigr).
$$
For $(x,\sigma)\in\mathcal I_r$, $s\in(-r,\pi-r)$, and
$p\in T_x\Sigma_\sigma^r$, set
$$
A_\rho((x,\sigma),s,p)
=
e^{(n-1)\rho(E_r(x,\sigma,s))}
j_r(s)^{n-2}
\sqrt{j_r(s)^2+|p|^2}
$$
and
$$
V_\rho((x,\sigma),s)
=
\int_{-r}^{s}
e^{n\rho(E_r(x,\sigma,t))}
j_r(t)^{n-1}\,dt.
$$
The formula for the Jacobian given in the previous section, together with the coarea formula gives the following expression
\begin{equation}\label{eqformulaA}
	\mathcal A(\rho,\Phi)(\sigma)
	=
	\int_{\Sigma_\sigma^r}
	\left[
	A_\rho
	\left(
	(x,\sigma),
	\Phi_\sigma(x),
	\nabla^{\Sigma_\sigma^r}\Phi_\sigma(x)
	\right)
	-
	h(r)V_\rho
	\left(
	(x,\sigma),
	\Phi_\sigma(x)
	\right)
	\right]
	dA_{\mathrm{can}}(x).
\end{equation}

We now study the Euler--Lagrange equation of this functional. To this end, we
define $\mathcal H(\rho,\Phi)$ as the unique smooth function on $\mathcal I_r$
such that, for every $\varphi\in C^\infty(\mathcal I_r)$,
\begin{equation}\label{eqD2A}
	\left.\frac{d}{dt}\right|_{t=0}
	\mathcal A(\rho,\Phi+t\varphi)(\sigma)
	=
	\int_{\Sigma_\sigma^r}
	\mathcal H(\rho,\Phi)(x,\sigma)
	\varphi_\sigma(x)\,
	dA_{\mathrm{can}}(x).
\end{equation}
Differentiating the expression in~\eqref{eqformulaA} and isolating $\varphi$ via integration by parts, we obtain
\begin{equation}\label{eqH}
	\begin{aligned}
		\mathcal H(\rho,\Phi)(x,\sigma)
		={}&
		-\mathrm{div}_{\Sigma_\sigma^r}
		\left[
		D_3A_\rho
		\left(
		(x,\sigma),
		\Phi_\sigma(x),
		\nabla^{\Sigma_\sigma^r}\Phi_\sigma(x)
		\right)
		\right] \\
		&+
		D_2A_\rho
		\left(
		(x,\sigma),
		\Phi_\sigma(x),
		\nabla^{\Sigma_\sigma^r}\Phi_\sigma(x)
		\right) \\
		&-
		h(r)D_2V_\rho
		\left(
		(x,\sigma),
		\Phi_\sigma(x)
		\right).
	\end{aligned}
\end{equation}
Here $D_2$ and $D_3$ denote differentiation with respect to $s$ and $p$,
respectively.

The operator $\mathcal H$ has the following geometric interpretation. Let
$ H(\rho,\Phi)(x,\sigma)$ denote the scalar mean curvature of
$\Sigma_\sigma^r(\Phi)$ in the metric $e^{2\rho}\mathrm{can}$, evaluated at
$\Sigma_\sigma^r(\Phi)(x)$. The first variation formula, together with equation~\eqref{eqformulaA}, gives
$$
\mathcal H(\rho,\Phi)(x,\sigma)
=
\left(
H(\rho,\Phi)(x,\sigma)-h(r)
\right)
e^{n\rho(\Sigma_\sigma^r(\Phi)(x))}
j_r(\Phi_\sigma(x))^{n-1}.
$$
In partiular, we obtain the following result.
\begin{lem}
	For each $\sigma\in\mathbb S^n$, the hypersphere
	$\Sigma_\sigma^r(\Phi)$ has constant mean curvature $h(r)$ if and only if
	$$
	\mathcal H(\rho,\Phi)(\cdot,\sigma)=0.
	$$
	Equivalently, the whole family consists of $h(r)$-CMC hyperspheres if and only
	if $\mathcal H(\rho,\Phi)=0$ on $\mathcal I_r$.
\end{lem}

The conformal transformation law for the mean-curvature vector yields the
following transformation formula for the Euler--Lagrange operator
$\mathcal H$.

\begin{lem}\label{confchangelawEL}
	For every $f\in C^\infty(\mathbb S^n)$ and $\sigma\in\mathbb S^n$,
	\begin{align*}
		\mathcal H(\rho+f,\Phi)(\cdot,\sigma)
		={}&
		e^{(n-1)f\circ\Sigma_\sigma^r(\Phi)}
		\mathcal H(\rho,\Phi)(\cdot,\sigma)
		\\
		&+
		h(r)j_r(\Phi_\sigma)^{n-1}
		e^{n\rho\circ\Sigma_\sigma^r(\Phi)}
		e^{(n-1)f\circ\Sigma_\sigma^r(\Phi)}
		\left(
		1-e^{f\circ\Sigma_\sigma^r(\Phi)}
		\right)
		\\
		&+
		(n-1)j_r(\Phi_\sigma)^{n-1}
		e^{(n-1)\rho\circ\Sigma_\sigma^r(\Phi)}
		e^{(n-1)f\circ\Sigma_\sigma^r(\Phi)}
		N_\sigma^r(\Phi)(f).
	\end{align*}
\end{lem}

\begin{proof}
	By the geometric identity above,
	$$
	\mathcal H(\rho,\Phi)(\cdot,\sigma)
	=
	\left(
	H(\rho,\Phi)(\cdot,\sigma)-h(r)
	\right)
	e^{n\rho\circ\Sigma_\sigma^r(\Phi)}
	j_r(\Phi_\sigma)^{n-1}.
	$$
	The conformal transformation law for the mean curvature is
	$$
	H(\rho+f,\Phi)(\cdot,\sigma)
	=
	e^{-f\circ\Sigma_\sigma^r(\Phi)}
	\left(
	H(\rho,\Phi)(\cdot,\sigma)
	+
	(n-1)e^{-\rho\circ\Sigma_\sigma^r(\Phi)}
	N_\sigma^r(\Phi)(f)
	\right).
	$$
	Substituting this formula into the preceding identity gives the result.
\end{proof}

Consequently, the study of $\mathcal H(\rho,\Phi)$ reduces to the case
$\rho=0$, for which we have the following explicit formula.

\begin{prop}
	For every $\Phi\in C^\infty(\mathcal I_r)$ sufficiently close to zero and
	every $\sigma\in\mathbb S^n$,
	\begin{equation}\label{H(0,phi)}
		\begin{aligned}
			\mathcal H(0,\Phi)(\cdot,\sigma)
			={}&
			j_r(\Phi_\sigma)^{n-3}j_r'(\Phi_\sigma)
			\frac{
				(n-1)j_r(\Phi_\sigma)^2
				+
				(n-2)
				\left|
				\nabla^{\Sigma_\sigma^r}\Phi_\sigma
				\right|^2
			}{
				\left(
				j_r(\Phi_\sigma)^2
				+
				\left|
				\nabla^{\Sigma_\sigma^r}\Phi_\sigma
				\right|^2
				\right)^{1/2}
			}
			\\
			&-
			\mathrm{div}_{\Sigma_\sigma^r}
			\left(
			\frac{
				j_r(\Phi_\sigma)^{n-2}
				\nabla^{\Sigma_\sigma^r}\Phi_\sigma
			}{
				\left(
				j_r(\Phi_\sigma)^2
				+
				\left|
				\nabla^{\Sigma_\sigma^r}\Phi_\sigma
				\right|^2
				\right)^{1/2}
			}
			\right)
			\\
			&-
			h(r)j_r(\Phi_\sigma)^{n-1}.
		\end{aligned}
	\end{equation}
\end{prop}

\subsectionnotoc{Representation of Local Perturbations}

The objective of this section is to study uniqueness of Zoll families for metrics in a neighborhood of the round metric. We specialize in two cases: First when $r \in (0, \pi/2)$ and $n=2$ or $3$, and second $r=\pi/2$, $n \geq 3$ and each element of the Zoll family has volume normalize equal to $\mathrm{Vol}_{\mathrm{can}}(\mathbb{S}^{n-1})$. More precisely, we prove the following proposition.  

\begin{prop}\label{prop:localperturbation}
	Assume that either $n=2$ or $3$ and $0<r<\pi/2$,
	or $n\geq2$ and $r=\pi/2$.
	For every integer $k\geq2$, there exist
	a neighborhood $\mathscr V_r$ of 0 in $C^\infty(\mathbb S^n)$, in the smooth topology,
	and a $C^k$ map
	$
	\psi:\mathscr V_r\longrightarrow C^k(\mathcal I_r),$ with
	$\psi(0)=0,$ such that the following holds. If
	$\rho\in\mathscr V_r\cap C^\infty(\mathbb S^n)$ satisfies
	\begin{itemize}
		\item $e^{2\rho}\mathrm{can}\in\mathcal Z_{h(r)}$,
		if $0<r<\pi/2$, or;
		
		\item $e^{2\rho}\mathrm{can}\in\mathcal Z'$,
		if $r=\pi/2$,
	\end{itemize}
	then every $h(r)$-Zoll family whose members, when $r=\pi/2$,
	are additionally assumed to have $e^{2\rho}\mathrm{can}$-area equal to
	$\mathrm{Area}_{\mathrm{can}}(\mathbb S^{n-1})$
	is given, after reparametrization of $\mathbb S^n$, by
	$\left\{\Sigma_\sigma^r\bigl(\psi(\rho)_\sigma\bigr)\right\}_{\sigma \in \mathbb{S}^n}$.
	Moreover, in the minimal case, the function $\psi(\rho)$ satisfies
	$\psi(\rho)(x,-\sigma)=-\psi(\rho)(x,\sigma)$.
\end{prop}

For the strictly CMC case this proposition follows by combining White's perturbation theorem
\cite[Theorem~3.2 and Section~7]{BWhite91} with the fact that every
$h(r)$-CMC hypersphere belongs to the Zoll family, for $n=2,3$. In dimension two,
this follows directly from the classical existence and uniqueness
theorem for the constant geodesic curvature ODE, while in dimension
three it is a consequence of Gálvez--Mira's completeness theorem for
transitive families \cite[Theorem~2.4]{GalvezMira20}. In the minimal case the proposition follows by the combination of White's perturbation theorem and a version of Allard’s regularity theorem.

For $j\geq 2$ and $0<\alpha<1$, let
$$
\mathscr E^{j,\alpha,+}
=
\mathrm{Emb}^{j,\alpha}(\mathbb S^{n-1},\mathbb S^n)
\big/
\mathrm{Diff}^{j,\alpha}_{+}(\mathbb S^{n-1}),
$$
denote the space of oriented unparametrized $C^{j,\alpha}$-embedded hyperspheres. We endow this space with the quotient topology, which can be show to be a Hausdorff space (for a more ample analysis of this space see \cite{AliasPiccione11}). The orientation of
$\Sigma\in\mathscr E^{j,\alpha,+}$ uniquely determine an connected $C^{j,\alpha}$ domain
$\Omega_\Sigma$ for which $\partial\Omega_\Sigma=\Sigma$ and
satisfies our orientation convention. Set
$L_r=\{\Sigma_\sigma^r:\sigma\in\mathbb S^n\}$. For each
$\sigma\in\mathbb S^n$, the hypersphere $\Sigma_\sigma^r$ is a
critical point of $\mathcal A_{h(r)}$ with nullity $n$, generated by
the perturbations within the family $\mathbb S^n$. Therefore, for $r \in (0,\pi/2)$, White's perturbation theorem can be stated as follows. The minor modifications needed to cover the minimal case are straightforward and left to the reader.

\begin{prop}[White's perturbation theorem]
	\label{prop:whiteperturbation}
	Fix $\alpha\in(0,1)$. For every integer $k\geq2$, there exist an integer $q\geq2$,
	a neighborhood $\mathcal V$ of $\mathrm{can}$ in
	$\mathrm{Met}^q(\mathbb S^n)$, a neighborhood $\mathcal U$
	of $L_r$ in $\mathscr E^{k,\alpha,+}$, and a $C^k$ map
	$
	\psi:
	\mathcal V
	\longrightarrow
	C^{k+1}(\mathcal I_r)
	$
	with the following properties:
	\begin{enumerate}
		\item $\psi(\mathrm{can})=0$;
		
		\item For every $g\in\mathcal V$, the map
		$
		\sigma
		\mapsto
		\Sigma_\sigma^r\bigl(\psi(g)_\sigma\bigr)
		$
		is an topological embedding of $\mathbb S^n$ into $\mathcal U$.
		
		\item For every $g\in\mathcal V$, define
		$$
		a_g(\sigma)
		=
		\mathcal A_{h(r)}
		\left(
		g,
		\Omega_\sigma^r\bigl(\psi(g)_\sigma\bigr)
		\right).
		$$
		Then $a_g$ is of class $C^k$, and
		$\Sigma_\sigma^r\bigl(\psi(g)_\sigma\bigr)$ is $h(r)$-CMC if
		and only if $\sigma$ is a critical point of $a_g$. Moreover,
		every oriented $h(r)$-CMC hypersphere in $\mathcal U$ is of
		this form for a unique $\sigma\in\mathbb S^n$.
	\end{enumerate}
\end{prop}

\begin{proof}[Proof of Proposition~\ref{prop:localperturbation}]
	Fix $k\geq2$, and choose $q,\mathcal V,\mathcal U$, and the 
	graphical map $\psi_{\mathrm W}$ as in 
	Proposition~\ref{prop:whiteperturbation}. 
	
	\noindent\underline{The CMC case:} Set 
	$$
	\mathscr V
	=
	\left\{
	\rho\in C^q(\mathbb S^n):
	e^{2\rho}\mathrm{can}\in\mathcal V
	\right\},
	\qquad
	\psi(\rho)
	=
	\psi_{\mathrm W}\bigl(e^{2\rho}\mathrm{can}\bigr).
	$$
	
	Let $\rho\in\mathscr V\cap C^\infty(\mathbb S^n)$, set
	$g=e^{2\rho}\mathrm{can}$, and let
	$\{N_\sigma\}_{\sigma\in\mathbb S^n}$ be an arbitrary
	$h(r)$-Zoll family for $g$. For each $\mu\in\mathbb S^n$, set
	$$
	\Sigma_\mu
	=
	\Sigma_\mu^r\bigl(\psi(\rho)_\mu\bigr),
	\qquad
	W_g
	=
	\{\Sigma_\mu:\mu\in\mathbb S^n\}.
	$$
	We prove that
	$$
	C_g
	=
	\left\{
	\sigma\in\mathbb S^n:
	N_\sigma=\Sigma_\mu
	\text{ for some }\mu\in\mathbb S^n
	\right\}
	$$
	is open, closed, and nonempty.
	
	We first show that $C_g$ is nonempty. Since $\mathbb S^n$ is compact,
	$a_g$ has a critical point $\mu_0$. Hence,
	Proposition~\ref{prop:whiteperturbation} implies that
	$\Sigma_{\mu_0}$ is $h(r)$-CMC. Since $g$ is smooth, elliptic
	regularity shows that this hypersphere is smooth.
	
	Suppose first that $n=2$. Fix a point of $\Sigma_{\mu_0}$ and let
	$N_{\sigma_0}$ be the unique member of the Zoll family with the same
	oriented tangent line at that point. Both curves solve the same
	constant geodesic curvature ODE with the same initial data. Uniqueness
	therefore gives
	$
	\Sigma_{\mu_0}=N_{\sigma_0}
	$. If $n=3$, the same conclusion follows from the completeness theorem of
	Gálvez--Mira \cite[Theorem~2.4]{GalvezMira20}. Thus, in either case,
	$\Sigma_{\mu_0}=N_{\sigma_0}$ for some
	$\sigma_0\in\mathbb S^n$, and hence $C_g\neq\varnothing$.
	
	The set $C_g$ is open. Indeed, if $\sigma_0\in C_g$, then
	$N_{\sigma_0}\in W_g\subset\mathcal U$. The smooth dependence of the
	Zoll family on its parameter implies that
	$N_\sigma\in\mathcal U$ whenever $\sigma$ is sufficiently close to
	$\sigma_0$. Since every $N_\sigma$ is $h(r)$-CMC, item~(3) of
	Proposition~\ref{prop:whiteperturbation} gives
	$N_\sigma\in W_g$, and therefore $\sigma\in C_g$.
	
	The set $C_g$ is also closed. Indeed, $W_g$ is compact and therefore
	closed in the Hausdorff space $\mathscr E^{k,\alpha,+}$. Since
	$$
	C_g
	=
	\left\{
	\sigma\in\mathbb S^n:N_\sigma\in W_g
	\right\},
	$$
	and the map $\sigma\mapsto N_\sigma$ is continuous, it follows that
	$C_g$ is closed. Consequently,
	$
	C_g=\mathbb S^n,
	$
	and, after reparametrization, the given Zoll family coincides with
	$\{\Sigma_\mu\}_{\mu\in\mathbb S^n}$, as desired.
	\\ 
	
	\noindent \underline{The minimal case:} Given a Riemannian metric $g\in\mathrm{Met}(\mathbb S^n)$, let
	$\mathcal S_g\subset \mathscr E^{k,\alpha,+}$ denote the set of oriented $g$-minimal
	hyperspheres in $\mathbb S^n$ with $g$-area equal to
	$\mathrm{Area}_{\mathrm{can}}(\mathbb S^{n-1})$. We prove the stronger statement that, after shrinking
	$\mathcal V \cap C^{\infty}(\mathbb{S}^n)$ if necessary, $\mathcal S_g\subset\mathcal U$, for every $g\in\mathcal V$. 
	
	The prove goes by contradiction. Suppose that there exist a sequence
	of metrics $g_\ell\to\mathrm{can}$ in the smooth topology and a sequence
	$$
	\Sigma_\ell\in\mathcal S_{g_\ell},
	\qquad
	\Sigma_\ell\notin\mathcal U .
	$$
	Since
	$
	\mathrm{area}_{g_\ell}(\Sigma_\ell)
	=
	\omega_{n-1}
	$
	and $g_\ell\to\mathrm{can}$ smoothly, we obtain $\sup_{\ell}\mathrm{Area}_{\mathrm{can}}(\Sigma_\ell) < \infty$. Moreover, since $H_{(\Sigma_\ell,g_\ell)}=0$, by the compactness theorem for integral varifolds
	\cite[Theorem~6.4]{Allard72} (see also \cite[Proposition 3]{ASong23}), after passing to a subsequence we obtain an integral
	stationary $(n-1)$-varifold $V$ in
	$(\mathbb S^n,\mathrm{can})$ such that
	$
	\Sigma_\ell\longrightarrow V
	$
	in the varifold sense. In particular,
	$$
	\|V\|(\mathbb S^n)=\omega_{n-1}.
	$$
	
	We claim that $V$ is a multiplicity-one equator. Indeed, by Allard's regularity theorem \cite[Theorem~8.19]{Allard72}, since $\|V\|(\mathbb S^n)=\omega_{n-1}$, the Euclidean cone over $V$ is smooth, therefore $V$ is a multiplicity-one equator, which we denote by $E$ (see also \cite[Claim 5.7]{Choi_Mantoulidis22}). 
	
	Since each $\Sigma_\ell$ is $g_\ell$-minimal and
	$g_\ell\to\mathrm{can}$ smoothly, a standard application of the
	monotonicity formula for minimal hypersurfaces, with constants
	uniform in $\ell$, shows that the varifold convergence
	$\Sigma_\ell\to|E|$ implies convergence of the supports in
	the Hausdorff topology.
	
	We may now apply White's version of Allard's regularity theorem
	\cite[Theorem~5.1]{White16}. To obtain that $\Sigma_\ell$ converge in the smooth graph topology to $E$ with multiplicity one. This
	contradicts $\Sigma_\ell\notin\mathcal U$ and proves the claim.
\end{proof}


\section{The Conformal rigidity for transitive families of CMC spheres }\label{section:maintheorem}

This section is organized as follows. In Section~\ref{subsection:Motivation}, we give a heuristic reduction of the proofs of Theorems~\ref{maintheoremA} and~\ref{maintheoremB} to bounds for the first and second variations of the functional
\begin{align*}
	\mathcal{A}_r: C^{\infty}(\mathbb{S}^n) \times C^{\infty}(\mathcal{I}_r)
	&\to C^{\infty}(\mathbb{S}^n), \\
	(\rho,\Phi)
	&\mapsto \mathcal A_r(e^{2\rho} \mathrm{can},\Phi).
\end{align*}
Still in Section~\ref{subsection:Motivation}, we state these bounds precisely and prove Theorems~\ref{maintheoremA} and~\ref{maintheoremB} assuming them. The proofs of the aforementioned bounds are then given in Sections~\ref{subsection:funk-radon} and~\ref{subsection:nonlinearbound}.

\subsection{Motivation and main arguments}\label{subsection:Motivation}
\subsubsection{\textbf{The infinitesimal rigidity for CMC-Zoll deformations}}
In order to motivate the proof of Theorem A, fix $n=2,3$, $r\in(0,\pi/2)$, set $h(r)=(n-1)\cot r$, and consider a smooth path $\rho_t\in C^\infty(\mathbb S^n)$ with $\rho_0=0$ and $\dot\rho_0=f$, such that $g_t=e^{2\rho_t}\mathrm{can}$ belongs to $\mathcal Z_{h(r)}$ for every sufficiently small $t$.

By the representation theorem of the preceding section, after reparametrizing each family, there exists a smooth path $\psi(\rho_t)\in C^\infty(\mathcal I_r)$ such that the $h(r)$-Zoll family associated with $g_t$ is given by $\{\Sigma_\sigma^r\bigl(\psi(\rho_t)_\sigma\bigr)\}_{\sigma \in \mathbb{S}^n}$.

From the $h(r)-CMC$ condition the function $\mathcal A\bigl(\rho_t,\psi(\rho_t)\bigr)$ is constant in $\mathbb{S}^n$, in particular its derivatives depend only on the time parameter. On the other hand, the time derivative of $\mathcal A\bigl(\rho_t,\psi(\rho_t)\bigr)$ can be computed, for each $\sigma \in \mathbb{S}^n$, as follows:
$$
\left.\frac{d}{dt}\right|_{t=0}
\mathcal A\bigl(\rho_t,\psi(\rho_t)\bigr)(\sigma)
=
\left.\frac{d}{dt}\right|_{t=0}
\mathcal A(\rho_t,0)(\sigma)
+
\left.\frac{d}{ds}\right|_{s=0}
\mathcal A\bigl(0,\psi(\rho_s)\bigr)(\sigma).
$$
The second term vanishes because the canonical family is critical for the weighted area functional. Moreover, the first term can be computed as
\begin{equation}\label{linofA}
	\left.\frac{d}{dt}\right|_{t=0}
	\mathcal A\bigl(\rho_t,\psi(\rho_t)\bigr)(\sigma)
	=
	(n-1)\left[
	\int_{\Sigma_\sigma^r}f\,dA_{\mathrm{can}}
	-
	\frac{n}{n-1}h(r)
	\int_{\Omega_\sigma^r}f\,dV_{\mathrm{can}}
	\right].
\end{equation}
Integrating \eqref{linofA} over $\sigma$ and applying the integral-geometric identities from Appendix~\ref{app:geometric.integral.formula}, we obtain
$$
\left.\frac{d}{dt}\right|_{t=0}
\mathcal A\bigl(\rho_t,\psi(\rho_t)\bigr)
=
(n-1)\,
\frac{
	\mathrm{area}_{\mathrm{can}}(\Sigma_\sigma^r)
	-\frac{n}{n-1} h(r)\mathrm{vol}_{\mathrm{can}}(\Omega_\sigma^r)
}{
	\mathrm{vol}_{\mathrm{can}}(\mathbb S^n)
}
\int_{\mathbb S^n}f\,dV_{\mathrm{can}}.
$$

Substituting this value back into equation~\eqref{linofA} and writing
$
\mathring f
=
f-\dashint_{\mathbb S^n}f\,dV_{\mathrm{can}}
$ for the mean-zero part of $f$, the equation can be rewritten as: 
$$
\int_{\Sigma_\sigma^r}\mathring f\,dA_{\mathrm{can}}
-
\frac{n}{n-1} h(r)
\int_{\Omega_\sigma^r}\mathring f\,dV_{\mathrm{can}}
=
0,
\qquad
\text{for every }\sigma\in\mathbb S^n.
$$
This identity is a necessary infinitesimal condition on the velocity of any smooth path through $\mathrm{can}$ in $\mathcal Z_{h(r)}$.

The relevant integral operator is a weighted combination of two classical integral-type operators: the spherical transform associated with the radius $r$, obtained by integrating over $\Sigma_\sigma^r$, and the spherical cap transform, obtained by integrating over $\Omega_\sigma^r$. For a general literature regarding integral transforms see \cite{Rubin00}, \cite{HelgasonBook1999} and \cite{Strichartz81}.

In the minimal equatorial setting, a similar linear obstruction is already present in the work of Ambrozio, Marques, and Neves. Indeed, the linearization of the area functional is given by the Funk transform described there. This description is essential to their analysis, since it allows them to classify conformal directions into those that can be integrated into Zoll deformations and those that are obstructed, see \cite[Theorem~A, Theorem~B and Lemma~A.1]{lucas_coda_andre}. This approach goes back to the works of Funk and Guillemin \cite{Funk13,Guillemin}.

Guided by this analogy, we define the weighted Funk--Radon transform associated with $r$ as the following integral operator:
$$
R_r(u)(\sigma)
=
\int_{\Sigma_\sigma^r}u\,dA_{\mathrm{can}}
-
\frac{n}{n-1}h(r)
\int_{\Omega_\sigma^r}u\,dV_{\mathrm{can}}.
$$

The identity $R_r(\mathring f)=0$ does not yet distinguish genuine
infinitesimal deformations from those generated by conformal
diffeomorphisms of the round sphere. In fact, since
$\mathcal Z_{h(r)}$ is invariant under pullback by diffeomorphisms, a
path of conformal transformations through the identity produces a path
in $\mathcal Z_{h(r)}$ whose infinitesimal conformal factor is a linear
function. Thus, the degree-one spherical harmonics correspond to the
trivial directions in the kernel of $R_r$.

In other words, denoting the space of degree-one spherical harmonics by
$$
\mathcal H^1(\mathbb S^n)
=
\left\{
\ell_w:\mathbb S^n\longrightarrow\mathbb R:
\ell_w(x)=\langle x,w\rangle,\quad w\in\mathbb R^{n+1}
\right\}.
$$
It follows that $\mathcal H^1(\mathbb S^{n})\subset\ker R_r$.

Let $f^\perp$ denote the $L^2$-orthogonal projection of $f$ onto
$
\left(
\mathcal H^0(\mathbb S^n)\oplus
\mathcal H^1(\mathbb S^n)
\right)^\perp,
$ where $\mathcal H^0(\mathbb S^n)$ denotes the space of constant functions. The preceding variational identity then implies that the velocity of any smooth path through $\mathrm{can}$ in $\mathcal Z_{h(r)}$ satisfies
\begin{equation}\label{R_rfperp}
	R_r(f^\perp)=0.
\end{equation}
The infinitesimal rigidity problem underlying Theorem~\ref{maintheoremA} therefore reduces to understanding the kernel of $R_r$ on the orthogonal complement of the constant and linear functions.

This formulation in terms of spherical harmonics is useful
because $R_r$ acts by scalar multiplication on each
spherical harmonic subspace. Indeed, the invariance of the family
$\{(\Omega_\sigma^r,\Sigma_\sigma^r)\}_{\sigma\in\mathbb S^n}$ under
the natural action of $\mathrm{SO}(n+1)$ implies that $R_r$ is
$\mathrm{SO}(n+1)$-equivariant. Therefore, our operator can be diagonalized via the spherical harmonic decomposition:
$$
L^2(\mathbb S^n)
=
\widehat{\bigoplus}_{k\geq0}\mathcal H^k(\mathbb S^n).
$$
In fact, since each $\mathcal H^k(\mathbb S^n)$ is an irreducible
$\mathrm{SO}(n+1)$-representation, Schur's lemma implies that, for every
$k\geq0$, there exists a constant $\mu_k(r)\in\mathbb R$ such that
\begin{equation}\label{structuralconstants}
	R_r(f)=\mu_k(r)f
	\qquad
	\text{for every }f\in\mathcal H^k(\mathbb S^n).
\end{equation}
For each fixed $r$, the constants $\{\mu_k(r)\}_{k\geq0}$ are called
the spectral multipliers of $R_r$. We observe that each constant $\mu_k(r)$ depend on the dimension of the sphere, however we omit this dependence because the dimension should be clear from the context. In particular, the previous argument is equivalent to the claim that $\mu_1(r)=0$, for every $r$ and dimension.

Once this diagonalization is established, the study of $R_r$ reduces
to the study of its spectral multipliers $\mu_k(r)$. Injectivity, surjectivity, continuity between Sobolev spaces, and
coercivity (that is, continuity of the inverse) are all translated
into conditions on the vanishing and asymptotic analysis of the constants
$\mu_k(r)$.

A direct computation shows that, for every $k\geq2$, the constants
$\mu_k(r)$ are explicitly given by
\begin{equation*}
	{\mu}_k(r)
	=
	\frac{-\omega_{n-1}(n+1)(n-1)
		\sin^{n+1}r}{
		\lambda_k(\mathbb S^n)\binom{n+k-2}{k}
	}
	C_{k-2}^{\frac{n+3}{2}}(\cos r),
\end{equation*}
see Lemma~\ref{lem:gegenbauer-representation}. Here $C_k^\lambda$
denotes the family of Gegenbauer polynomials, also known as
ultraspherical polynomials, which we will introduce in the next section with more details. We refer to \cite{Szego75,Vilenkin68} for a more extensive discussion. 

Although
there is an extensive literature devoted to these polynomials, to the
best of the authors' knowledge, no complete explicit classification of
the zeros relevant to this problem is known, even in low dimensions;
see \cite{BilykMastrianniSteinerberger24}. With this in mind, we approach the equation~\eqref{R_rfperp} from a qualitative perspective. For that we introduce the following notation:

\begin{equation}\label{Ginj}
	\mathcal G^n_{\mathrm{inj}}
	=
	\left\{
	r\in\left(0,\frac{\pi}{2}\right):
	R_r\text{ is injective on }
	\bigl(
	\mathcal H^0(\mathbb S^n)
	\oplus
	\mathcal H^1(\mathbb S^n)
	\bigr)^\perp
	\right\}.
\end{equation}

For each fixed $k\geq 2$, define
\begin{equation}\label{ZkandZ}
	Z_k
	=
	\left\{
	r\in\left(0,\frac{\pi}{2}\right):\mu_k(r)=0
	\right\},
	\qquad
	Z=\bigcup_{k\geq2}Z_k.
\end{equation}

The change of variable $s=\cos r$ transforms the equation $\mu_k(r)=0$ into the vanishing of a nonzero Gegenbauer polynomial. Thus each $Z_k$ is finite and $Z$ is at most countable, and $\mathcal G^n_{\mathrm{inj}}$ is residual and has full measure in $(0,\pi/2)$.

The conclusion of this discussion is the following infinitesimal rigidity for CMC perturbations of transitive families. 

\begin{lem}\label{lem:infinitesimal-rigidity}
	Let $n=2$ or $3$, and $r\in\mathcal G_{\mathrm{inj}}^n$. Suppose that
	$\rho_t$ is a smooth path in $C^\infty(\mathbb S^n)$ satisfying
	$\rho_0=0$ and
	$
	e^{2\rho_t}\mathrm{can}\in\mathcal Z_{h(r)}
	$
	for every sufficiently small $t$. Then
	$
	\dot\rho_0
	\in
	\mathcal H^0(\mathbb S^n)\oplus\mathcal H^1(\mathbb S^n).
	$
\end{lem}

The previous lemma is the infinitesimal counterpart of Theorem~\ref{maintheoremA}. The local statement does not follow at this point, since we could have a non trivial sequence of metrics $g_\ell \in \mathcal Z_{h(r)}$ converging to $\mathrm {can}$ that does not integrate into a smooth path. To avoid the existence of such non trivial sequence we will need stronger properties of the inverse of our integral operator.     
\subsubsection{\textbf{The local rigidity for CMC Zoll deformations}}
In order to improve the preceding infinitesimal rigidity result to the local rigidity statement of Theorem A, we follow a strategy developed by Kiyohara in the context of classical Zoll metrics~\cite{Kiyohara87} and by Ambrozio, Marques, and Neves in the context of minimal Zoll families~\cite{lucas_coda_andre}. Its main difficulty is concentrated in the fact that injectivity of the linearization of the weighted area functional gives no quantitative control on its inverse. In particular, even after removing the natural symmetries, it does not prevent the existence of sequences of distinct metrics $g_k\in\mathcal Z_{h(r)}$ converging to $\mathrm{can}$, as previously mentioned. In the following exposition we explain the mechanism to bypass this issue.

For this purpose, consider a metric
$g=e^{2\rho}\mathrm{can}\in\mathcal Z_{h(r)}$
sufficiently close to $\mathrm{can}$, and let $\psi(\rho)$ be the graphical representation of its Zoll family given by Whites' implicit function theorem (see Proposition~\ref{prop:localperturbation}). Formally, the Taylor expansion to first order with quadratic remainder around the canonical family is
$$
\mathcal A\bigl(\rho,\psi(\rho)\bigr)(\sigma)
=
\mathcal A(0,0)
+
(n-1)R_r(\rho)(\sigma)
+
\mathcal Q_r\bigl(\rho,\psi(\rho)\bigr)(\sigma).
$$
Here $\mathcal Q_r(\rho,\Phi)$ denotes a remainder that vanishes to second order in the pair $(\rho,\Phi)$ at $(0,0)$. 

The Zoll condition implies that
$\mathcal A\bigl(\rho,\psi(\rho)\bigr)$
is constant as a function of $\sigma$. Let $\rho^\perp$ denote the
$L^2$-orthogonal projection of $\rho$ onto
$
\left(
\mathcal H^0(\mathbb S^n)
\oplus
\mathcal H^1(\mathbb S^n)
\right)^\perp.
$
Subtracting the average of the preceding expansion and using
$\mathcal H^1(\mathbb S^n)\subset\ker R_r$ gives
$$
R_r(\rho^\perp)
=
-\frac1{n-1}
\left(
\mathcal Q_r\bigl(\rho,\psi(\rho)\bigr)
-
\dashint_{\mathbb S^n}
\mathcal Q_r\bigl(\rho,\psi(\rho)\bigr)(\tau)\,
dV_{\mathrm{can}}(\tau)
\right).
$$

Therefore, if local rigidity were to fail, there would exist a sequence
$\rho_k\to0$ such that
$e^{2\rho_k}\mathrm{can}\in\mathcal Z_{h(r)}$
and $R_r(\rho_k^\perp)$ decays at the quadratic rate prescribed by the
remainder
$\mathcal Q_r\bigl(\rho_k,\psi(\rho_k)\bigr)$. A priori, this does not contradict the injectivity of $R_r$. To exclude
this behavior we combine a coercive estimate for $R_r$, which forces a decay of linear rate, with a compatible
upper bound for $\mathcal Q_r$. Here compatible means that the norm in which the remainder have quadratic decay must be tamed by the norm in which the inverse of $R_r$ is continuous.

We make this discussion precise by collecting the estimates for the
weighted Funk--Radon transform proved in
Section~\ref{subsection:funk-radon} and the estimates for the quadratic error proved in Section~\ref{subsection:nonlinearbound}.

For $\alpha\in\mathbb R$, define:
\begin{equation}\label{Gcoerc}
	\mathcal G_\alpha^n
	=
	\left\{
	r\in\left(0,\frac{\pi}{2}\right):
	\begin{gathered}
		\exists\,C_{r,n,\alpha}>0\text{ such that}\\
		\|f\|_{H^s(\mathbb S^n)}
		\leq
		C_{r,n,\alpha}
		\|R_r(f)\|_{H^{s+\alpha}(\mathbb S^n)},\\
		\forall\,s\in\mathbb R,\quad
		\forall\,f\in C^\infty(\mathbb S^n)
		\text{ with}\\
		f\perp
		\bigl(
		\mathcal H^0(\mathbb S^n)
		\oplus
		\mathcal H^1(\mathbb S^n)
		\bigr)
	\end{gathered}
	\right\}.
\end{equation}
Thus, $\mathcal G_\alpha^n$ is the set of radii for which $R_r$ is $\alpha$-coercive on the orthogonal
complement of
$\mathcal H^0(\mathbb S^n)\oplus\mathcal H^1(\mathbb S^n)$, for any Sobolev order $s \in \mathbb{R}$. In the next section we provide a description for the set $\mathcal{G}^n_\alpha$ in terms of the asymptotical behavior of the spectral multipliers associated to $R_r$.

Combining the results of the following papers 
\cite{Rubin00,BrandoliniGariboldiGigante25,HanRahimi26} on the study of zeros and asymptotic behavior of the family of Gegenbauer polynomials we obtain the following proposition.

\begin{prop}\label{prop:coercivity-main-theorem}
	Let $n\geq2$.
	
	\begin{enumerate}
		\item If $\alpha=(n+1)/2$, then $\mathcal G_\alpha^n$ has
		Lebesgue measure zero. Moreover, if $r/\pi\in\mathbb Q$, then
		$$
		r\in\mathcal G_{\frac{n+1}{2}}^n
		\quad\Longleftrightarrow\quad
		r\in\mathcal G_{\mathrm{inj}}^n.
		$$
		
		\item If $\alpha>(n+1)/2$, then $\mathcal G_\alpha^n$ has full
		Lebesgue measure in $(0,\pi/2)$.
	\end{enumerate}
\end{prop}
\begin{proof}
	See Proposition~\ref{prop:coercivity.regimes}.
\end{proof}
In the next section, we prove this proposition and study the sets
$\mathcal G_\alpha^n$ in more detail for the different values of
$\alpha$. Next we turn to the estimate of the quadratic remainder $\mathcal Q_r$, which can be summarized in the following proposition.   
\begin{prop}\label{estimativaprincipal}
	Fix $r\in(0,\pi/2)$ and $n=2$ or $3$. There exist a neighborhood $\mathscr V_r$ of
	$\mathrm{can}$ in the space of Riemannian metrics of $\mathbb{S}^n$ endowed with the smooth topology and a
	constant $C_r>0$ with the following property. For every
	$\rho\in C^\infty(\mathbb S^n)$ such that
	$e^{2\rho}\mathrm{can}\in\mathscr V_r \cap \mathcal Z_{h(r)}$, then:
	\begin{equation}\label{w24boundeq}
		\|R_r(\mathring\rho)\|_{W^{1,2}(\mathbb S^n)}
		\leq
		C_r\|\rho\|_{W^{2,4}(\mathbb S^n)}^2,
	\end{equation}
	where $\mathring{\rho}$ denote the mean-zero part of $\rho$.
\end{prop}
\begin{proof}
	Combine Theorem~\ref{thm:nonlinear-estimate} and Remark~\ref{rmk:nonlinear-estimate}.
\end{proof}
Schematically the proof is given as follows. Applying Taylor expansion to the defining function $A$ of the operator $\mathcal A$ (see
equation~\eqref{eqformulaA}) we can derive the dependence of the quadratic error $\mathcal{Q}_r$ in terms of the graphical defining function $\psi(\rho)$ and its derivatives. Which reduces the work to estimate derivatives of $\psi(\rho)$ in terms of derivatives of $\rho$.

These bounds are provided by the CMC equation $\mathcal H\bigl(\rho,\psi(\rho)\bigr)=0$. More precisely, the conformal transformation
law for the Euler--Lagrange operator $\mathcal H$ rewrites this condition
as an equation for $\mathcal H\bigl(0,\psi(\rho)\bigr)$. Applying again Taylor in the graphical variable we obtain a equation of the form:
$$
\mathcal J^r_\sigma \psi(\rho)_\sigma
=
-D_\rho\mathcal H(0,0)_\sigma[\rho]
+
\mathcal E_\sigma\bigl(\rho,\psi(\rho)\bigr),
$$
Where $\mathcal J^r_\sigma$ is the Jacobi operator of $\Sigma^r_\sigma$ and $\mathcal E_\sigma$ is quadratic. Finally, elliptic
estimates for the equation above and a standard absorption argument,
gives an estimate for $\psi(\rho)$ in terms of $\rho$.

Combining Proposition~\ref{estimativaprincipal},
Proposition~\ref{prop:coercivity-main-theorem}, we can now prove
Theorem~\ref{maintheoremA}. The strategy of the proof follows the
interpolation and absorption argument introduced by Kiyohara in his
study of the rigidity of Zoll deformations \cite{Kiyohara87}, adapted
here to the weighted Funk--Radon transform.

\begin{proof}[Proof of Theorem~\ref{maintheoremA}]
	Let $n=2$ or $3$ and set
	$
	\mathcal G^n=\mathcal G_{\frac{n+3}{2}}^n.
	$
	By Proposition~\ref{prop:coercivity-main-theorem}, the set
	$\mathcal G^n$ has full Lebesgue measure in $(0,\pi/2)$. Fix
	$r\in\mathcal G^n$, and $k\geq 8$ and let $\mathscr{V}_r \subset C^\infty(\mathbb{S}^n)$ be a neighborhood of zero such that the graphical perturbation $\psi(\rho)$ constructed in Proposition~\ref{prop:localperturbation} is $C^k$-continuous. 
	
	First, we recall the well know $(n+1)$-parameter family of conformal
	diffeomorphisms of $\mathbb S^n$, indexed by the open unit ball $\mathbb B^{n+1} = \{a \in \mathbb{R}^{n+1} : |a|<1 \}$. Following \cite{Schoen_Yau94} Section \S8,  for each
	$a\in\mathbb B^{n+1}$ the following is a conformal transformation
	$$
	\varphi_a(x)
	=
	\frac{
		(1-|a|^2)x+2(1+\langle a,x\rangle)a
	}{
		1+2\langle a,x\rangle+|a|^2
	},
	\qquad
	x\in\mathbb S^n.
	$$
	Moreover, $\varphi_0=\mathrm{id}$ and
	$$
	\varphi_a^*\mathrm{can}
	=
	e^{2\omega_a}\mathrm{can},
	\qquad
	\omega_a(x)
	=
	\log
	\left(
	\frac{1-|a|^2}
	{1+2\langle a,x\rangle+|a|^2}
	\right).
	$$
	Consequently, if $g=e^{2\rho}\mathrm{can}$, then
	$$
	\varphi_a^*g
	=
	e^{2\Lambda(a,\rho)}\mathrm{can},
	\qquad
	\Lambda(a,\rho)
	=
	\rho\circ\varphi_a+\omega_a.
	$$
	
	Let $\Pi_1$ denote the $L^2$-orthogonal projection onto
	$\mathcal H^1(\mathbb S^n)$ and define
	$$
	\mathcal B(a,\rho)
	=
	\Pi_1\bigl(\rho\circ\varphi_a+\omega_a\bigr).
	$$
	From the construction we have $\mathcal B(0,0)=0$ and that
	$
	D_a\mathcal B(0,0):
	\mathbb R^{n+1}\to\mathcal H^1(\mathbb S^n)
	$
	is an isomorphism. The inverse function theorem allows us to shrink the  neighborhood $\mathscr{V}_r$ of $0$ in order to obtain a smooth map $\rho\mapsto a(\rho)$ such that
	$$
	\varphi_{a(\rho)}^*
	\bigl(e^{2\rho}\mathrm{can}\bigr)
	=
	e^{2\lambda(\rho)}\mathrm{can},
	\qquad
	\lambda(\rho) \doteq \Lambda(a(\rho),\rho) 
	\in
	\mathcal H^1(\mathbb S^n)^{\perp}.
	$$
	
	We now define the part of $\lambda(\rho)$ with zero mean
	value by
	$$
	u(\rho)
	=
	\Lambda(a(\rho),\rho)
	-
	\dashint_{\mathbb S^n}
	\Lambda(a(\rho),\rho)\,dV_{\mathrm{can}}, \qquad u(\rho) \perp \mathcal{H}^0(\mathbb{S}^n) \oplus \mathcal{H}^1(\mathbb{S}^n).
	$$ 
	And notice that the mean of $\lambda(\rho)$ and $u(\rho)$ are related by the following formula:
	\begin{equation*}
		\dashint_{\mathbb S^n} \lambda(\rho)\, dV_{\mathrm{can}} = - \frac{1}{n} \log\left[ \frac{\mathrm{vol}_\mathrm{can}(\mathbb{S}^n)}{\mathrm{vol}_g(\mathbb{S}^n)}\dashint_{\mathbb{S}^n} e^{n\,u(\rho)}\,dV_{\mathrm{can}} \right]
	\end{equation*}
	Using the fact that $u(\rho)$ has zero mean and shrinking the neighborhood $\mathscr{V}_r$ if necessary, we can assume the existence of a constant $C>0$, such that
	$$\left|\dashint_{\mathbb S^n} \lambda(\rho)\, dV_{\mathrm{can}} \right| \leq C \left(\,\left|\log \left[ \frac{\mathrm{vol}_\mathrm{can}(\mathbb{S}^n)}{\mathrm{vol}_g(\mathbb{S}^n)}\right]\right| + \|u(\rho)\|_{L^2}^2\right). $$  
	
	Therefore, combining Proposition~\ref{prop:coercivity-main-theorem}, applied with
	$s=-(n+1)/2$ and $\alpha=(n+3)/2$, and
	Proposition~\ref{estimativaprincipal}, with the volume normalization hypotheses $ {\mathrm{vol}_\mathrm{can}(\mathbb{S}^n)}={\mathrm{vol}_g(\mathbb{S}^n)}$, we have
	$$
	\|u(\rho)\|_{H^{-\frac{n+1}{2}}}
	\leq
	C_r\|R_r(u(\rho))\|_{H^1}
	\leq
	C_r\|\lambda(\rho)\|_{W^{2,4}}^2
	\leq
	C_r\|u(\rho)\|_{W^{2,4}}^2.
	$$
	For a positive constant $C_r>0$. Recalling that $n\leq 3$  Theorem~4 and Remarks~11 and~13 of \cite{Triebel86}, together with the
	Sobolev embedding theorem, provide the following Sobolev interpolation inequality: 
	$$\|u(\rho)\|_{W^{2,4}}^2 \lesssim \|u(\rho)\|_{H^{-\frac{n+1}{2}}} \|u(\rho)\|_{C^8}.$$
	Therefore, shrinking $\mathscr V_r$ if necessary, we obtain $u(\rho)=0$
	whenever $\rho\in\mathscr V_r$ and
	$g=e^{2\rho}\mathrm{can}\in\mathcal Z_{h(r)}$ satisfies the volume
	normalization. Consequently, $\lambda(\rho)=0$ and
	$g=\bigl(\varphi_{a(\rho)}^{-1}\bigr)^*\mathrm{can}$, as desired. 
	
\end{proof}

\subsubsection{\textbf{The local rigidity for Minimal Zoll deformations}}
The proof of Theorem~\ref{maintheoremB} follows the same argument as the
proof of Theorem~\ref{maintheoremA}. We quickly discuss the counterparts of Propositions~\ref{prop:coercivity-main-theorem} and~\ref{estimativaprincipal} for the minimal case. 

For $r=\pi/2$ the integral operator $R_r$ reduces to the classical Funk transform
$$
\mathcal F(f)(v)
:=
\int_{\Sigma_v^{\pi/2}} f\,dA_{\mathrm{can}}.
$$
Which has the following coercive estimate, see for instance \cite[Lemma~4.3]{Strichartz81}. 

\begin{prop}\label{prop:funk-coercivity}
	Let $\mathcal F$ denote the Funk transform on $\mathbb S^n$.
	
	\begin{enumerate}
		\item The kernel of $\mathcal F$ on $C^\infty(\mathbb S^n)$ is precisely
		the space of odd functions.
		
		\item For every $s\in\mathbb R$, there exists a constant $C=C(n,s)>0$ such that, 
		$$
		C\|f\|_{H^s(\mathbb S^n)}
		\leq
		\|\mathcal F(f)\|_{H^{s+\frac{n-1}{2}}(\mathbb S^n)}
		$$
		for every $f\in C_{\mathrm{even}}^\infty(\mathbb{S}^n)$.
	\end{enumerate}
\end{prop}

Regarding the counterpart of Proposition~\ref{estimativaprincipal}, we first observe that the dimensional restriction there is imposed uniquely to ensure the existence of a graphical parametrization of the Zoll family near the canonical one. When $r=\pi/2$ and $n\geq3$, Proposition~\ref{prop:localperturbation} provides such a parametrization under the assumption that $g\in\mathcal Z'$. Moreover, the same normalization allows us to work directly with $\rho$, rather than restricting to its mean-zero part in equation~\ref{w24boundeq}. The proof of the previous claim is given in Theorem~\ref{thm:nonlinear-estimate}.

\begin{proof}[Proof of Theorem~\ref{maintheoremB}]
	Combining Proposition~\ref{prop:funk-coercivity}, with $s=(3-n)/2$, and Theorem~\ref{thm:nonlinear-estimate}, together with the normalization defining $\mathcal Z'$, there exists an open neighborhood $\mathcal U$ of the canonical metric such that, if $\rho$ is even and $e^{2\rho}\mathrm{can}\in\mathcal U\cap\mathcal Z'$, then for a constant $C>0$
	$$
	C^{-1}\|\rho\|_{H^{\frac{3-n}{2}}}
	\leq
	\|\mathcal F(\rho)\|_{H^1}
	\leq C
	\|\rho\|_{W^{2,4}}^2
	\leq
	C\|\rho\|_{H^{\frac{3-n}{2}}}\,
	\|\rho\|_{H^{n+\frac52}}.
	$$
	The last inequality follows from the interpolation given in Theorem~4 and Remarks~11 and~13 of~\cite{Triebel86}. After shrinking $\mathcal U$, the last factor is sufficiently small to be absorbed, and hence $\rho=0$, as desired.
\end{proof}

\subsection{The Funk--Radon transform}\label{subsection:funk-radon}
The goal of this section is to establish Proposition~\ref{prop:coercivity-main-theorem} and its consequences. At this point we do not impose any restriction on the dimension. Recall from Section~\ref{section:maintheorem} that for
$n\geq2$ and $0<r<\pi/2$, the Funk--Radon transform is defined by
$$
R_r(f)(\sigma)
=
\int_{\Sigma_\sigma^r}f\,dA_{\mathrm{can}}
-
n\cot r
\int_{\Omega_\sigma^r}f\,dV_{\mathrm{can}}.
$$
We study $R_r$ through its spherical harmonic decomposition and the
asymptotic behavior of its spectral multipliers. We begin by extending
it to a bounded $\mathrm{SO}(n+1)$-equivariant operator on
$L^2(\mathbb S^n)$, so that Schur's lemma applies.

\begin{prop}\label{L2steprop}
	For every fixed $r\in (0,\pi/2)$, there exists a constant $C_r>0$ such that 
	$$
	\| R_r(f) \|_{L^2(\mathbb{S}^n)} \leq C_r \| f \|_{L^2(\mathbb{S}^n)}
	$$
	for every $f\in C^\infty(\mathbb{S}^n)$. In particular, the operator $R_r : C^\infty(\mathbb{S}^n) \to C^\infty(\mathbb{S}^n)$ extends uniquely to a $\mathrm{SO}(n+1)$-equivariant bounded linear operator
	$$
	R_r : L^2(\mathbb{S}^n) \to L^2(\mathbb{S}^n).
	$$
	
\end{prop}
\begin{proof}
	By H\"older's inequality, for every $v\in\mathbb{S}^n$ we have
	$$
	\biggl|\int_{\Sigma_v^r} f \, dA_{\mathrm{can}}\biggr|^2
	\leq \mathrm{area}_{\mathrm{can}}(\Sigma_v^r) \int_{\Sigma_v^r} |f|^2 \, dA_{\mathrm{can}},
	\qquad
	\biggl|\int_{\Omega_v^r} f \, dV_{\mathrm{can}}\biggr|^2
	\leq \mathrm{vol}_{\mathrm{can}}(\Omega_v^r) \int_{\Omega_v^r} |f|^2 \, dV_{\mathrm{can}} .
	$$
	
	Integrate these inequalities over $v\in\mathbb{S}^n$ and use the geometric integral formulas
	(Lemma~\ref{lem:geometric.integral.formula.boundary})
	to obtain
	$$
	\int_{\mathbb{S}^n} \biggl|\int_{\Sigma_v^r} f \, dA_{\mathrm{can}}\biggr|^2 dV_{\mathrm{can}}(v)
	\leq \mathrm{area}_{\mathrm{can}}(\Sigma_v^r)^2 \int_{\mathbb{S}^n} |f|^2 \, dV_{\mathrm{can}},
	$$
	$$
	\int_{\mathbb{S}^n} \biggl|\int_{\Omega_v^r} f \, dV_{\mathrm{can}}\biggr|^2 dV_{\mathrm{can}}(v)
	\leq \mathrm{vol}_{\mathrm{can}}(\Omega_v^r)^2 \int_{\mathbb{S}^n} |f|^2 \, dV_{\mathrm{can}} .
	$$
	
	From the expression $R_r(f)(v) = \int_{\Sigma_v^r} f \, dA_{\mathrm{can}} - n\cot r \int_{\Omega_v^r} f \, dV_{\mathrm{can}}$ we obtain:
	$$
	\|R_r(f)\|_{L^2(\mathbb{S}^n)}
	\leq \Bigl( \mathrm{area}_{\mathrm{can}}(\Sigma_v^r) + n\cot r\, \mathrm{vol}_{\mathrm{can}}(\Omega_v^r) \Bigr) \|f\|_{L^2(\mathbb{S}^n)} .
	$$
	This is the required estimate with $C_r = \mathrm{area}_{\mathrm{can}}(\Sigma_v^r) + n\cot r\, \mathrm{vol}_{\mathrm{can}}(\Omega_v^r)$.
\end{proof}

Consequently, Schur's lemma yields the following diagonalization of the
Funk--Radon transform.

\begin{cor}\label{cor:funk-radon-diagonalization}
	For every fixed $r\in(0,\pi/2)$, there exists a bounded sequence of real
	numbers $\{\mu_k(r)\}_{k\geq0}$ such that
	$$R_r(f)=\mu_k(r)f$$ 
	for every $f\in\mathcal H^k(\mathbb S^n)$. In particular, if
	$f=\sum_{k=0}^{\infty}f_k$ is the spherical harmonic decomposition of
	$f\in L^2(\mathbb S^n)$, then
	$$R_r(f)=\sum_{k=0}^{\infty}\mu_k(r)f_k,$$ 
	with convergence in $L^2(\mathbb S^n)$.
\end{cor}

For each fixed $r\in(0,\pi/2)$, the sequence
$\{\mu_k(r)\}_{k\geq0}$ is called the sequence of \emph{spectral
	multipliers} of $R_r$.

To compute these multipliers, it suffices to evaluate $R_r$ on a
convenient nonzero element of each space
$\mathcal H^k(\mathbb S^n)$. As is usual in the computation of
such eigenvalues, the most convenient choice is the most symmetric
element of each eigenspace. On the
sphere, these are the zonal, or ultraspherical, harmonics, which are
invariant under the stabilizer of a fixed pole
$\sigma\in\mathbb S^n$, and therefore can be understood as functions of the the height function determined by this pole. For each $k\geq0$, the space of degree $k$ spherical harmonics with this invariance is one-dimensional
and is generated by a family of orthogonal polynomials known as the
ultraspherical, or Gegenbauer, polynomials. In what follows we present its definition, for a more detailed discussion see \cite{Szego75} and \cite{Vilenkin68}. 

For
$\lambda>0$, the Gegenbauer polynomials $C_k^\lambda$ are
defined by the generating function
\begin{equation}\label{eq:gegenbauer-generating-function}
	\frac{1}{(1-2xt+t^2)^\lambda}
	=
	\sum_{k=0}^{\infty}C_k^\lambda(x)t^k.
\end{equation}
For spherical harmonics on $\mathbb S^n$, the
relevant parameter is $\lambda=(n-1)/2$. Thus, after fixing a pole $\sigma\in\mathbb S^n$, 
the function $\xi_{k,\sigma}(x) := C_k^{\frac{n-1}{2}}\bigl(\langle x,\sigma\rangle\bigr)$ 
satisfies
$$
\xi_{k,\sigma}\in\mathcal H^k(\mathbb S^n).
$$
Once given the definition of ultraspherical
harmonics the spectral multipliers can be computed as follows. We write $\omega_{n-1}=\mathrm{area}_{\mathrm{can}}(\mathbb S^{n-1})$.

\begin{lem}\label{lem:gegenbauer-representation}
	For every $n\geq2$ and every $r\in(0,\pi/2)$, the spectral multipliers
	of $R_r$ satisfy the following identities:
	\begin{enumerate}
		\item[(a)]
		$
		\mu_0(r)
		=
		\omega_{n-1}\sin^{n-1}r
		-
		n\omega_{n-1}\cot r
		\int_0^r\sin^{n-1}t\,dt.
		$
		
		\item[(b)] $\mu_1(r)=0$.
		
		\item[(c)] For every $k\geq2$,
		\begin{equation}\label{eq:mu-gegenbauer}
			\mu_k(r)
			=
			-\omega_{n-1}(n+1)(n-1)
			\sin^{n+1}r\,\widetilde{\mu}_k(r),
		\end{equation}
		where
		\begin{equation}\label{eq:reduced-multiplier-gegenbauer}
			\widetilde{\mu}_k(r)
			=
			\frac{1}{
				\lambda_k(\mathbb S^n)\binom{n+k-2}{k}
			}
			C_{k-2}^{\frac{n+3}{2}}(\cos r),
		\end{equation}
		and $\lambda_k(\mathbb S^n)=k(k+n-1)$.
	\end{enumerate}
\end{lem}
\begin{proof}
	Item~(a) follows directly from the definition of $R_r$. To prove item~(b), 
	fix a pole $\sigma\in\mathbb S^n$ and consider the linear spherical 
	harmonic $C_1^{\frac{n-1}{2}}(x)=(n-1)x$
	\cite[Chapter~IX, Section~3.3, equation~(2)]{Vilenkin68}. Therefore,
	$$
	\begin{aligned}
		\mu_1(r)
		&=
		\frac{\omega_{n-1}}{C_1^{\frac{n-1}{2}}(1)}
		\left[
		\sin^{n-1}r\;C_1^{\frac{n-1}{2}}(\cos r)
		-
		n\cot r
		\int_0^r
		C_1^{\frac{n-1}{2}}(\cos t)\sin^{n-1}t\,dt
		\right] \\
		&=
		\frac{\omega_{n-1}}{n-1}
		\left[
		(n-1)\sin^{n-1}r\cos r
		-
		n\cot r
		\int_0^r
		(n-1)\cos t\,\sin^{n-1}t\,dt
		\right] \\
		&=
		\frac{\omega_{n-1}}{n-1}
		\left[
		(n-1)\sin^{n-1}r\cos r
		-
		(n-1)\cot r\,\sin^n r
		\right] \\
		&=0.
	\end{aligned}
	$$
	
	To prove item~(c), fix a pole $\sigma\in\mathbb S^n$ and consider the
	zonal harmonic
	$
	\xi_{k,\sigma}(x)
	=
	C_k^{\frac{n-1}{2}}\bigl(\langle x,\sigma\rangle\bigr).
	$
	By Corollary~\ref{cor:funk-radon-diagonalization},
	$$
	\mu_k(r)
	=
	\frac{R_r(\xi_{k,\sigma})(\sigma)}
	{C_k^{\frac{n-1}{2}}(1)}.
	$$
	
	Let $N_\sigma^r$ denote the outward unit normal to $\Sigma_\sigma^r$.
	Since
	$
	-\Delta_{\mathrm{can}}\xi_{k,\sigma}
	=
	\lambda_k(\mathbb S^n)\xi_{k,\sigma},
	$
	the divergence theorem gives
	\begin{align*}
		\int_{\Omega_\sigma^r}\xi_{k,\sigma}\,dV_{\mathrm{can}}
		&=
		-\frac{1}{\lambda_k(\mathbb S^n)}
		\int_{\Sigma_\sigma^r}
		\left\langle
		N_\sigma^r,\nabla\xi_{k,\sigma}
		\right\rangle
		\,dA_{\mathrm{can}} \\
		&=
		\frac{\sin r}{\lambda_k(\mathbb S^n)}
		\frac{d}{dx}C_k^{\frac{n-1}{2}}(\cos r)
		\mathrm{area}_{\mathrm{can}}(\Sigma_\sigma^r) \\
		&=
		\frac{\omega_{n-1}}{\lambda_k(\mathbb S^n)}
		\sin^n r\,
		\frac{d}{dx}C_k^{\frac{n-1}{2}}(\cos r).
	\end{align*}

	On the other hand, since $\langle p,\sigma \rangle = \cos r$ on $\Sigma_\sigma^r$, we have 
	$$
	\int_{\Sigma_\sigma^r}
	\xi_{k,\sigma}\,dA_{\mathrm{can}}
	=
	\omega_{n-1}\sin^{n-1}r\,
	C_k^{\frac{n-1}{2}}(\cos r).
	$$
	It follows from the definition of $R_r$ that
	$$
	\begin{aligned}
		R_r(\xi_{k,\sigma})(\sigma)
		=
		\omega_{n-1}\sin^{n-1}r
		\left[
		C_k^{\frac{n-1}{2}}(\cos r)
		-
		\frac{n\cos r}{\lambda_k(\mathbb S^n)}
		\frac{d}{dx}C_k^{\frac{n-1}{2}}(\cos r)
		\right].
	\end{aligned}
	$$
	
	To finish the proof we recollect some properties of the Gegenbauer polynomials. First, $C_k^\lambda$ 
	satisfies the following differential equation
	\cite[Chapter~IX, Section~3.2, equation~(3)]{Vilenkin68}:
	$$
	(1-x^2)y''(x)
	-
	(2\lambda+1)xy'(x)
	+
	k(k+2\lambda)y(x)
	=
	0.
	$$
	Taking $\lambda=(n-1)/2$ and recalling that
	$\lambda_k(\mathbb S^n)=k(k+n-1)$, we obtain
	$$
	C_k^{\frac{n-1}{2}}(x)
	-
	\frac{nx}{\lambda_k(\mathbb S^n)}
	\frac{d}{dx}C_k^{\frac{n-1}{2}}(x)
	=
	-\frac{1-x^2}{\lambda_k(\mathbb S^n)}
	\frac{d^2}{dx^2}C_k^{\frac{n-1}{2}}(x).
	$$
	
	Also, the family of Gegenbauer polynomials is closed under
	differentiation
	\cite[Chapter~IX, Section~3.2, equation~(4)]{Vilenkin68}. More precisely,
	$$
	\frac{d}{dx}C_k^\lambda(x)
	=
	2\lambda C_{k-1}^{\lambda+1}(x).
	$$
	Applying this identity twice gives
	$$
	\frac{d^2}{dx^2}C_k^{\frac{n-1}{2}}(x)
	=
	(n-1)(n+1)
	C_{k-2}^{\frac{n+3}{2}}(x).
	$$
	Consequently,
	$$
	R_r(\xi_{k,\sigma})(\sigma)
	=
	-\frac{\omega_{n-1}(n-1)(n+1)}
	{\lambda_k(\mathbb S^n)}
	\sin^{n+1}r\,
	C_{k-2}^{\frac{n+3}{2}}(\cos r).
	$$
	
	Finally, the normalization induced by
	equation~\eqref{eq:gegenbauer-generating-function} gives
	\cite[Chapter~IX, Section~3.3, equation~(9)]{Vilenkin68}
	\begin{equation}\label{C_k^n(1)}
		C_k^{\frac{n-1}{2}}(1)
		=
		\binom{n+k-2}{k}.
	\end{equation}
	This proves item~(c) and completes the proof.
\end{proof}

It is well understood that the continuity, injectivity, and
coercivity of the Funk--Radon transform can be described in terms of the nonvanishing and the asymptotic behavior of these constants, see for instance \cite[Section~4]{Strichartz81}.

The asymptotic behavior of orthogonal polynomials
has been extensively studied. In particular, the representation
obtained above allows us to specialize the $L^2$ extension Proposition~\ref{L2steprop} to Sobolev spaces of real order. We refer to
\cite[Section~4]{Strichartz83} for the definition of the real order Sobolev spaces used below.   

\begin{prop}\label{prop:Rr.regularizing}
	For every fixed $r\in(0,\pi/2)$ and every $s\in\mathbb R$, there
	exists a constant $C_{r,n}>0$ such that
	$$
	\|R_r(f)\|_{H^{s+\frac{n-1}{2}}(\mathbb S^n)}
	\leq
	C_{r,n}\|f\|_{H^s(\mathbb S^n)}
	$$
	for every $f\in C^\infty(\mathbb S^n)$. In particular, the operator
	$$
	R_r:C^\infty(\mathbb S^n)\longrightarrow
	C^\infty(\mathbb S^n)
	$$
	extends uniquely to a bounded linear operator
	$$
	R_r:H^s(\mathbb S^n)
	\longrightarrow
	H^{s+\frac{n-1}{2}}(\mathbb S^n).
	$$
\end{prop}
\begin{proof}
	From the Darboux estimate for Gegenbauer polynomials
	\cite[(8.21.14)]{Szego75} and the normalization~\eqref{C_k^n(1)}, we conclude the existence of a constant $C_{r,n} >0$ such that: 
	\begin{equation*}
		\left| \frac{C_k^{\frac{n+3}{2}}(\cos r)} {C_k^{\frac{n+3}{2}}(1)} \right| \leq C_{r,n} k^{-\frac{n+3}{2}}. 
	\end{equation*}
	In particular, from Lemma~\ref{lem:gegenbauer-representation} it hold true that: 
	$$|\mu_k(r)| \leq C_{r,n} k^{-\frac{n-1}{2}}.$$
	Therefore, the desired result follows from classical arguments.
\end{proof}

\begin{rmk}
	The generalized Funk--Radon transform can be related to other integral
	transforms considered in the literature. In particular, Rubin
	\cite{Rubin00} introduced the generalized Minkowski--Funk transform
	$$
	(M_r^\alpha f)(\sigma)
	=
	\frac{c_{n,\alpha}}
	{(\sin^2 r)^{\alpha-1+\frac n2}}
	\int_{\Omega_\sigma^r}
	\bigl(\langle \sigma,x\rangle-\cos r\bigr)^{\alpha-1}
	f(x)\,dV_{\mathrm{can}}(x),
	$$
	where $\mathrm{Re}\alpha>0$, $r\in(0,\pi)$, and
	$
	c_{n,\alpha}
	=
	2^{\alpha-1}\pi^{-\frac n2}
	\frac{\Gamma\left(\alpha+\frac n2\right)}{\Gamma(\alpha)}.
	$
	
	This operator is also $\mathrm{SO}(n+1)$-equivariant and therefore
	diagonalizes in terms of spherical harmonics. Taking $\alpha=2$ and
	denoting its spectral multipliers by $m_k(r)$, we obtain from
	\cite[Lemma~2.1]{Rubin00} the relation
	$$
	\mu_k(r)
	=
	-\frac{\omega_{n-1}}{n(n+2)}
	\sin^{n+1}r\,
	\bigl(\lambda_k(\mathbb S^n)-n\bigr)m_k(r),
	\qquad k\ge2.
	$$
	In particular, the Sobolev continuity established in the preceding
	proposition can also be recovered from \cite[Theorem~2.9]{Rubin00}.
	More importantly for our purposes, the results of Rubin and subsequent
	authors on the zeros of Gegenbauer polynomials will be useful in our
	study of the injectivity and coercivity of $R_r$.
\end{rmk}

To study injectivity and coercivity, we return to the notation
introduced in Section~\ref{section:maintheorem}. By Lemma~\ref{lem:gegenbauer-representation} (see also~\eqref{Ginj} and~\eqref{ZkandZ}),
the exceptional sets are given by
$$
Z_k
=
\left\{
r\in(0,\pi/2):
C_{k-2}^{\frac{n+3}{2}}(\cos r)=0
\right\},
\qquad k\geq2,
$$
and
$$
Z=\bigcup_{k\geq2}Z_k,
\qquad
\mathcal G_{\mathrm{inj}}^n=(0,\pi/2)\setminus Z.
$$
Since each $Z_k$ is finite, we obtain the following.

\begin{lem}\label{lem:generic-injectivity}
	The set $Z$ is countable. Consequently,
	$\mathcal G_{\mathrm{inj}}^n$ is a residual dense subset of
	$(0,\pi/2)$ and has full Lebesgue measure.
\end{lem}

Similarly, for $\alpha\in\mathbb R$, the coercivity estimate defining
$\mathcal G_\alpha^n$ is equivalent to the condition
$$
\inf_{k\geq2}
|\mu_k(r)|
\left(1+\lambda_k(\mathbb S^n)\right)^{\alpha/2}
>0.
$$
Using equation~\eqref{eq:reduced-multiplier-gegenbauer}, we therefore obtain
\begin{equation}\label{Galphaset}
	\mathcal G_\alpha^n
	=
	\left\{
	r\in\mathcal G_{\mathrm{inj}}^n:
	\inf_{k\geq1}
	k^{\alpha-n}
	\left|C_k^{\frac{n+3}{2}}(\cos r)\right|
	>0
	\right\}.
\end{equation}
We say that $R_r$ is \emph{$\alpha$-coercive} when
$r\in\mathcal G_\alpha^n$. The following characterization of coercivity in terms of the spectral
multipliers is classical.

\begin{prop}\label{prop:Rr.coercivity}
	Let $\alpha,s\in\mathbb R$ and suppose that $r\in\mathcal G_\alpha^n$.
	Then there exists a constant $C_{r,n,\alpha}>0$ such that
	$$
	\|f\|_{H^s(\mathbb S^n)}
	\leq
	C_{r,n,\alpha}
	\|R_r(f)\|_{H^{s+\alpha}(\mathbb S^n)}
	$$
	for every
	$
	f\in C^\infty(\mathbb S^n)
	\cap
	\bigl(
	\mathcal H^0(\mathbb S^n)
	\oplus
	\mathcal H^1(\mathbb S^n)
	\bigr)^\perp.
	$ In particular, if $r\in\mathcal G_{\frac{n-1}{2}}^n$, then
	$$
	R_r:
	H^s(\mathbb S^n)
	\cap
	\bigl(
	\mathcal H^0(\mathbb S^n)
	\oplus
	\mathcal H^1(\mathbb S^n)
	\bigr)^\perp
	\longrightarrow
	H^{s+\frac{n-1}{2}}(\mathbb S^n)
	\cap
	\bigl(
	\mathcal H^0(\mathbb S^n)
	\oplus
	\mathcal H^1(\mathbb S^n)
	\bigr)^\perp
	$$
	is an isomorphism of Hilbert spaces.
\end{prop}

\begin{proof}
	Write $f=\sum_{k\geq2}f_k$ from to the spherical harmonic
	decomposition. If
	$$
	c_{r,\alpha}
	=
	\inf_{k\geq2}
	|\mu_k(r)|
	\left(1+\lambda_k(\mathbb S^n)\right)^{\alpha/2},
	$$
	then
	$$
	\begin{aligned}
		\|R_r(f)\|_{H^{s+\alpha}(\mathbb S^n)}^2
		&=
		\sum_{k\geq2}
		\left(1+\lambda_k(\mathbb S^n)\right)^{s+\alpha}
		|\mu_k(r)|^2
		\|f_k\|_{L^2(\mathbb S^n)}^2  \\
		&\geq
		c_{r,\alpha}^2
		\sum_{k\geq2}
		\left(1+\lambda_k(\mathbb S^n)\right)^s
		\|f_k\|_{L^2(\mathbb S^n)}^2  \\
		&=
		c_{r,\alpha}^2
		\|f\|_{H^s(\mathbb S^n)}^2.
	\end{aligned}
	$$
	This proves the coercive estimate. The last assertion of the statement follows from
	Proposition~\ref{prop:Rr.regularizing}.
\end{proof}

Although $\mathcal G_{\mathrm{inj}}^n$ is large both topologically and
from the measure-theoretic point of view, as we will see below, for several choices 
of $\alpha$ and $n$, the set $\mathcal G_\alpha^n$ is empty. One very relevant example of this phenomena is the case $n=3$ and $\alpha=\frac{n-1}{2}=1$, for which the set $\mathcal G_1^3$ is empty. 
Consequently, for every $r\in\mathcal G_{\mathrm{inj}}^3$ and every $s\in\mathbb R$, the injective operator
$$
R_r:
H^s(\mathbb S^3)\cap
\bigl(\mathcal H^0(\mathbb S^3)\oplus\mathcal H^1(\mathbb S^3)\bigr)^\perp
\longrightarrow
H^{s+1}(\mathbb S^3)\cap
\bigl(\mathcal H^0(\mathbb S^3)\oplus\mathcal H^1(\mathbb S^3)\bigr)^\perp
$$
is not an isomorphism. This is in sharp contrast with the minimal case;
compare with \cite[Section~4, Lemma~4.3]{Strichartz81}.

Most strikingly, generic arguments in the sense of Baire are also not available. In every dimension the set
$
\bigcup_{\alpha\in\mathbb R}\mathcal G_\alpha^n
$ is topologically small. This phenomenon is a consequence of the 
distribution of the zeros of the family of Gegenbauer polynomials
$C_k^{\frac{n+3}{2}}$.

\begin{prop}\label{prop:gegenbauer.residual.small.values}
	Let $\lambda>0$ and $\alpha\in\mathbb R$. Then the set
	$$
	\left\{
	r\in\left(0,\frac{\pi}{2}\right):
	\liminf_{k\to\infty}
	k^\alpha
	\left|C_k^\lambda(\cos r)\right|
	=0
	\right\}
	$$
	is residual in $\left(0,\frac{\pi}{2}\right)$. In particular,
	$
	\bigcup_{\alpha\in\mathbb R}\mathcal G_\alpha^n
	$
	is meager.
\end{prop}

\begin{proof}
	By \cite[Theorem~6.11.2, equation~(6.11.12), and the subsequent
	observation]{Szego75}, for every
	$N\geq1$ the set
	$$
	\bigcup_{k\geq N}
	\left\{
	r\in\left(0,\frac{\pi}{2}\right):
	C_k^\lambda(\cos r)=0
	\right\}
	$$
	is dense in $\left(0,\frac{\pi}{2}\right)$. Consequently, for $m,N\geq1$
	$$
	\mathcal U_{m,N}
	=
	\bigcup_{k\geq N}
	\left\{
	r\in\left(0,\frac{\pi}{2}\right):
	k^\alpha
	\left|C_k^\lambda(\cos r)\right|
	<
	\frac1m
	\right\},
	$$
	is open and  dense. Therefore, 
	$$
	\left\{
	r\in\left(0,\frac{\pi}{2}\right):
	\liminf_{k\to\infty}
	k^\alpha
	\left|C_k^\lambda(\cos r)\right|
	=0
	\right\} = \bigcap_{m=1}^{\infty}
	\bigcap_{N=1}^{\infty}
	\mathcal U_{m,N}
	$$
	is residual, as claimed. Consequently the set
	$\mathcal G_\alpha^n$ is meager. In particular, since
	$
	\bigcup_{\alpha\in\mathbb R}\mathcal G_\alpha^n
	=
	\bigcup_{m\in\mathbb N}\mathcal G_m^n
	$ is a countable union of meager sets is therefore meager.
\end{proof}

As previously commented, to the best of the authors' knowledge, an explicit characterization of the radii $r$ for which $R_r$ is injective and $\alpha$-coercive for some $\alpha\in\mathbb R$ is not presently known. For a wider literature in the area see \cite{Rubin00,BilykMastrianniSteinerberger24}.
Nevertheless, several relevant results in this direction are available and will be recalled below. From a perturbative point
of view, Proposition~\ref{prop:gegenbauer.residual.small.values} shows
that coercivity is not robust with respect to the choice of the radius.
This leads us to change perspective and investigate the size of the set
of coercive radii from the measure-theoretic point of view.

Before stating the general results, we discuss the case $n=3$. In this
dimension, the relevant Gegenbauer polynomials admit elementary
trigonometric representations, making explicit computations possible.
This example provides useful intuition for the general results
presented afterward.

\begin{example}\label{example:dimension-three-coercivity}\label{formulaparamu_kemdim3}
	Assume that $n=3$. To understand the decay of the spectral
	multipliers, it is convenient to keep separate the two geometric
	components of the Funk--Radon transform. Define
	$$
	\mathcal T_r(f)(\sigma)
	=
	\int_{\Sigma_\sigma^r}f\,dA_{\mathrm{can}},
	\qquad
	\mathcal B_r(f)(\sigma)
	=
	\int_{\Omega_\sigma^r}f\,dV_{\mathrm{can}}.
	$$
	Then
	$$
	R_r=\mathcal T_r-3\cot r\,\mathcal B_r.
	$$
	
	Both operators are $\mathrm{SO}(4)$-equivariant. Hence, by Schur's
	lemma, there exist spectral multipliers $\tau_k(r)$ and $\beta_k(r)$
	such that
	$$
	\mathcal T_r(f)=\tau_k(r)f,
	\qquad
	\mathcal B_r(f)=\beta_k(r)f,
	\qquad
	f\in\mathcal H^k(\mathbb S^3).
	$$
	In particular,
	$$
	\mu_k(r)=\tau_k(r)-3\cot r\,\beta_k(r).
	$$
	
	By \cite{Vilenkin68},
	$
	C_k^1(\cos t)
	=
	\frac{\sin((k+1)t)}{\sin t}.
	$
	Therefore, setting $m=k+1$, evaluation on the zonal harmonic of
	degree $k$ gives
	$$
	\tau_{m-1}(r)
	=
	\frac{4\pi\sin r}{m}\sin(mr),
	\quad
	\beta_{m-1}(r)
	=
	\frac{4\pi}{m(m^2-1)}
	\left[
	\cos r\,\sin(mr)
	-
	m\sin r\,\cos(mr)
	\right].
	$$
	Consequently,
	\begin{equation}\label{eq:mu-dimension-three}
		\begin{aligned}
			\mu_{m-1}(r)
			= {}&
			\frac{4\pi\sin r}{m}\sin(mr)\\
			&-
			\frac{12\pi\cot r}{m(m^2-1)}
			\left[
			\cos r\,\sin(mr)
			-
			m\sin r\,\cos(mr)
			\right],
			\qquad m\geq3.
		\end{aligned}
	\end{equation}
	
	The first main observation is that these integral operators display two different regularizing orders:
	$$
	|\tau_{m-1}(r)|\lesssim m^{-1},
	\qquad
	|\beta_{m-1}(r)|\lesssim m^{-2}.
	$$
	Thus, the upper regularity estimate for $R_r$ is governed by the
	spherical transform $\mathcal T_r$. However, on the harmonic degrees
	for which $\tau_k(r)=0$, the multiplier of $R_r$ is determined by the
	more regularizing ball transform $\mathcal B_r$. These degrees may
	therefore obstruct coercivity at the regularizing order predicted by
	$\mathcal T_r$.
	
	The radius $r=\pi/3$ already illustrates this phenomenon. A direct
	inspection of the residue classes of $m$ modulo $6$ in
	\eqref{eq:mu-dimension-three} shows that
	$$
	\mu_{m-1}\left(\frac{\pi}{3}\right)\neq0
	\qquad
	\text{for every }m\geq3.
	$$
	Hence $\pi/3\in\mathcal G_{\mathrm{inj}}^3$. On the other hand, if
	$m=3\ell$, then
	$$
	\tau_{3\ell-1}\left(\frac{\pi}{3}\right)=0, \quad
	\mu_{3\ell-1}\left(\frac{\pi}{3}\right)
	=
	\frac{6\pi(-1)^\ell}{9\ell^2-1}.
	$$
	Since
	$
	1+\lambda_{3\ell-1}(\mathbb S^3)=9\ell^2,
	$
	we obtain
	$$
	\left(1+\lambda_{3\ell-1}(\mathbb S^3)\right)^{1/2}
	\left|
	\mu_{3\ell-1}\left(\frac{\pi}{3}\right)
	\right|
	\longrightarrow0.
	$$
	Therefore
	$
	\frac{\pi}{3}\notin\mathcal G_1^3.
	$ Thus, although $R_{\pi/3}$ is injective, it is not an isomorphism
	between the Sobolev spaces $H^s$ and $H^{s+1}$. This loss of coercivity is not particular of $\pi/3$, in fact the same argument prevails whenever $r/\pi\in\mathbb Q$. In fact, suppose that
	$r/\pi=p/q\in\mathbb Q$, with $(p,q)=1$. Taking $m=\ell q$ in \eqref{eq:mu-dimension-three}, and observing that
	$1+\lambda_{m-1}(\mathbb S^3)=m^2$, we obtain
	$$
	\left(1+\lambda_{m-1}(\mathbb S^3)\right)^{\frac{1+\nu}{2}}
	|\mu_{m-1}(r)|
	=
	\frac{12\pi\cos r\,m^{1+\nu}}{m^2-1}
	\longrightarrow0
	$$
	for every $0\leq\nu<1$. Thus no radius $r \in (0,\pi/2)$ for which $r/\pi \in \mathbb{Q}$  is such that $R_r$ is
	$(1+\nu)$-coercive.     
	
	If we insist on obtaining coercivity of order $1+\nu$, for $0\leq\nu<1$, avoiding radius $r \in (0,\pi/2)$ for which $r/\pi \in \mathbb{Q}$ is not enough. Let $r/\pi\in\mathbb R\setminus\mathbb Q$. By
	\eqref{eq:mu-dimension-three},
	$$
	m^{1+\nu}\left|\mu_{m-1}(r)\right|
	\leq
	C_r\left(
	m^\nu|\sin(mr)|+m^{\nu-1}
	\right).
	$$
	Arguing as in \cite[Theorem~3.7]{Rubin00}, Dirichlet's approximation
	theorem provides a sequence $m_j\to\infty$ such that
	$$
	\mathrm{dist}
	\left(
	\frac{m_jr}{\pi},\mathbb Z
	\right)
	\leq
	\frac1{m_j}.
	$$
	Since
	$$
	|\sin(m_jr)|
	\leq
	\pi\mathrm{dist}
	\left(
	\frac{m_jr}{\pi},\mathbb Z
	\right),
	$$
	we obtain, after increasing $C_r$ if necessary,
	$$
	m_j^{1+\nu}\left|\mu_{m_j-1}(r)\right|
	\leq
	C_r m_j^{\nu-1}
	\longrightarrow0.
	$$
	Therefore, $\mathcal G_{1+\nu}^3$ is empty for every
	$0\leq\nu<1$.
	
	A first discrepancy appears at the endpoint $\alpha=2$. If
	$r/\pi\in\mathbb Q$ and $r\in\mathcal G_{\mathrm{inj}}^3$, the
	periodicity of $\sin(mr)$ and $\cos(mr)$ in
	\eqref{eq:mu-dimension-three} gives
	$$
	|\mu_k(r)|\geq c_r k^{-2},
	\qquad k\geq2,
	$$
	and hence $r\in\mathcal G_2^3$. 
	
	As already observed by Rubin, ordinary rational approximation is no
	longer sufficient to treat the irrational case at the endpoint
	$\alpha=2$. Indeed, according to
	\eqref{eq:mu-dimension-three}, the term that must be controlled is
	$$
	m\sin r\,\sin(mr)
	+
	3\cos r\,\cos(mr).
	$$
	After dividing by $\sin r$, we must find a sequence $m_j\to\infty$
	such that
	$$
	m_j\sin(m_jr)
	+
	3\cot r\,\cos(m_jr)
	\longrightarrow0.
	$$
	Thus, controlling the values of $m\sin(mr)$ alone is not sufficient:
	the two oscillatory terms must be controlled simultaneously. It is
	enough, for instance, to find a sequence satisfying
	$$
	m_j\sin(m_jr)\longrightarrow3\cot r,
	\qquad
	\cos(m_jr)\longrightarrow-1.
	$$
	
	This problem was recently resolved for almost every radius by Han--Rahimi~\cite{HanRahimi26}. In fact, we can deduce
	from \cite[Lemma~3.1]{HanRahimi26} (see also Theorem~2.6)
	that, for almost every irrational $r/\pi$, there exists a sequence
	$m_j\to\infty$ satisfying the simultaneous limits above.. Hence
	$$
	m_j^2|\mu_{m_j-1}(r)|
	\longrightarrow0.
	$$
	Therefore, almost every irrational radius fails to be $2$-coercive.
	Since the rational radii form a countable set, we conclude that $\mathcal G_2^3$
	has null measure. 
	
	Notice that, for $\alpha<(n+1)/2=2$, the set
	$\mathcal G_\alpha^3$ has measure zero. For $\alpha>(n+1)/2$,
	however, we will see that this set has full measure. This sharp
	change in its measure-theoretic behavior can be geometrically
	explained by the distribution and spacing of the zeros of the
	Gegenbauer polynomials.
	
	Indeed, fix $0<\delta<\pi/4$ and set
	$I_\delta=[\delta,\pi/2-\delta]$. From
	\cite[Lemma~14]{BrandoliniGariboldiGigante25}
	(see also Appendix~\ref{appendix:gegenbauer-sublevels}), there exist
	constants $c_{\delta,n}>0$ and $k_0(\delta,n)\geq2$ such that
	$$
	|\mu_k(r)|
	\geq
	c_{\delta,n}
	k^{-\frac{n-1}{2}}
	\min
	\left\{
	1,
	k\Theta_k(r)
	\right\},
	\qquad
	r\in I_\delta,
	\quad
	k\geq k_0(\delta,n),
	$$
	where
	$$
	\Theta_k(r)
	=
	\mathrm{dist}(r,Z_k).
	$$
	
	Given $\eps>0$, define
	$$
	E_k(\delta,\eps)
	=
	\left\{
	r\in I_\delta:
	\Theta_k(r)<k^{-2-\eps}
	\right\}.
	$$
	Since $Z_k$ contains at most $k$ roots, the set
	$E_k(\delta,\eps)$ is contained in the union of at most $k$
	intervals of radius $k^{-2-\eps}$. Consequently,
	$$
	|E_k(\delta,\eps)|
	\leq
	C_n k\,k^{-2-\eps}
	=
	C_n k^{-1-\eps}.
	$$
	It follows that
	$$
	\sum_{k=2}^{\infty}|E_k(\delta,\eps)|<\infty.
	$$
	
	By the Borel--Cantelli lemma, for almost every $r\in I_\delta$ there
	exists $K_{r,\delta,\eps}\geq2$ such that
	$$
	\Theta_k(r)
	\geq
	k^{-2-\eps}
	\qquad
	\text{for every }k\geq K_{r,\delta,\eps}.
	$$
	For these indices, the preceding estimate gives
	$$
	|\mu_k(r)|
	\geq
	c_{\delta,n}
	k^{-\frac{n-1}{2}}
	\min
	\left\{
	1,
	k^{-1-\eps}
	\right\}
	=
	c_{\delta,n}
	k^{-\frac{n+1}{2}-\eps}.
	$$
	
	Excluding the countable set $\bigcup_{k\geq2}Z_k$ and absorbing the
	finitely many remaining indices into the constant, we conclude that,
	for almost every $r\in I_\delta$, there exists
	$c_{r,\delta,n,\eps}>0$ such that
	$$
	|\mu_k(r)|
	\geq
	c_{r,\delta,n,\eps}
	k^{-\frac{n+1}{2}-\eps}
	\qquad
	\text{for every }k\geq2.
	$$
	Finally, exhausting $(0,\pi/2)$ by the intervals $I_{\pi/(8j)}$, we
	obtain a set of full measure in $(0,\pi/2)$ on which
	$$
	|\mu_k(r)|
	\geq
	c_{r,n,\eps}
	k^{-\frac{n+1}{2}-\eps}
	\qquad
	\text{for every }k\geq2.
	$$
	In dimension $n=3$, this proves that $\mathcal G_{2+\eps}^3$ has full
	Lebesgue measure for every $\eps>0$.
\end{example}

For arbitrary dimensions, the change in coercivity behavior at the
critical exponents $(n-1)/2$ and $(n+1)/2$ persists. We emphasize,
however, that in the intermediate regime
$(n-1)/2\leq\alpha<(n+1)/2$, the set $\mathcal G_\alpha^n$ depends on
the dimension and, differently from the case $n=3$, does not need to be empty.
More precisely, we have the following result.

\begin{prop}\label{prop:coercivity.regimes}
	Let $n\geq2$.
	
	\begin{enumerate}[label=(\alph*),ref=\alph*]
		\item\label{item:coercivity.below-natural}
		If $\alpha<(n-1)/2$, then $\mathcal G_\alpha^n=\varnothing$.
		
		\item\label{item:coercivity.subcritical}
		Suppose that $(n-1)/2\leq\alpha<(n+1)/2$.
		\begin{enumerate}[label=(\roman*),ref=\theenumi.\roman*]
			\item\label{item:coercivity.subcritical-irrational}
			No irrational radius belongs to $\mathcal G_\alpha^n$.
			
			\item\label{item:coercivity.subcritical-rational}
			If $r/\pi=p/q\in\mathbb Q$, with $(p,q)=1$, then
			$$
			r\in\mathcal G_\alpha^n
			\quad\Longleftrightarrow\quad
			r\in\mathcal G_{\mathrm{inj}}^n
			\quad\text{and}\quad
			\frac{q(n+1)-2p(n-1)}4\notin\mathbb Z.
			$$
			In particular, if $n\equiv3\pmod4$, then
			$\mathcal G_\alpha^n=\varnothing$.
		\end{enumerate}
		
		\item\label{item:coercivity.critical}
		Suppose that $\alpha=(n+1)/2$.
		\begin{enumerate}[label=(\roman*),ref=\theenumi.\roman*]
			\item\label{item:coercivity.critical-rational}
			If $r/\pi\in\mathbb Q$, then
			$$
			r\in\mathcal G_{\frac{n+1}{2}}^n
			\quad\Longleftrightarrow\quad
			r\in\mathcal G_{\mathrm{inj}}^n.
			$$
			
			\item\label{item:coercivity.critical-irrational}
			The set of irrational radii belonging to
			$\mathcal G_{\frac{n+1}{2}}^n$ has Lebesgue measure zero.
			Consequently,
			$$
			\left|\mathcal G_{\frac{n+1}{2}}^n\right|=0.
			$$
		\end{enumerate}
		
		\item\label{item:coercivity.supercritical}
		If $\alpha>(n+1)/2$, then $\mathcal G_\alpha^n$ has full Lebesgue measure in
		$(0,\pi/2)$.
	\end{enumerate}
\end{prop}
\begin{proof}
	
	Item~\ref{item:coercivity.below-natural} follows directly from
	the Darboux estimate for Gegenbauer polynomials
	\cite[(8.21.14)]{Szego75}. Items
	\ref{item:coercivity.subcritical-irrational} and
	\ref{item:coercivity.subcritical-rational} are respectively
	\cite[Theorems~3.7 and~3.2]{Rubin00}.
	
	Item~\ref{item:coercivity.critical-rational} also follows from
	\cite[Theorem~3.2]{Rubin00}, while
	item~\ref{item:coercivity.critical-irrational} follows from
	\cite[Theorem~2.6]{HanRahimi26}. Finally, item~\ref{item:coercivity.supercritical} follows exactly as
	in Example~\ref{example:dimension-three-coercivity} following \cite[Lemma~14]{BrandoliniGariboldiGigante25}.
\end{proof}

\subsection{The non-linear bound}\label{subsection:nonlinearbound}

In what follows, we establish Proposition~\ref{estimativaprincipal}. To this end, we recall that, for $n\geq 2$ and any integer $k\geq2$, Proposition~\ref{prop:whiteperturbation} guarantees the existence of a neighborhood $\mathscr V_r$ of zero in the smooth topology of $C^\infty(\mathbb S^n)$ and a map
$$
\psi:\mathscr V_r\longrightarrow C^k(\mathcal I_r),
$$
continuous in the $C^k$ topology such that $\psi(0)=0$.  With this notation in place, we prove the following.

\begin{prop}\label{prop:nonlinearalln}
	Let $n\geq2$, $r\in(0,\pi/2]$, and $k\geq3$, and suppose that $\mathscr V_r$ and $\psi$ are as above. After shrinking $\mathscr V_r$ if necessary, there exists a constant $C_r>0$ such that, for every $\rho\in\mathscr V_r$ for which
	$
	\{\Sigma_\sigma^r(\psi(\rho)_\sigma)\}_{\sigma\in\mathbb S^n}
	$
	is  a $h(r)$-CMC Zoll family in $(\mathbb S^n,e^{2\rho}\mathrm{can})$, we have
	$$
	\bigl\|R_r(\rho) + c_r(\rho)\bigr\|_{H^1(\mathbb S^n)}
	\leq
	C_r\|\rho\|_{W^{2,4}(\mathbb S^n)}^2,
	$$
	where
	$
	c_r(\rho)=\frac{1}{n-1}\bigl(\mathcal A(0,0)-
	\mathcal A(\rho,\psi(\rho))\bigr).
	$
\end{prop}

Under the regime of Proposition~\ref{prop:localperturbation} where we have local uniqueness for CMC-Zoll families we obtain the following consequence of the previous proposition.

\begin{thm}\label{thm:nonlinear-estimate}
	Assume that either $n=2$ or $3$ and $0<r<\pi/2$,
	or $n\geq3$ and $r=\pi/2$.
	There exist a neighborhood $\mathscr V_r$ of zero in the
	smooth topology of $C^\infty(\mathbb S^n)$ and a constant
	$C_r>0$ such that every $\rho\in\mathscr V_r$ satisfying
	\begin{itemize}
		\item $e^{2\rho}\mathrm{can}\in \mathcal Z_{h(r)}$, if $0<r<\pi/2$, or; 
		\item $e^{2\rho}\mathrm{can}\in \mathcal Z'$, if $r=\pi/2$,
	\end{itemize}
	it holds true that
	\begin{equation}\label{R+ceq}
		\bigl\|R_r(\rho) + c_r(\rho)\bigr\|_{H^1(\mathbb S^n)}
		\leq
		C_r\|\rho\|_{W^{2,4}(\mathbb S^n)}^2,
	\end{equation}
	where
	$
	c_r(\rho)=\frac{1}{n-1}\bigl(\mathcal A(0,0)-
	\mathcal A(\rho,\psi(\rho))\bigr).
	$
\end{thm}
\begin{rmk}\label{rmk:nonlinear-estimate}
	If $e^{2\rho} \mathrm{can} \in \mathcal {Z}_{h(r)}$ the function $c_r(\rho)$ is constant, therefore from the integral formulas in Appendix~\ref{app:geometric.integral.formula} we obtain the identity
	$$R_r(\mathring{\rho})
	=
	R_r(\rho)+c_r(\rho)
	-
	\dashint_{\mathbb S^n}
	\bigl(R_r(\rho)+c_r(\rho)\bigr)\,
	dV_{\mathrm{can}}.$$
	As a consequence, equation~\ref{R+ceq} implies the weaker bound, 
	$$\bigl\|R_r(\mathring\rho)\bigr\|_{H^1(\mathbb S^n)}
	\leq
	C_r\|\rho\|_{W^{2,4}(\mathbb S^n)}^2,$$
	where $\mathring \rho = \rho - \dashint_{\mathbb{S}^n} \rho \, dV_{\mathrm{can}}$. 
\end{rmk}

Let $k\geq 3$, under the dimensional hypothesis of Theorem~\ref{thm:nonlinear-estimate}, we observe that Proposition~\ref{prop:localperturbation} provides
a neighborhood $\mathscr V_r$ of zero in the smooth topology of
$C^\infty(\mathbb S^n)$, and a map
$\psi:\mathscr V_r\to C^k(\mathcal I_r)$, continuous in the $C^k$
topology with $\psi(0)=0$, such that every Zoll family
under consideration, with the area normalization in
the minimal case, can be written, after reparametrization,
as
$\{\Sigma_\sigma^r(\psi(\rho)_\sigma)\}_{\sigma\in\mathbb S^n}$.
Moreover, if $r=\pi/2$, the following symmetries hold true $$\psi(\rho)(x,-\sigma)=-\psi(\rho)(x,\sigma), \qquad \Sigma^{\pi/2}_{-\sigma}(\psi(\rho))=\Sigma^{\pi/2}_{\sigma}(\psi(\rho)).$$
Therefore, the previous theorem reduces to a special case of Proposition~\ref{prop:nonlinearalln} as claimed.

The remainder of this section is devoted to the proof of
Proposition~\ref{prop:nonlinearalln} for arbitrary $n\geq2$
and $r\in(0,\pi/2]$. We fix an integer $k\geq3$ and organize
the argument into a sequence of lemmas.

Set $\beta(\rho)=\rho\circ\Sigma(\psi(\rho))$ and, for $m=1,2$, write
$\nabla^m\beta(\rho)=(\nabla^m\rho)\circ\Sigma(\psi(\rho))$, where
derivatives are taken with respect to the canonical metric. We write
$A\lesssim B$ if $A\leq CB$, where $C$ depends only on $n$, $r$, and
the fixed neighborhood, but not on $\rho$ or $\sigma$. Also will be convenient to define the following operator, for 
$\phi\in C^1(\mathcal I_r)$ we set
$
\partial_u\phi(x,\sigma)
=
d\phi_{(x,\sigma)}
\bigl(
\cos(r)u-\langle x,u\rangle\sigma,
u
\bigr).
$
In particular, identifying $\rho\in C^1(\mathbb S^n)$ with its pullback
to $\mathcal I_r$ by the first projection, we have
$
\partial_u\rho(x,\sigma)
=
\left\langle
\nabla\rho(x),
\cos(r)u-\langle x,u\rangle\sigma
\right\rangle.
$

We begin with the following Taylor expansion of the weighted area functional
$\mathcal A\bigl(\rho,\psi(\rho)\bigr)$ around the canonical family.

\begin{lem}\label{taylorexpforA}
	Let $n\geq 2$ and $r\in(0,\pi/2]$. After possibly shrinking
	$\mathscr V_r$, for every $\rho\in\mathscr V_r$ there exists a smooth
	map $A_\rho:\mathcal I_r\times\mathbb R^3\to\mathbb R$ such that, for
	every $\sigma\in\mathbb S^n$,
	\begin{align}
		\mathcal A\bigl(\rho,\psi(\rho)\bigr)(\sigma)
		-{}&
		\mathcal A(0,0)(\sigma)
		=
		(n-1)
		\int_{\Sigma_\sigma^r}
		\beta(\rho)\,
		dA_{\mathrm{can}}
		-
		h(r)
		\int_{\Omega_\sigma^r}
		\bigl(e^{n\rho}-1\bigr)\,
		dV_{\mathrm{can}}
		\nonumber\\
		&+
		h(r)
		\int_{\Sigma_\sigma^r}
		\bigl(1-e^{n\rho}\bigr)
		\psi(\rho)\,
		dA_{\mathrm{can}}
		+
		\int_{\Sigma_\sigma^r}
		A_\rho\left(
		\beta(\rho),
		\psi(\rho),
		\left|
		\nabla^{\Sigma_\sigma^r}\psi(\rho)
		\right|^2
		\right)
		dA_{\mathrm{can}}.
		\label{1oraproofA}
	\end{align}
	Moreover, the map $A_\rho$ satisfies the following estimates:
	$$
	|A_\rho(p,z,u,y)|
	+
	|D_pA_\rho(p,z,u,y)|
	\lesssim
	z^2+u^2+|y|,
	$$
	$$
	|\partial_zA_\rho(p,z,u,y)|
	+
	|\partial_uA_\rho(p,z,u,y)|
	\lesssim
	|z|+|u|+|y|,
	$$
	and
	$$
	|\partial_yA_\rho(p,z,u,y)|
	\lesssim 1,
	$$
	for every $p\in\mathcal I_r$ and every $(z,u,y)$ sufficiently close
	to $(0,0,0)$. The implicit constants are uniform for
	$\rho\in\mathscr V_r$.
\end{lem}
\begin{proof}
	For $p=(x,\sigma)\in\mathcal I_r$ and $(z,u,y)$ near the origin,
	define
	$$
	\begin{aligned}
		B_\rho(p,z,u,y)
		=e^{(n-1)z}j_r(u)^{n-2}
		\sqrt{j_r(u)^2+y} -
		h(r)\int_{-r}^{u}
		e^{n\rho(E_r(p,s))}
		j_r(s)^{n-1}\,ds.
	\end{aligned}
	$$
	The explicit formula for the weighted area functional gives
	$$
	\mathcal A\bigl(\rho,\psi(\rho)\bigr)(\sigma)
	=
	\int_{\Sigma_\sigma^r}
	B_\rho\left(
	(x,\sigma),
	\beta(\rho)(x,\sigma),
	\psi(\rho)(x,\sigma),
	\left|
	\nabla^{\Sigma_\sigma^r}\psi(\rho)(x,\sigma)
	\right|^2
	\right)
	dA_{\mathrm{can}}(x).
	$$
	At $(z,u,y)=(0,0,0)$, direct differentiation gives
	$$
	\partial_zB_\rho(p,0,0,0)=n-1,
	\qquad
	\partial_yB_\rho(p,0,0,0)=\frac12,
	$$
	and
	$$
	\begin{aligned}
		\partial_uB_\rho(p,0,0,0)
		&=
		(n-1)\cot r-h(r)e^{n\rho(x)}
		\\
		&=
		h(r)\bigl(1-e^{n\rho(x)}\bigr).
	\end{aligned}
	$$
	We therefore set
	$$
	\begin{aligned}
		A_\rho(p,z,u,y)
		=
		B_\rho(p,z,u,y)
		-
		B_\rho(p,0,0,0)
		-
		(n-1)z
		-
		h(r)\bigl(1-e^{n\rho(x)}\bigr)u.
	\end{aligned}
	$$
	The coarea formula yields
	$$
	\int_{\Sigma_\sigma^r}
	B_\rho\bigl((x,\sigma),0,0,0\bigr)\,
	dA_{\mathrm{can}}(x)
	=
	\mathcal A(0,0)(\sigma)
	-
	h(r)
	\int_{\Omega_\sigma^r}
	\bigl(e^{n\rho}-1\bigr)\,
	dV_{\mathrm{can}}.
	$$
	Substituting the definition of $A_\rho$ into the expression above for
	the weighted area functional gives the desired expansion. The stated
	estimates follow by applying Taylor's theorem to $A_\rho$ at the
	origin and using the uniform bounds for the derivatives of $B_\rho$.
\end{proof}

From the Taylor expansion in
Lemma~\ref{taylorexpforA}, we obtain
$$
\begin{aligned}
	(n-1)\bigl(R_r(\rho)(\sigma)+c_r(\rho)\bigr)
	={}&
	(n-1)
	\int_{\Sigma_\sigma^r}
	\bigl(\rho(x)-\beta(\rho)(x,\sigma)\bigr)\,
	dA_{\mathrm{can}}(x)
	\\
	&+
	h(r)
	\int_{\Omega_\sigma^r}
	\bigl(e^{n\rho}-1-n\rho\bigr)\,
	dV_{\mathrm{can}}
	\\
	&-
	h(r)
	\int_{\Sigma_\sigma^r}
	\bigl(1-e^{n\rho(x)}\bigr)
	\psi(\rho)_\sigma(x)\,
	dA_{\mathrm{can}}(x)
	\\
	&-
	\int_{\Sigma_\sigma^r}
	A_\rho
	\left(
	(x,\sigma),
	\beta(\rho)(x,\sigma),
	\psi(\rho)(x),
	\left|
	\nabla^{\Sigma_\sigma^r}\psi(\rho)(x)
	\right|^2
	\right)
	dA_{\mathrm{can}}(x).
\end{aligned}
$$
The Taylor expansion of the exponential function gives
$$|1-e^{n\rho}|\lesssim|\rho|,\qquad \mbox{and} \qquad
|e^{n\rho}-1-n\rho|\lesssim|\rho|^2.$$
We also observe that
$$|\beta(\rho)(x,\sigma)-\rho(x)|
\lesssim|\nabla\rho(x)|^2+
|\psi(\rho)_\sigma(x)|^2.$$
Therefore, after shrinking $\mathscr V_r$ if necessary, the estimates
for $A_\rho$ and Young's inequality give the following pointwise
quadratic bound:
\begin{equation}\label{eq:pointwise-quadratic-Rr}
	\begin{aligned}
		\left|R_r(\rho)(\sigma)+c_r(\rho)\right|
		\lesssim
		\|\rho\|_{L^2(\Sigma_\sigma^r)}^2
		+
		\|\nabla\rho\|_{L^2(\Sigma_\sigma^r)}^2
		+
		\|\beta(\rho)\|_{L^2(\Sigma_\sigma^r)}^2
		+
		\|\psi(\rho)\|_{H^1(\Sigma_\sigma^r)}^2
		+
		\|\rho\|_{L^2(\Omega_\sigma^r)}^2.
	\end{aligned}
\end{equation}

A similar bound also hold for the derivative of the function $R_r(\rho) + c_r(\rho)$. In fact, for $\sigma\in\mathbb S^n$ and unitary $u\in T_\sigma\mathbb S^n$, we have
$$
\begin{aligned}
	\left\langle\nabla R_r(\rho)(\sigma),u\right\rangle
	={}&
	\int_{\Sigma_\sigma^r}
	\left\langle
	\nabla\rho(x),
	\cos(r)u-\langle x,u\rangle\sigma
	\right\rangle
	dA_{\mathrm{can}}(x)
	\\
	&-
	\frac{n}{n-1}h(r)
	\int_{\Omega_\sigma^r}
	\left\langle
	\nabla\rho(x),
	\langle x,\sigma\rangle u-\langle x,u\rangle\sigma
	\right\rangle
	dV_{\mathrm{can}}(x).
\end{aligned}
$$


Once again, the Zoll condition implies that
$\mathcal A\bigl(\rho,\psi(\rho)\bigr)$ is constant on $\mathbb S^n$.
Therefore, differentiating equation~\eqref{1oraproofA} in the direction
$u\in T_\sigma\mathbb S^n$ and using the estimates for $A_\rho$, we
obtain
\begin{multline}\label{derivadafunk13cmc}
	\left|
	\partial_u
	\left[
	(n-1)
	\int_{\Sigma_\sigma^r}
	\beta(\rho)\,dA_{\mathrm{can}}
	-
	h(r)
	\int_{\Omega_\sigma^r}
	\left(e^{n\rho}-1\right)\,dV_{\mathrm{can}}
	\right]
	\right|
	\\
	\lesssim
	\|\nabla\rho\|_{L^2(\Sigma_\sigma^r)}^2
	+
	\|\beta(\rho)\|_{L^2(\Sigma_\sigma^r)}^2
	+
	\|\partial_u\beta(\rho)\|_{L^2(\Sigma_\sigma^r)}^2
	+
	\|\partial_u\psi(\rho)\|_{H^1(\Sigma_\sigma^r)}^2
	+
	\|\psi(\rho)\|_{H^1(\Sigma_\sigma^r)}^2.
\end{multline}
Comparing the expression inside brackets with the definition of $R_r$
gives
\begin{equation}\label{eq:gradient-Rr-comparison}
	\begin{split}
		\left|
		\left\langle
		\nabla R_r(\rho)(\sigma),u
		\right\rangle
		\right|
		\lesssim{}&
		\left|
		\int_{\Sigma_\sigma^r}
		\bigl(
		\partial_u\rho-\partial_u\beta(\rho)
		\bigr)\,
		dA_{\mathrm{can}}
		\right|
		+
		\left|
		\partial_u
		\int_{\Omega_\sigma^r}
		\left(e^{n\rho}-1-n\rho\right)\,
		dV_{\mathrm{can}}
		\right|
		\\
		&+
		\left|
		\partial_u
		\left[
		(n-1)
		\int_{\Sigma_\sigma^r}
		\beta(\rho)\,dA_{\mathrm{can}}
		-
		h(r)
		\int_{\Omega_\sigma^r}
		\left(e^{n\rho}-1\right)\,dV_{\mathrm{can}}
		\right]
		\right|.
	\end{split}
\end{equation}

The second term in equation~\eqref{eq:gradient-Rr-comparison} is directly
estimated by
\begin{equation}\label{volume.nonlinear.estimate}
	\left|
	\partial_u
	\int_{\Omega_\sigma^r}
	\left(e^{n\rho}-1-n\rho\right)\,
	dV_{\mathrm{can}}
	\right|
	\lesssim
	\|\rho\|_{L^2(\Omega_\sigma^r)}^2
	+
	\|\nabla\rho\|_{L^2(\Omega_\sigma^r)}^2.
\end{equation}
On the other hand, the first term can be estimated as follows. For
$t\in[0,1]$, define
$
f_t(x,\sigma)
=
\partial_u
\left(
\rho\circ
\Sigma_\sigma^r\bigl(t\psi(\rho)\bigr)
\right)(x).
$ Then $f_0=\partial_u\rho$ and
$f_1=\partial_u\beta(\rho)$. Hence,
\begin{equation}\label{boundary.comparison}
	\left|
	\left\langle
	\nabla
	\left(
	\int_{\Sigma_\sigma^r}
	\rho\,dA_{\mathrm{can}}
	\right),
	u
	\right\rangle
	-
	\int_{\Sigma_\sigma^r}
	\partial_u\beta(\rho)\,
	dA_{\mathrm{can}}
	\right|
	\leq
	\int_{\Sigma_\sigma^r}
	|f_0-f_1|\,
	dA_{\mathrm{can}}.
\end{equation}

Direct differentiation gives
$$
\left|
\left.\frac{d}{dt}\right|_{t=0}f_t
\right|
\lesssim
\bigl(|\nabla\rho|+|\nabla^2\rho|\bigr)|\psi(\rho)|
+
|\nabla\rho|\,|\partial_u\psi(\rho)|.
$$
Moreover, after shrinking $\mathscr V_r$ if necessary,
$$
\sup_{t\in[0,1]}
\left|
\frac{d^2}{dt^2}f_t
\right|
\lesssim
|\psi(\rho)|^2
+
|\partial_u\psi(\rho)|^2.
$$
Combining this estimate with equation~\eqref{boundary.comparison}, we obtain the quadratic bound
\begin{equation}\label{boundary.difference.estimate}
	\begin{aligned}
		\left|
		\int_{\Sigma_\sigma^r}
		\bigl(
		\partial_u\rho-\partial_u\beta(\rho)
		\bigr)\,
		dA_{\mathrm{can}}
		\right|
		\lesssim{}&
		\|\nabla\rho\|_{L^2(\Sigma_\sigma^r)}^2
		+
		\|\nabla^2\rho\|_{L^2(\Sigma_\sigma^r)}^2
		\\
		&+
		\|\psi(\rho)\|_{L^2(\Sigma_\sigma^r)}^2
		+
		\|\partial_u\psi(\rho)\|_{L^2(\Sigma_\sigma^r)}^2.
	\end{aligned}
\end{equation}

Combining equation~\eqref{eq:pointwise-quadratic-Rr} with
equations~\eqref{derivadafunk13cmc},
\eqref{eq:gradient-Rr-comparison},
\eqref{volume.nonlinear.estimate}, and
\eqref{boundary.difference.estimate}, we obtain the following lemma.

\begin{lem}\label{lem:pointwise-Rr-estimates}
	After shrinking $\mathscr V_r$ if necessary, for every
	$\rho\in\mathscr V_r$ satisfying
	$\mathcal H\bigl(\rho,\psi(\rho)\bigr)=0$,
	every $\sigma\in\mathbb S^n$, and every unit vector
	$u\in T_\sigma\mathbb S^n$, the following estimates hold:
	\begin{enumerate}
		\item
		$$
		\begin{aligned}
			\left|R_r(\rho)(\sigma)+c_r(\rho)\right|
			\lesssim{}&
			\|\rho\|_{L^2(\Sigma_\sigma^r)}^2
			+
			\|\nabla\rho\|_{L^2(\Sigma_\sigma^r)}^2
			+
			\|\beta(\rho)\|_{L^2(\Sigma_\sigma^r)}^2
			\\
			&+
			\|\psi(\rho)\|_{H^1(\Sigma_\sigma^r)}^2
			+
			\|\rho\|_{L^2(\Omega_\sigma^r)}^2.
		\end{aligned}
		$$
		
		\item
		$$
		\begin{aligned}
			\left|\partial_uR_r(\rho)(\sigma)\right|
			\lesssim{}&
			\|\nabla\rho\|_{L^2(\Sigma_\sigma^r)}^2
			+
			\|\nabla^2\rho\|_{L^2(\Sigma_\sigma^r)}^2
			+
			\|\beta(\rho)\|_{L^2(\Sigma_\sigma^r)}^2
			+
			\|\nabla\beta(\rho)\|_{L^2(\Sigma_\sigma^r)}^2
			\\
			&+
			\|\psi(\rho)\|_{H^1(\Sigma_\sigma^r)}^2
			+
			\|\partial_u\psi(\rho)\|_{H^1(\Sigma_\sigma^r)}^2
			+
			\|\rho\|_{L^2(\Omega_\sigma^r)}^2
			+
			\|\nabla\rho\|_{L^2(\Omega_\sigma^r)}^2.
		\end{aligned}
		$$
	\end{enumerate}
\end{lem}

We can therefore reduce the problem to estimate the terms involving
$\beta(\rho)$ and $\psi(\rho)$. We estimate the former using perturbed
versions of the geometric integral formulas established in
Appendix~\ref{app:geometric.integral.formula}. On the other hand, the terms that involve $\psi(\rho)$ are controlled by the relation
between $\rho$ and $\psi(\rho)$ given by the CMC equation
$\mathcal H(\rho,\psi(\rho))=0$. We now provide these estimates, and for this objective we will need the following result.

\begin{lem}\label{kernJestmts.cmc}
	There exist $\alpha>0$ and smooth functions $V_i: \mathcal I_r\times\mathbb R\times\mathbb R \to\mathbb R$, for $	i=1,\ldots,n+1$, such that, after shrinking $\mathscr V_r$ if necessary, for every
	$\rho\in\mathscr V_r$ satisfying
	$e^{2\rho}\mathrm{can}\in\mathcal Z_{h(r)}$,
	$\sigma\in\mathbb S^n$, and $i=1,\ldots,n+1$,
	\begin{equation}\label{orthogonal.condition.cmc}
		\int_{\Sigma_\sigma^r}
		\psi(\rho)
		\bigl(x_i-\cos(r)\sigma_i\bigr)\,
		dA_{\mathrm{can}}(x)
		=
		\int_{\Sigma_\sigma^r}
		V_i\Bigl(
		(x,\sigma),
		\psi(\rho),
		\left|
		\nabla^{\Sigma_\sigma^r}
		\psi(\rho)
		\right|^2
		\Bigr)\,
		dA_{\mathrm{can}}(x).
	\end{equation}
	Moreover, for every $\omega\in\mathcal I_r$, $s,t\in\mathbb R$, and
	$i=1,\ldots,n+1$,
	\begin{equation}\label{constant.quadratic.cmc}
		|V_i(\omega,s,t)|
		+
		|D_\omega V_i(\omega,s,t)|
		\leq
		\alpha\bigl(|s|^2+|t|\bigr),
	\end{equation}
	\begin{equation}\label{constant.quadratic2.cmc}
		|\partial_sV_i(\omega,s,t)|
		\leq
		\alpha\bigl(|s|+|t|\bigr),
		\qquad
		|\partial_tV_i(\omega,s,t)|
		\leq\alpha.
	\end{equation}
\end{lem}
\begin{proof}
	We start with the case $r \in (0,\pi/2)$. Let $L_r=\{\Sigma_\sigma^r:\sigma\in\mathbb S^n\}$, and consider also the following map defined on the space of embedded hypersurfaces
	$$
	\Omega(S)=\int_S x\,dA_{\mathrm{can}}.
	$$
	We recall that its restriction to $L_r$ is an embedding for every
	$r\in(0,\pi/2)$ and that its image $\Omega(L_r)$ is a round sphere.
	Let $P$ be the tubular projection onto this sphere. Without loss of
	generality, we may assume that the graphical perturbation map $\psi$
	constructed above satisfies, as in
	\cite[Theorem~3.2]{BWhite91},
	\begin{equation}\label{Omegagauge}
		P\left(
		\Omega\bigl(\Sigma_\sigma^r(\psi(\rho))\bigr)
		\right)
		=
		\Omega(\Sigma_\sigma^r).
	\end{equation}
	
	Consequently,
	$\Omega(\Sigma_\sigma^r(\psi(\rho)))-\Omega(\Sigma_\sigma^r)$
	is orthogonal to
	$T_{\Omega(\Sigma_\sigma^r)}\Omega(L_r)=\sigma^\perp$.
	Moreover, for $w\in\sigma^\perp$ and a graphical perturbation
	$\varphi$, the first variation of $\Omega$ is
	$$
	\begin{aligned}
		\left\langle
		\left.\frac{d}{dt}\right|_{t=0}
		\Omega\bigl(\Sigma_\sigma^r(t\varphi)\bigr),
		w
		\right\rangle
		&=
		\int_{\Sigma_\sigma^r}
		\varphi
		\left\langle
		N_\sigma^r+h(r)x,w
		\right\rangle
		\,dA_{\mathrm{can}}
		\\
		&=
		n\cot r
		\int_{\Sigma_\sigma^r}
		\varphi
		\left\langle
		x-\cos(r)\sigma,w
		\right\rangle
		\,dA_{\mathrm{can}}.
	\end{aligned}
	$$
	
	Notice that if $w_i=e_i-\sigma_i\sigma$, we have smooth functions $\widetilde V_i$ given by the expansion
	$$
	\begin{aligned}
		0
		=
		\left\langle
		\Omega\bigl(\Sigma_\sigma^r(\psi(\rho))\bigr)
		-
		\Omega(\Sigma_\sigma^r),
		w_i
		\right\rangle &=
		n\cot r
		\int_{\Sigma_\sigma^r}
		\psi(\rho)
		\bigl(x_i-\cos(r)\sigma_i\bigr)
		\,dA_{\mathrm{can}}
		\\
		&\qquad+
		\int_{\Sigma_\sigma^r}
		\widetilde V_i\left(
		(x,\sigma),
		\psi(\rho),
		\left|
		\nabla^{\Sigma_\sigma^r}\psi(\rho)
		\right|^2
		\right)
		dA_{\mathrm{can}}.
	\end{aligned}
	$$
	
	The first-variation formula gives
	$\widetilde V_i(\omega,0,0)
	=\partial_s\widetilde V_i(\omega,0,0)=0$,
	and the same identities hold after differentiation in $\omega$.
	Therefore, Taylor's theorem and the compactness of $\mathcal I_r$
	show that
	$$
	V_i(\omega,s,t)
	=
	-\frac{1}{n\cot r}\widetilde V_i(\omega,s,t)
	$$
	satisfies \eqref{constant.quadratic.cmc} and
	\eqref{constant.quadratic2.cmc}. The preceding expansion then gives
	\eqref{orthogonal.condition.cmc}.
	
	In the minimal case $r=\pi/2$, the map $\Omega$ defined above vanishes on the family of equators and hence degenerates. Following the construction in \cite{lucas_coda_andre}, we instead set
	$$
	\Omega(S)=\left(\int_S x_ix_j\,dA_{\mathrm{can}}\right)_{i,j}
	\in \mathbb{M}(n+1,n+1),
	$$
	defined onto the space of square matrices of dimension $n+1$. Its restriction to the family of unoriented equators
	$ \bar{L}_{\pi/2} = \{\Sigma_{[\sigma]}^{\pi/2} : \sigma\in \mathbb{S}^{n} \}$ is a smooth embedding. As before, we can assume that the graphical family satisfies the gauge~\eqref{Omegagauge} where $P$ denotes the tubular projection onto $\Omega( \bar{L}_{\pi/2})$. From this point the proof follows as before without any significant change.    
\end{proof}

\begin{prop}\label{prop:psi-estimates}
	Let $n\geq2$. After shrinking $\mathscr V_r$ if necessary,
	the following estimates hold for every $\rho\in\mathscr V_r$
	satisfying $\mathcal H\bigl(\rho,\psi(\rho)\bigr)=0$:
	\begin{equation}\label{psiH^2}
		\|\psi(\rho)\|_{H^2(\Sigma_\sigma^r)}
		\lesssim
		\|\beta(\rho)\|_{L^2(\Sigma_\sigma^r)}
		+
		\|\nabla\beta(\rho)\|_{L^2(\Sigma_\sigma^r)}.
	\end{equation}
	Moreover, for every $\sigma \in \mathbb{S}^n$ and $u\in T_\sigma \mathbb{S}^n$ unitary,  
	\begin{equation}\label{partial_upsiH^2}
		\|\partial_u\psi(\rho)\|_{H^2(\Sigma_\sigma^r)}\lesssim
		\|\beta(\rho)\|_{L^2(\Sigma_\sigma^r)}
		+
		\|\nabla\beta(\rho)\|_{L^2(\Sigma_\sigma^r)}
		+
		\|\nabla^2\beta(\rho)\|_{L^2(\Sigma_\sigma^r)}.
	\end{equation}
\end{prop}
\begin{proof}
	Fix $\rho\in\mathscr V_r$ satisfying
	$\mathcal H\bigl(\rho,\psi(\rho)\bigr)=0$.
	Applying the conformal transformation law in
	Lemma~\ref{confchangelawEL}, on each $\Sigma_\sigma^r$
	we have:
	$$
	-\mathcal H\bigl(0,\psi(\rho)\bigr)
	=
	j_r\bigl(\psi(\rho)\bigr)^{n-1}
	\left[
	h(r)\bigl(1-e^{\beta(\rho)}\bigr)
	+
	(n-1)N_\sigma^r\bigl(\psi(\rho)\bigr)(\rho)
	\right].
	$$
	
	On the other hand, the linearization of
	$\mathcal H(0,\cdot)$ at zero is the Jacobi operator
	$\mathcal {J}^r_\sigma$ of the CMC submanifold $\Sigma_{\sigma}^r$. Therefore, setting the quadratic error term: 
	$$E_\sigma\bigl(\psi(\rho)\bigr)
	=
	\mathcal {J}^r_\sigma\psi(\rho)
	-
	\mathcal H\bigl(0,\psi(\rho)\bigr),$$
	we see that $\psi(\rho)$ satisfies the following equation: 
	\begin{equation}\label{lineareqtrapsierho}
		\mathcal {J}^r_\sigma\psi(\rho)(\cdot,\sigma)
		=
		E_\sigma\bigl(\psi(\rho)\bigr)
		-
		j_r\bigl(\psi(\rho)\bigr)^{n-1}
		\left[
		h(r)\bigl(1-e^{\beta(\rho)}\bigr)
		+
		(n-1)N_\sigma^r\bigl(\psi(\rho)\bigr)(\rho)
		\right].
	\end{equation}
	
	We claim that after shrinking $\mathscr V_r$ once more if necessary, the error and its
	derivative with respect to the parameter satisfy the following uniform
	quadratic estimates 
	\begin{equation}\label{eq:quadratic-error-L2}
		\left\|E_\sigma\bigl(\psi(\rho)\bigr)\right\|_{L^2(\Sigma_\sigma^r)}
		\lesssim
		\|\psi(\rho)\|_{W^{1,4}(\Sigma_\sigma^r)}^2
		+
		\|\psi(\rho)\|_{L^\infty(\Sigma_\sigma^r)}
		\|\psi(\rho)\|_{H^2(\Sigma_\sigma^r)}.
	\end{equation}
	and, for every unit vector $u\in T_\sigma\mathbb S^n$,
	\begin{equation}\label{eq:quadratic-error-parameter-L2}
		\begin{aligned}
			\left\|\partial_uE_\sigma\bigl(\psi(\rho)\bigr)\right\|_{L^2(\Sigma_\sigma^r)}
			\lesssim{}&
			\|\psi(\rho)\|_{W^{1,4}(\Sigma_\sigma^r)}
			\|\partial_u\psi(\rho)\|_{W^{1,4}(\Sigma_\sigma^r)}
			+
			\|\partial_u\psi(\rho)\|_{L^\infty(\Sigma_\sigma^r)}
			\|\psi(\rho)\|_{H^2(\Sigma_\sigma^r)}
			\\
			&+
			\|\psi(\rho)\|_{L^\infty(\Sigma_\sigma^r)}
			\|\partial_u\psi(\rho)\|_{H^2(\Sigma_\sigma^r)}.
		\end{aligned}
	\end{equation}
	For clarity we move on with the argument and verify these estimates at the end of the proof. Applying standard elliptic estimates to equation~\eqref{lineareqtrapsierho} gives
	$$
	\begin{aligned}
		\|\psi(\rho)\|_{H^2(\Sigma_\sigma^r)}
		\lesssim{}&
		\|\beta(\rho)\|_{L^2(\Sigma_\sigma^r)}
		+
		\|\nabla\beta(\rho)\|_{L^2(\Sigma_\sigma^r)} + \|\psi(\rho)^T\|_{L^2(\Sigma_\sigma^r)}
		\\
		&+
		\|\psi(\rho)\|_{W^{1,4}(\Sigma_\sigma^r)}^2
		+
		\|\psi(\rho)\|_{L^\infty(\Sigma_\sigma^r)}
		\|\psi(\rho)\|_{H^2(\Sigma_\sigma^r)}.
	\end{aligned}
	$$
	Here $\psi(\rho)^T$ denotes the
	$L^2(\Sigma_\sigma^r)$-orthogonal projection of $\psi(\rho)$ onto
	$\ker\mathcal {J}^r_\sigma$. Lemma~\ref{kernJestmts.cmc}, together with the standard interpolation and absorption argument provide equation~\ref{psiH^2}. 
	
	Similarly, differentiating equation~\eqref{lineareqtrapsierho} and
	applying the standard elliptic estimate to $\partial_u\psi(\rho)$ gives
	$$
	\begin{aligned}
		\|\partial_u\psi(\rho)\|_{H^2(\Sigma_\sigma^r)}
		\lesssim{}&
		\|\partial_u\mathcal H\bigl(0,\psi(\rho)\bigr)\|_{L^2(\Sigma_\sigma^r)}
		+
		\|\partial_uE_\sigma\bigl(\psi(\rho)\bigr)\|_{L^2(\Sigma_\sigma^r)}
		+
		\|(\partial_u\psi(\rho))^T\|_{L^2(\Sigma_\sigma^r)}
		\\
		\lesssim{}&
		\|\beta(\rho)\|_{L^2(\Sigma_\sigma^r)}
		+
		\|\nabla\beta(\rho)\|_{L^2(\Sigma_\sigma^r)}
		+
		\|\nabla^2\beta(\rho)\|_{L^2(\Sigma_\sigma^r)}
		+
		\|(\partial_u\psi(\rho))^T\|_{L^2(\Sigma_\sigma^r)}
		\\
		&+
		\|\psi(\rho)\|_{W^{1,4}(\Sigma_\sigma^r)}
		\|\partial_u\psi(\rho)\|_{W^{1,4}(\Sigma_\sigma^r)}
		+
		\|\partial_u\psi(\rho)\|_{L^\infty(\Sigma_\sigma^r)}
		\|\psi(\rho)\|_{H^2(\Sigma_\sigma^r)}
		\\
		&+
		\|\psi(\rho)\|_{L^\infty(\Sigma_\sigma^r)}
		\|\partial_u\psi(\rho)\|_{H^2(\Sigma_\sigma^r)}.
	\end{aligned}
	$$
	Using equation~\eqref{psiH^2} and absorbing the remaining nonlinear
	terms, we obtain
	$$
	\|\partial_u\psi(\rho)\|_{H^2(\Sigma_\sigma^r)}
	\lesssim
	\|\beta(\rho)\|_{L^2(\Sigma_\sigma^r)}
	+
	\|\nabla\beta(\rho)\|_{L^2(\Sigma_\sigma^r)}
	+
	\|\nabla^2\beta(\rho)\|_{L^2(\Sigma_\sigma^r)}
	+
	\|(\partial_u\psi(\rho))^T\|_{L^2(\Sigma_\sigma^r)}.
	$$
	
	To estimate the term
	$\|(\partial_u\psi(\rho))^T\|_{L^2(\Sigma_\sigma^r)}$, we recall that
	Lemma~\ref{kernJestmts.cmc} provide functions $V_i$ such that
	$$
	\int_{\Sigma_\sigma^r}
	\psi(\rho)\bigl(x_i-\cos(r)\sigma_i\bigr)\,dA_{\mathrm{can}}
	=
	\int_{\Sigma_\sigma^r}
	V_i\left(
	(x,\sigma),
	\psi(\rho),
	\left|\nabla^{\Sigma_\sigma^r}\psi(\rho)\right|^2
	\right)\,dA_{\mathrm{can}}.
	$$
	Differentiating this identity in the direction of a unit tangent vector
	$u\in T_\sigma\mathbb S^n$ gives
	$$
	\begin{aligned}
		&\int_{\Sigma_\sigma^r}
		\bigl(x_i-\cos(r)\sigma_i\bigr)
		\partial_u\psi(\rho)\,dA_{\mathrm{can}}
		= \int_{\Sigma_\sigma^r}
		\sigma_i\langle x,u\rangle\psi(\rho)\,dA_{\mathrm{can}}
		\\
		&\qquad+
		\int_{\Sigma_\sigma^r}
		\Bigl(
		D_\omega V_i
		\bigl[
		(\cos(r)u-\langle x,u\rangle\sigma,u)
		\bigr]
		+
		D_sV_i\,\partial_u\psi(\rho)
		+
		D_tV_i\,\partial_u
		\left|\nabla^{\Sigma_\sigma^r}\psi(\rho)\right|^2
		\Bigr)\,dA_{\mathrm{can}},
	\end{aligned}
	$$
	where the arguments of $V_i$ and its derivatives are omitted. Finally,
	recalling that $V_i$ satisfies
	\eqref{constant.quadratic.cmc} and
	\eqref{constant.quadratic2.cmc}, we obtain
	$$
	\begin{aligned}
		\|(\partial_u\psi(\rho))^T\|_{L^2(\Sigma_\sigma^r)}
		&\lesssim
		\|\psi(\rho)^T\|_{L^2(\Sigma_\sigma^r)}
		+
		\|\psi(\rho)\|_{H^1(\Sigma_\sigma^r)}^2
		+
		\|\partial_u\psi(\rho)\|_{H^1(\Sigma_\sigma^r)}^2
		\\
		&\lesssim
		\|\psi(\rho)\|_{H^1(\Sigma_\sigma^r)}^2
		+
		\|\partial_u\psi(\rho)\|_{H^1(\Sigma_\sigma^r)}^2.
	\end{aligned}
	$$
	Substituting this estimate into the preceding inequality and repeating
	the same absorption argument proves
	equation~\eqref{partial_upsiH^2}.
	
	To prove equations~\eqref{eq:quadratic-error-L2} and
	\eqref{eq:quadratic-error-parameter-L2}, we follow the ideas of
	White~\cite{White87Parametric}. Notice that, by
	equation~\eqref{H(0,phi)} and the identity
	$\mathcal J^r_\sigma\Phi=-\Delta_{\Sigma_\sigma^r}\Phi-(n-1)(1+\cot^2r)\Phi$,
	the error can be decomposed as
	$$
	E_\sigma\bigl(\psi(\rho)\bigr)=-Q_r\left(\psi(\rho),\left|\nabla^{\Sigma_\sigma^r}\psi(\rho)\right|^2\right)-\widehat G_r\left(\psi(\rho),\nabla^{\Sigma_\sigma^r}\psi(\rho),(\nabla^{\Sigma_\sigma^r})^2\psi(\rho)\right).
	$$
	
	Here $Q_r$ and $\widehat G_r$ are given by
	$$
	\begin{aligned}
		&Q_r(s,t)
		=
		j_r(s)^{n-3}j_r'(s)
		\frac{(n-1)j_r(s)^2+(n-2)t}{\sqrt{j_r(s)^2+t}}
		-
		h(r)j_r(s)^{n-1}
		+
		(n-1)\bigl(1+\cot^2r\bigr)s,
		\\
		&G_r(s,t)
		=
		\frac{j_r(s)^{n-2}}{\sqrt{j_r(s)^2+t}}-1,
		\\
		&\widehat G_r(s,L,A)
		=
		-G_r(s,|L|^2)\mathrm{tr}A
		-
		\partial_sG_r(s,|L|^2)|L|^2
		-
		2\partial_tG_r(s,|L|^2)A(L,L).
	\end{aligned}
	$$
	
	In these formulas, $s,t\in\mathbb R$, $L\in\mathbb R^{n+1}$, and
	$A\in\mathrm{Sym}^2(\mathbb R^{n+1})$. At each
	$x\in\Sigma_\sigma^r$, the intrinsic Hessian is viewed as an element
	of $\mathrm{Sym}^2(\mathbb R^{n+1})$ by extending to zero on the orthogonal of $T_x\Sigma_\sigma^r$. With this convention, its trace is the intrinsic Laplacian.
	
	Differentiating $E_\sigma$, we also obtain the following representation:
	$$
	\begin{aligned}
		\partial_uE_\sigma\bigl(\psi(\rho)\bigr)
		={}&
		-\partial_sQ_r\,
		\partial_u\psi(\rho)
		-
		2\partial_tQ_r
		\left\langle
		\nabla^{\Sigma_\sigma^r}\psi(\rho),
		\nabla^{\Sigma_\sigma^r}\partial_u\psi(\rho)
		\right\rangle
		\\
		&-
		\partial_s\widehat G_r\,
		\partial_u\psi(\rho)
		-
		D_L\widehat G_r
		\left[
		\nabla^{\Sigma_\sigma^r}\partial_u\psi(\rho)
		\right]
		-
		D_A\widehat G_r
		\left[
		(\nabla^{\Sigma_\sigma^r})^2\partial_u\psi(\rho)
		\right].
	\end{aligned}
	$$
	
	We now note that, at the origin,
	$$
	\begin{gathered}
		Q_r(0,0)=G_r(0,0)=0,
		\qquad
		\partial_sQ_r(0,0)=0,
		\qquad
		\partial_tQ_r(0,0)=\frac{n-3}{2}\cot r,
		\\
		\partial_sG_r(0,0)=(n-3)\cot r,
		\qquad
		\partial_tG_r(0,0)=-\frac12.
	\end{gathered}
	$$
	
	Thus, in a sufficiently small neighborhood of the origin,
	$$
	\left|G_r(s,|L|^2)\right|
	\lesssim
	|s|+|L|^2,\,
	\left|\partial_sG_r(s,|L|^2)\right|
	+
	\left|\partial_tG_r(s,|L|^2)\right|
	\lesssim 1,\,
	\left|Q_r(s,|L|^2)\right|
	\lesssim
	|s|^2+|L|^2.
	$$
	
	After shrinking $\mathscr V_r$ if necessary, the previous decomposition of the error function gives the pointwise estimate
	$$
	\left|E_\sigma\bigl(\psi(\rho)\bigr)\right|
	\lesssim
	|\psi(\rho)|^2
	+
	\left|\nabla^{\Sigma_\sigma^r}\psi(\rho)\right|^2
	+
	|\psi(\rho)|
	\left|(\nabla^{\Sigma_\sigma^r})^2\psi(\rho)\right|.
	$$
	Consequently, Hölder's inequality gives
	$$
	\left\|E_\sigma\bigl(\psi(\rho)\bigr)\right\|_{L^2(\Sigma_\sigma^r)}
	\lesssim
	\|\psi(\rho)\|_{W^{1,4}(\Sigma_\sigma^r)}^2
	+
	\|\psi(\rho)\|_{L^\infty(\Sigma_\sigma^r)}
	\|\psi(\rho)\|_{H^2(\Sigma_\sigma^r)},
	$$
	which is equation~\eqref{eq:quadratic-error-L2}.
	
	Similarly, near the origin,
	$$
	\left|\partial_sQ_r(s,|L|^2)\right|
	\lesssim
	|s|+|L|^2,
	\qquad
	\left|\partial_tQ_r(s,|L|^2)\right|
	\lesssim 1.
	$$
	The second derivatives of $G_r$ are also bounded there. Combining these
	estimates with the chain-rule formula above and shrinking $\mathscr V_r$
	if necessary, we obtain
	$$
	\begin{aligned}
		\left|\partial_uE_\sigma\bigl(\psi(\rho)\bigr)\right|
		\lesssim{}&
		\left(
		|\psi(\rho)|
		+
		\left|\nabla^{\Sigma_\sigma^r}\psi(\rho)\right|
		\right)
		\left(
		|\partial_u\psi(\rho)|
		+
		\left|\nabla^{\Sigma_\sigma^r}\partial_u\psi(\rho)\right|
		\right)
		\\
		&+
		|\partial_u\psi(\rho)|
		\left|
		\left(\nabla^{\Sigma_\sigma^r}\right)^2\psi(\rho)
		\right|
		\\
		&+
		\left(
		|\psi(\rho)|
		+
		\left|\nabla^{\Sigma_\sigma^r}\psi(\rho)\right|^2
		\right)
		\left|
		\left(\nabla^{\Sigma_\sigma^r}\right)^2
		\partial_u\psi(\rho)
		\right|.
	\end{aligned}
	$$
	Using the uniform $C^3$ control of $\psi(\rho)$, together with Hölder's inequality and interpolation, we obtain the desired bound. 
\end{proof}

We now have all the ingredients needed to prove
Proposition~\ref{prop:nonlinearalln}.

\begin{proof}[Proof of Proposition~\ref{prop:nonlinearalln}]
	Lemma~\ref{lem:pointwise-Rr-estimates} and
	Proposition~\ref{prop:psi-estimates} gives
	$$
	\begin{aligned}
		\|R_r(\rho) + c_r(\rho)\|_{H^1(\mathbb S^n)}^2
		\lesssim
		\int_{\mathbb S^n}
		\sum_{m=0}^2\Bigg[
		\|\nabla^m\rho\|_{L^2(\Sigma_\sigma^r)}^4
		+
		\|\nabla^m\beta(\rho)\|_{L^2(\Sigma_\sigma^r)}^4\Bigg]
		+
		\sum_{m=0}^1
		\|\nabla^m\rho\|_{L^2(\Omega_\sigma^r)}^4
		dV_{\mathrm{can}}(\sigma).
	\end{aligned}
	$$
	
	An application of Lemma~\ref{lem:perturbed.geometric.integral.formula} to $f=|\nabla^m\rho|$, yields, after shrinking $\mathscr V_r$ if necessary, to
	$$
	\int_{\mathbb S^n}
	\|\nabla^m\beta(\rho)\|_{L^2(\Sigma_\sigma^r)}^4\,
	dV_{\mathrm{can}}(\sigma)
	\leq
	C\|\nabla^m\rho\|_{L^4(\mathbb S^n)}^4,
	\qquad m=0,1,2.
	$$
	
	Combining these estimates with the integral geometric formula in Appendix~\ref{app:geometric.integral.formula} we obtain
	$$
	\|R_r(\rho) + c_r(\rho)\|_{H^1(\mathbb S^n)}^2
	\leq
	C\sum_{m=0}^2
	\|\nabla^m\rho\|_{L^4(\mathbb S^n)}^4
	\leq
	C\|\rho\|_{W^{2,4}(\mathbb S^n)}^4.
	$$
	As claimed. 
\end{proof}
\section{The inverse problem}

In what follows we prove Theorem~\ref{maintheoremC}. Let $(M^n,\bar g)$, $n\geq 2$, be a connected Riemannian manifold of
constant sectional curvature $\kappa\in\mathbb R$, such that
$$
r^*
:=
\inf_{p\in M}\mathrm{inj}_{\bar g}(p)
>0.
$$
Fix $0<r<r^*$ (if $\kappa>0$, assume also $r<{\pi}/{2\sqrt{\kappa}}$). The conditions over $r$ guarantee that every geodesic sphere is CMC but not minimal. We recall, once more, that the minimal case was studied in~\cite{lucas_coda_andre}.   

From this point forward we fix $(M,\bar {g})$ and $0< r< r^*$ satisfying the above hypothesis. The first step into the proof is to describe our geodesic spheres as the preimage of a suitable family of functions. 

\begin{lem}\label{lem:radial-potential}
	There exist $\varepsilon>0$ and a smooth function $ \rho:\mathbb R\longrightarrow\mathbb R $ such that:
	
	\begin{enumerate}
		\item\label{lem:radial-potential1}
		$
		\rho'(t)>0,\; 
		\text{for every }t\in(0,r+\varepsilon);
		$
		
		\item\label{lem:radial-potential2} The function
		$$
		\phi:M\times M\longrightarrow\mathbb R,
		\qquad
		\phi(p,\sigma)
		:=
		\rho\bigl(d_{\bar g}(p,\sigma)\bigr),
		$$
		is smooth;
		
		\item\label{lem:radial-potential3} Setting $ \phi(\sigma)(p):=\phi(p,\sigma) $, it holds true that
		$$
		\mathrm{Hess}_{\bar g}\phi(\sigma)
		=
		\bigl(a+b \rho(t)\bigr)\bar g,
		$$
		where: 
		$$
		(a,b)=
		\begin{cases}
			(0,-\kappa), & \text{if } \kappa=-1,\\
			(1,0),      & \text{if } \kappa=0,\\
			(0,-\kappa), & \text{if } \kappa=1,
		\end{cases}
		$$
		whenever $ t=d_{\bar g}(p,\sigma)<r+\varepsilon$. 
	\end{enumerate}
\end{lem}

\begin{proof}
	As is classical from the theory of space forms, we define $\rho:[0,r+ 2\epsilon) \to \mathbb{R}$, for $\epsilon >0$ such that $r+ 4\epsilon < r^*$, as: 
	$$
	\rho(t)
	=
	\begin{cases}
		-\cos(\sqrt{\kappa}\,t),
		& \kappa>0,\\[0.4em]
		\dfrac{t^2}{2},
		& \kappa=0,\\[0.6em]
		\cosh(\sqrt{-\kappa}\,t),
		& \kappa<0.
	\end{cases}
	$$
	And we extend the function $\rho$ to the line via a cutoff function that is zero on $[r+ 3\epsilon, \infty)$.  
\end{proof}

Recall that throughout the paper we compute the mean curvature of the
boundary of a region with respect to the inward pointing unit normal
vector field. Under this convention, we collect the following facts about the geometry of geodesic spheres.

\begin{lem}\label{lem:geodesic-sphere-geometry}
	Let $\rho$ and $\phi$ as in Lemma~\ref{lem:radial-potential}. For every $\sigma\in M$, the following
	properties hold.
	
	\begin{enumerate}
		\item\label{lem:geodesic-sphere-geometry1} The $\bar g$-geodesic sphere of radius $r$ centered at $\sigma$
		is contained in a level set of $\phi(\sigma)$, namely
		$$
		\Sigma_\sigma^r 
		=
		\phi(\sigma)^{-1}(c) \cap \{p : |d_{\bar{g}}(p,\sigma) - r|< \epsilon  \},
		\qquad c:=\rho(r).
		$$
		
		\item\label{lem:geodesic-sphere-geometry2} Given $p\in M$ and $w\in \mathbb{S}_{\bar{g}}(T_pM)$, set
		$\sigma=\exp_p^{\bar g}(rw)$. Then
		$$
		\nabla^{\bar g}\phi(\sigma)(p)
		=
		-\rho'(r)w.
		$$
		
		\item\label{lem:geodesic-sphere-geometry3} The mean curvature of $\Sigma_\sigma^r$ with respect to
		$\bar g$ is given by
		$$
		h^{*} \doteq H_{(\Sigma_\sigma^r,\bar g)}
		=
		(n-1)\frac{(a+b \rho(r))}{\rho'(r)}.
		$$
		In particular, $h^*>0$. 
	\end{enumerate}
\end{lem}
\begin{proof}
	The first assertion follows immediately from the definition of
	$\phi(\sigma)$ and the fact that $\rho$ is strictly increasing in a
	neighborhood of $r$. For the second assertion, let
	$\gamma(t)=\exp_p^{\bar g}(tw)$. Since $\gamma(r)=\sigma$, the Gauss Lemma gives
	$\nabla^{\bar g}d_\sigma(p)=-w$. Hence
	$$
	\nabla^{\bar g}\phi(\sigma)(p)
	=
	\rho'(r)\nabla^{\bar g}d_\sigma(p)
	=
	-\rho'(r)w.
	$$
	
	Finally, along $\Sigma_\sigma^r$ we have
	$|\nabla^{\bar g}\phi(\sigma)|_{\bar g}=\rho'(r)$ and
	$
	\mathrm{Hess}_{\bar g}\phi(\sigma)
	=
	\bigl(a+b \rho(r)\bigr)\bar g.
	$ Which immediately provides the desired result.
\end{proof}

Inspired by the arguments of~\cite{Hangan96} and \cite{lucas_survey26}, we define the compatibility three-tensor $\mathcal{T}_g: \mathfrak{X}(M)\times \mathfrak{X}(M) \times \mathfrak{X}(M) \to \mathbb{R}$, by: 
$$
\mathcal T_g(X,Y,Z) 
:= 
g\bigl(\nabla^g_XY-\bar\nabla_XY,Z\bigr). 
$$ 
Under this notation, the fundamental equation that a metric satisfying the hypotheses of Theorem~\ref{maintheoremC} has to satisfy is given by the following proposition.

\begin{prop}\label{prop:fundamental-equation}
	Let $g$ be a Riemannian metric on $M$ such that
	$H_{(\Sigma_\sigma^r,g)} = \restr{f}{\Sigma_\sigma^r}$, for every $\sigma \in M$. Denote by
	$h^*>0$ the mean curvature of $\Sigma_\sigma^r$ with respect to
	$\bar g$. Then, for every $p\in M$ and $w\in \mathbb{S}_{\bar{g}}(T_pM)$, setting
	$\sigma_w:=\exp_p^{\bar g}(rw)$, the following equations
	hold true:
	\begin{equation}\label{eq:fundamental-cmc-identity2}
		\left[\mathcal T_g - g \otimes \mathrm{tr}^{12}_g(\mathcal T_g)\right] 
		\bigl(\nabla^g\phi(\sigma_w)(p),\nabla^g\phi(\sigma_w)(p),\nabla^g\phi(\sigma_w)(p)\bigr) 
		=0
	\end{equation} 
	and
	\begin{equation}\label{eq:fundamental-cmc-identity}
		\frac{n-1}{\rho'(r)}\frac{f(p)}{h^*}
		\left\|\nabla^g\phi(\sigma_w)(p)\right\|_g^3
		+
		\left\|\nabla^g\phi(\sigma_w)(p)\right\|_{\bar g}^2
		-
		\mathrm{tr}_g(\bar g)(p)
		\left\|\nabla^g\phi(\sigma_w)(p)\right\|_g^2
		=0.
	\end{equation}
\end{prop}
\begin{proof}
	Given an arbitrary function $V\in C^\infty(M)$, 
	$$
	\mathrm{Hess}_gV(X,Y) 
	= 
	X(YV)-(\nabla^g_XY)(V). 
	$$ 
	Therefore, the difference of the Hessian with respect to the metric $g$ and $\bar g$ can be computed via the tensor $\mathcal{T}_g$ as:    
	$$
	\mathcal T_g(X,Y,\nabla^gV)  
	= 
	\mathrm{Hess}_{\bar g}V(X,Y) 
	-\mathrm{Hess}_gV(X,Y) 
	. 
	$$
	
	Moreover, if $c$ is a regular value of $V$ and 
	$\Sigma=V^{-1}(c)$, its mean curvature can be written as 
	$$
	H_{(\Sigma,g)} 
	= 
	\frac{1}{|\nabla^gV|_g} 
	\left( 
	\Delta_gV-\mathrm{Hess}_gV(N,N) 
	\right). 
	$$ 
	From the previous relations we obtain the following equation along $\Sigma$: 
	$$
	\begin{aligned} 
		&\mathcal T_g 
		\bigl(\nabla^g V,\nabla^g V,\nabla^g V\bigr) 
		- 
		\|\nabla^g V\|_g^2\, 
		\mathrm{tr}^{12}_g(\mathcal T_g) 
		\bigl(\nabla^g V\bigr) \\ 
		&\qquad = 
		H_{(\Sigma,g)}\|\nabla^g V\|_g^3 
		+ 
		\mathrm{Hess}_{\bar g}V 
		\bigl(\nabla^g V,\nabla^g V\bigr) 
		- 
		\|\nabla^g V\|_g^2\, 
		\mathrm{tr}_g 
		\bigl(\mathrm{Hess}_{\bar g}V\bigr). 
	\end{aligned} 
	$$
	
	Given $p\in M$, $w\in \mathbb{S}_{\bar{g}}(T_pM)$, and 
	$\sigma_w=\exp_p^{\bar g}(rw)$, we specialize the previous formula to 
	$V=\phi(\sigma_w)$. Setting 
	$X_w:=\nabla^g\phi(\sigma_w)(p)$, and recalling that $H_{(\Sigma_v^r,g)}(p)=f(p)$, we obtain 
	\begin{equation}
		\label{eq:fundamental-specialized}
		\begin{aligned} 
			&\mathcal T_g(X_w,X_w,X_w) 
			- 
			\|X_w\|_g^2\mathrm{tr}^{12}_g(\mathcal T_g)(X_w) \\ 
			&\qquad = 
			f(p)\|X_w\|_g^3 
			+ 
			\mathrm{Hess}_{\bar g}\phi(\sigma_w)(X_w,X_w) 
			- 
			\|X_w\|_g^2 
			\mathrm{tr}_g 
			\bigl(\mathrm{Hess}_{\bar g}\phi(\sigma_w)\bigr). 
		\end{aligned} 
	\end{equation}
	
	From item~\ref{lem:radial-potential3} of Lemma~\ref{lem:radial-potential} the right-hand side can be simplified to
	\begin{equation}
		\label{eq:fundamental-even-part}
		f(p)\|X_w\|_g^3 
		+ 
		(a+b \rho(r))\|X_w\|_{\bar g}^2 
		- 
		(a+b \rho(r)) 
		\mathrm{tr}_g(\bar g)\|X_w\|_g^2. 
	\end{equation}
	
	To finish the proof we show that equation~\eqref{eq:fundamental-specialized} splits into two independent equations. Indeed, item~\ref{lem:geodesic-sphere-geometry2} of the Lemma~\ref{lem:geodesic-sphere-geometry} implies that for any Riemannian metric 
	$$
	\nabla^g\phi(\sigma_{-w})(p) 
	= 
	-\nabla^g\phi(\sigma_w)(p), 
	$$ 
	and hence $X_{-w}=-X_w$. Substituting $-w$ into equation~\eqref{eq:fundamental-specialized}, we observe that the 
	left-hand side is odd in $X_w$, while its right-hand side is even, as can be seen from its
	simplified form~\eqref{eq:fundamental-even-part}. Consequently, both sides vanish individually. 
	
	Since  $
	h^* 
	= 
	\frac{n-1}{\rho'(r)}(a+b \rho(r)), 
	$  we obtain from equation~\eqref{eq:fundamental-even-part} that: 
	$$
	\frac{n-1}{\rho'(r)}\frac{f(p)}{h^*}\|X_w\|_g^3 
	+ 
	\|X_w\|_{\bar g}^2 
	- 
	\mathrm{tr}_g(\bar g)\|X_w\|_g^2 
	=0, 
	$$ 
	as claimed.
\end{proof}

Recently, Guajardo~\cite{Guajardo26} extended the
tensorial approach of Hangan's ideas developed by Ambrozio,
Marques, and Neves, proving that the minimality of the canonical
family of geodesic spheres in $\mathbb{CP}^n$ or $\mathbb{HP}^n$,
$n\geq 2$, characterizes the Fubini--Study metric up to scaling.
As in Proposition~\ref{prop:fundamental-equation}, his proof is based in study the tensor $\mathcal{T}_g$ and exploits the parity properties of the minimality equations to establish rigidity.

We are set to prove Theorem~\ref{maintheoremC}. 
\begin{proof}[Proof of Theorem~\ref{maintheoremC}]
	We first prove that $g$ is conformal to $\bar{g}$. Fix $p\in M$ and let $\{e_j\}_{j=1}^n$ be a $\bar g$-orthonormal basis of
	$T_pM$ which diagonalizes $g$, that is, there exist $\lambda_j>0$ such that
	$$
	g(e_i,e_j)=\lambda_i\delta_{ij}.
	$$
	Using the notation of Proposition~\ref{prop:fundamental-equation}, set
	$$
	X_i:=\nabla^g\phi(\sigma_{e_i})(p).
	$$
	From item~\ref{lem:geodesic-sphere-geometry2} of
	Lemma~\ref{lem:geodesic-sphere-geometry}, we have
	$$
	X_i
	=
	-\frac{\rho'(r)}{\lambda_i}e_i.
	$$
	Consequently,
	$$
	\|X_i\|_g^2
	=
	\frac{\rho'(r)^2}{\lambda_i},
	\qquad
	\|X_i\|_{\bar g}^2
	=
	\frac{\rho'(r)^2}{\lambda_i^2}, \qquad \mathrm{tr}_g(\bar g)(p)
	=
	\sum_{j=1}^n\frac{1}{\lambda_j}.
	$$
	
	Returning the vector $X_i$ to equation~\eqref{eq:fundamental-cmc-identity} from Proposition~\ref{prop:fundamental-equation}, we obtain the following system of equation on the variables $\{\lambda_j\}_{j=1}^n$:
	$$
	(n-1)\frac{f(p)}{h^*}\frac{1}{\lambda_i^{3/2}}
	+
	\frac{1}{\lambda_i^2}
	-
	\left(\sum_{j=1}^n\frac{1}{\lambda_j}\right)\frac{1}{\lambda_i}
	=0,\qquad i=1,\ldots,n.
	$$
	A simple algebraic manipulation shows that
	$\lambda_i=(h^*/f(p))^2$ for every $i=1,\ldots,n$. Therefore,
	$g=e^{2u}\bar g$, for $u=\log (h^*/f)$, as desired. 
	
	Next we prove that $u$ is constant. First we notice that, by item~\ref{lem:geodesic-sphere-geometry2} of Lemma~\ref{lem:geodesic-sphere-geometry}, we have for every $w \in \mathbb{S}_{\bar{g}}(T_pM)$, that 
	$$\nabla^g\phi(\sigma_w)(p)=-\rho'(r)e^{-2u(p)}w.$$
	Therefore, from equation~\eqref{eq:fundamental-cmc-identity2}, it holds that
	$$ [\mathcal{T}_g - g \otimes \mathrm{tr}^{12}_g(\mathcal T_g)](w,w,w)=0,$$
	for every $w \in \mathbb{S}_{\bar{g}}(T_pM)$. However, from the conformal transformation law for the connection we conclude: 
	$$0 = [\mathcal{T}_g - g \otimes \mathrm{tr}^{12}_g(\mathcal T_g)](w,w,w)= (n-1)\|w\|^2_g\, g(\nabla^g u, w), $$
	for every $w \in \mathbb{S}_{\bar{g}}(T_pM)$. In particular, $u$ is constant and $g$ is homothetic to $\bar{g}$ as claimed. 
\end{proof}

\begin{rmk}
	The main tool to reach formula~\eqref{eq:fundamental-cmc-identity} is the existence of the function $\rho$ in Lemma~\ref{lem:geodesic-sphere-geometry}. Hence one may ask whether is possible to generalize the analysis given here to more general scenarios were such existence is given. However, we recall that the existence of such function for a Riemannian manifold $(M, \bar {g})$ implies constant sectional curvature.    
\end{rmk}

\appendix
\counterwithout{thm}{subsection}
\counterwithin{thm}{section}
\renewcommand{\thethm}{\thesection.\arabic{thm}}

\section{A geometric integral formula}\label{app:geometric.integral.formula}

In this appendix we record the geometric integral identities related to the Radon transform. Fix $r\in (0,\pi)$. For each $v\in \mathbb{S}^n$, let
$$
\Sigma^{r}_v=\{x\in \mathbb{S}^n : d_{\mathbb{S}^n}(x,v)=r\},
\qquad
\Omega^{r}_v=\{x\in \mathbb{S}^n : d_{\mathbb{S}^n}(x,v)\leq r\}.
$$
The following is a consequence of the integral formula for double fibration of homogeneous spaces, see for instance \cite[Section \S2, Proposition 2.2]{HelgasonBook1999}. 

\begin{lem}\label{lem:geometric.integral.formula.boundary}
	For every $f\in C^0(\mathbb{S}^n)$, the following integral formulae holds true:
	$$
	\dashint_{\mathbb{S}^n} \left( \dashint_{\Sigma^{r}_v} f \, dA_{\mathrm{can}} \right) dV_{\mathrm{can}}(v)
	= \dashint_{\mathbb{S}^n} f \, dV_{\mathrm{can}}, \qquad  \dashint_{\mathbb{S}^n} \left( \dashint_{\Omega^{r}_v} f \, dV_{\mathrm{can}} \right) dV_{\mathrm{can}}(v)
	= \dashint_{\mathbb{S}^n} f \, dV_{\mathrm{can}} .
	$$
\end{lem}

The following perturbative version of the previous formulas is obtain by a simple modification of~\cite[Lemma 9.5]{lucas_coda_andre}.
\begin{lem}\label{lem:perturbed.geometric.integral.formula}
	Fix $r\in(0,\pi)$. There exist $\varepsilon_{n,r}>0$ and
	$C_{n,r}>0$ such that the following holds. Let
	$\Phi\in C^2(\mathcal I_r)$ be such that
	$
	\|\Phi\|_{C^2(\mathcal I_r)}<\varepsilon_{n,r}.
	$
	For $f\in C^0(\mathbb S^n)$, set
	$
	f(\Phi)(x,\sigma)
	=
	f\bigl(\Sigma_\sigma^r(\Phi)(x)\bigr).
	$ Then
	$$
	\int_{\mathbb S^n}
	\|f(\Phi)(\cdot,\sigma)\|_{L^2(\Sigma_\sigma^r)}^2\,
	dV_{\mathrm{can}}(\sigma)
	\leq
	C_{n,r}\|f\|_{L^2(\mathbb S^n)}^2,
	$$
	and
	$$
	\int_{\mathbb S^n}
	\|f(\Phi)(\cdot,\sigma)\|_{L^2(\Sigma_\sigma^r)}^4\,
	dV_{\mathrm{can}}(\sigma)
	\leq
	C_{n,r}\|f\|_{L^4(\mathbb S^n)}^4.
	$$
\end{lem}
\section{Distance to the Gegenbauer zeros and coercivity}
\label{appendix:gegenbauer-sublevels}

This appendix is devoted to make explicit the relation between the distance of the
spectral multipliers to zero and the distance to the zeros of the corresponding
Gegenbauer polynomials. The proof of the distance estimate is a geometric adaptation of the argument given in \cite[Lemma~14]{BrandoliniGariboldiGigante25}.

Fix $n\geq2$. For $k\geq2$, recall that
$$
\mu_k(r)
=
-\omega_{n-1}(n+1)(n-1)\sin^{n+1}(r)
\widetilde{\mu}_k(r),
$$
where
$$
\widetilde{\mu}_k(r)
=
\frac{1}{
	\lambda_k(\mathbb S^n)\binom{n+k-2}{k}
}
C_{k-2}^{\frac{n+3}{2}}(\cos r).
$$
Let
$$
Z_k
=
\left\{
r\in\left(0,\frac{\pi}{2}\right):
\mu_k(r)=0
\right\}
$$
and define
$$
\Theta_k(r)
=
\mathrm{dist}(r,Z_k).
$$
For $\lambda>0$, we denote the angular zeros of
$C_k^\lambda$ by
$$
0<\theta_{1,k}^\lambda<\cdots<\theta_{k,k}^\lambda<\pi,
\qquad
C_k^\lambda(\cos\theta_{j,k}^\lambda)=0.
$$

The following estimate for the slope of the Gegenbauer polynomial at a angular zero is classical. 

\begin{lem}\label{lem:appendix.gegenbauer.estimates}
	Let $0<\delta<\pi/4$ and
	$I_\delta=[\delta,\pi/2-\delta]$. There exist
	$k_0(\delta,n)\geq2$ and constants
	$m_{\delta,n},M_{\delta,n}>0$ such that, for every
	$k\geq k_0(\delta,n)$,
	$$
	\left|
	\frac{d}{dr}\widetilde{\mu}_k(\widehat\theta)
	\right|
	\geq
	m_{\delta,n}k^{\frac{3-n}{2}}
	$$
	for every $\widehat\theta\in Z_k\cap I_{\delta/2}$, and
	$$
	\sup_{r\in I_{\delta/2}}
	\left|
	\frac{d^2}{dr^2}\widetilde{\mu}_k(r)
	\right|
	\leq
	M_{\delta,n}k^{\frac{5-n}{2}}.
	$$
	Moreover, after increasing $k_0(\delta,n)$, every interval between
	consecutive zeros that meets $I_\delta$ has both endpoints in
	$I_{\delta/2}$.
\end{lem}

\begin{proof}
	The asymptotic formula for the slope of the angular zero associated
	with the Gegenbauer polynomials is presented in the following paper \cite{ForsterPetras90}. More precisely, uniformly for roots contained
	in $I_{\delta/2}$,
	$$
	\left|
	\frac{d}{dx}
	C_{k-2}^{\frac{n+3}{2}}(\cos\widehat\theta)
	\right|
	\geq
	c_{\delta,n}k^{\frac{n+3}{2}}.
	$$
	Since
	$\lambda_k(\mathbb S^n)\binom{n+k-2}{k}$ is of the order $k^n$,
	differentiating $\widetilde{\mu}_k$ at a root gives the first estimate.
	
	The identity
	$\frac{d}{dx}C_k^\lambda(x)=2\lambda C_{k-1}^{\lambda+1}(x)$,
	together with the classical interior estimate
	\cite[Theorem~8.21.10]{Szego75}
	$\sup_{I_{\delta/2}}|C_k^\lambda(\cos r)|
	\leq C_{\delta,\lambda}k^{\lambda-1}$, gives
	$$
	\sup_{r\in I_{\delta/2}}
	\left|
	\frac{d^2}{dr^2}
	C_{k-2}^{\frac{n+3}{2}}(\cos r)
	\right|
	\leq
	C_{\delta,n}k^{\frac{n+5}{2}}.
	$$
	Dividing by $\lambda_k(\mathbb S^n)\binom{n+k-2}{k}$ gives the second estimate.
	
	Finally, the spacing between the zeros of the Gegenbauer polynomials is controlled by the following relation, given in \cite[Theorem~6.11.1]{Szego75}
	$$
	\theta_{j+1,k-2}^{\frac{n+3}{2}}
	-
	\theta_{j,k-2}^{\frac{n+3}{2}}
	\leq
	M_n\frac{\log(k-2)}{k-2}.
	$$
	Since the right-hand side converges to zero, every interval between
	consecutive roots that meets $I_\delta$ has both endpoints in
	$I_{\delta/2}$ for all sufficiently large $k$.
\end{proof}

The preceding estimates give the desired comparison with the distance
to the zero set.

\begin{prop}\label{prop:multiplier.distance.to.zeros}
	Let $0<\delta<\pi/4$. There exist
	$c_{\delta,n}>0$ and $k_0(\delta,n)\geq2$ such that
	$$
	|\mu_k(r)|
	\geq
	c_{\delta,n}
	k^{-\frac{n-1}{2}}
	\min\left\{
	1,
	k\Theta_k(r)
	\right\}
	$$
	for every $r\in I_\delta$ and every
	$k\geq k_0(\delta,n)$.
\end{prop}

\begin{proof}
	Apply Lemma~\ref{lem:appendix.gegenbauer.estimates} on
	$I_{\delta/2}$ and set
	$$
	S_k
	=
	m_{\delta,n}k^{\frac{3-n}{2}},
	\qquad
	B_k
	=
	M_{\delta,n}k^{\frac{5-n}{2}}.
	$$
	If $\widehat\theta\in Z_k\cap I_{\delta/2}$, then
	$|\widetilde{\mu}_k'(\widehat\theta)|\geq S_k$ and
	$\sup_{I_{\delta/2}}|\widetilde{\mu}_k''|\leq B_k$. Hence, whenever
	$|r-\widehat\theta|\leq S_k/(2B_k)$,
	$$
	|\widetilde{\mu}_k'(r)|
	\geq
	\frac{S_k}{2}
	$$
	and therefore
	$$
	|\widetilde{\mu}_k(r)|
	\geq
	\frac{S_k}{2}|r-\widehat\theta|.
	$$
	
	Choose
	$$
	\rho_{\delta,n}
	=
	\min\left\{
	\frac12,
	\frac{m_{\delta,n}}{4M_{\delta,n}}
	\right\}
	$$
	and set $\rho_k=\rho_{\delta,n}k^{-1}$. The preceding estimate is valid on the neighborhood of size $\rho_k$ of each root, therefore within each of this neighborhoods it holds true that 
	$$
	|\widetilde{\mu}_k(r)|
	\geq
	c_{\delta,n}k^{\frac{3-n}{2}}
	|r-\widehat\theta|.
	$$
	Now we argue that outside these neighborhoods of the roots we have the following lower bound: 
	$$
	|\widetilde{\mu}_k(r)|
	\geq
	c_{\delta,n}k^{\frac{1-n}{2}}
	$$
	Indeed, at the endpoints of these neighborhoods we have
	$$
	|\widetilde{\mu}_k(r)|
	\geq
	c_{\delta,n}k^{\frac{1-n}{2}}.
	$$
	Moreover, the roots of $\widetilde{\mu}_k$ are simple, and the zeros of its
	derivative interlace them. Consequently, between two consecutive
	roots the absolute value of $\widetilde{\mu}_k$ increases from zero to
	its unique maximum and then decreases to zero. Which conclude the proof of the previous affirmation.
	
	Combining the estimates inside and outside these neighborhoods gives
	$$
	|\widetilde{\mu}_k(r)|
	\geq
	c_{\delta,n}
	k^{-\frac{n-1}{2}}
	\min\left\{
	1,
	k\Theta_k(r)
	\right\}.
	$$
	Finally, since $r\in I_\delta$,
	$\sin^{n+1}r\geq\sin^{n+1}\delta$. The relation between
	$\mu_k$ and $\widetilde{\mu}_k$ therefore proves the result, after
	changing the constant $c_{\delta,n}$.
\end{proof}

\bibliographystyle{halpha}
\bibliography{refs}

\begin{thebibliography}{AMN25}
\expandafter\ifx\csname url\endcsname\relax
  \def\url#1{\texttt{#1}}\fi
\expandafter\ifx\csname doi\endcsname\relax
  \def\doi#1{\burlalt{doi:#1}{http://dx.doi.org/#1}}\fi
\expandafter\ifx\csname urlprefix\endcsname\relax\def\urlprefix{URL }\fi
\expandafter\ifx\csname href\endcsname\relax
  \def\href#1#2{#2}\fi
\expandafter\ifx\csname burlalt\endcsname\relax
  \def\burlalt#1#2{\href{#2}{#1}}\fi

\bibitem[AB22]{Asselle_Benedetti22}
Luca Asselle and Gabriele Benedetti.
\newblock Normal forms for strong magnetic systems on surfaces: trapping
  regions and rigidity of {Z}oll systems.
\newblock {\em Ergodic Theory Dynam. Systems}, 42(6):1871--1897, 2022.
\newblock \doi{10.1017/etds.2021.11}.

\bibitem[AL20]{Asselle_Lange20}
Luca Asselle and Christian Lange.
\newblock On the rigidity of {Z}oll magnetic systems on surfaces.
\newblock {\em Nonlinearity}, 33(7):3173--3194, 2020.
\newblock \doi{10.1088/1361-6544/ab839c}.

\bibitem[All72]{Allard72}
William~K. Allard.
\newblock On the first variation of a varifold.
\newblock {\em Ann. of Math. (2)}, 95(3):417--491, 1972.
\newblock \doi{10.2307/1970868}.

\bibitem[Amb26]{lucas_survey26}
Lucas Ambrozio.
\newblock Spheres with minimal equators.
\newblock {\em S\~{a}o Paulo J. Math. Sci.}, 20(2):Paper No. 25, 31, 2026.
\newblock \doi{10.1007/s40863-026-00548-0}.

\bibitem[AMN24]{AmbrozioMarquesNeves24}
Lucas Ambrozio, Fernando~C. Marques, and Andr\'{e} Neves.
\newblock Rigidity theorems for the area widths of {R}iemannian manifolds.
\newblock arXiv preprint, 2024,
  \burlalt{2408.14375}{http://arxiv.org/abs/2408.14375}.

\bibitem[AMN25]{lucas_coda_andre}
Lucas Ambrozio, Fernando~C. Marques, and Andr\'{e} Neves.
\newblock Riemannian metrics on the sphere with {Z}oll families of minimal
  hypersurfaces.
\newblock {\em J. Differential Geom.}, 130(2):269--341, 2025.
\newblock \doi{10.4310/jdg/1747156792}.

\bibitem[AP11]{AliasPiccione11}
Luis~J. Al\'{\i}as and Paolo Piccione.
\newblock On the manifold structure of the set of unparameterized embeddings
  with low regularity.
\newblock {\em Bull. Braz. Math. Soc. (N.S.)}, 42(2):171--183, 2011.
\newblock \doi{10.1007/s00574-011-0009-4}.

\bibitem[BGG25]{BrandoliniGariboldiGigante25}
Luca Brandolini, Bianca Gariboldi, and Giacomo Gigante.
\newblock Irregularities of distribution on two-point homogeneous spaces.
\newblock In {\em The mathematical heritage of Guido Weiss}, Appl. Numer.
  Harmon. Anal., pages 101--125. Birkh\"{a}user/Springer, Cham, 2025.
\newblock \doi{10.1007/978-3-031-76793-7\_5}.

\bibitem[BM06]{BurnsMatveev06}
Keith Burns and Vladimir~S. Matveev.
\newblock On the rigidity of magnetic systems with the same magnetic geodesics.
\newblock {\em Proc. Amer. Math. Soc.}, 134(2):427--434, 2006.
\newblock \doi{10.1090/S0002-9939-05-08196-7}.

\bibitem[BMS24]{BilykMastrianniSteinerberger24}
Dmitriy Bilyk, Michelle Mastrianni, and Stefan Steinerberger.
\newblock Single radius spherical cap discrepancy via gegenbadly approximable
  numbers.
\newblock {\em Adv. Math.}, 452:Paper No. 109812, 18, 2024.
\newblock \doi{10.1016/j.aim.2024.109812}.

\bibitem[BP26]{BohrPaternain26}
Jan Bohr and Gabriel~P. Paternain.
\newblock {Zoll} magnetic structures and ruled surfaces, 2026,
  \burlalt{2606.25173}{http://arxiv.org/abs/2606.25173}.
\newblock \doi{10.48550/arXiv.2606.25173}.

\bibitem[CM22]{Choi_Mantoulidis22}
Kyeongsu Choi and Christos Mantoulidis.
\newblock Ancient gradient flows of elliptic functionals and {Morse} index.
\newblock {\em Amer. J. Math.}, 144(2):541--573, 2022.
\newblock \doi{10.1353/ajm.2022.0010}.

\bibitem[CZ21]{ChengZhou21}
Da~Rong Cheng and Xin Zhou.
\newblock Existence of curves with constant geodesic curvature in a
  {Riemannian} 2-sphere.
\newblock {\em Transactions of the American Mathematical Society},
  374(12):9007--9028, 2021.
\newblock \doi{10.1090/tran/8510}.

\bibitem[CZ23]{ChengZhou23}
Da~Rong Cheng and Xin Zhou.
\newblock Existence of constant mean curvature 2-spheres in {Riemannian}
  3-spheres.
\newblock {\em Communications on Pure and Applied Mathematics},
  76(11):3374--3436, 2023.
\newblock \doi{10.1002/cpa.22114}.

\bibitem[FP90]{ForsterPetras90}
Klaus-J\"{u}rgen F\"{o}rster and Knut Petras.
\newblock On estimates for the weights in gaussian quadrature in the
  ultraspherical case.
\newblock {\em Math. Comp.}, 55(191):243--264, 1990.
\newblock \doi{10.2307/2008803}.

\bibitem[Fun13]{Funk13}
P.~Funk.
\newblock \"{U}ber {F}l\"{a}chen mit lauter geschlossenen geod\"{a}tischen
  {L}inien.
\newblock {\em Math. Ann.}, 74(2):278--300, 1913.
\newblock \doi{10.1007/BF01456044}.

\bibitem[GM20]{GalvezMira20}
Jos\'{e}~A. G\'{a}lvez and Pablo Mira.
\newblock Uniqueness of immersed spheres in three-manifolds.
\newblock {\em J. Differential Geom.}, 116(3):459--480, 2020.
\newblock \doi{10.4310/jdg/1606964415}.

\bibitem[Gre63]{Green63}
L.~W. Green.
\newblock Auf wiedersehensfl{\"a}chen.
\newblock {\em Annals of Mathematics}, 78(2):289--299, 1963.
\newblock \doi{10.2307/1970344}.

\bibitem[Gua26]{Guajardo26}
Diego Guajardo.
\newblock Oriented {Zoll} families of minimal hypersurfaces in
  {$\mathbb{CP}^n$} and {$\mathbb{HP}^n$}, 2026,
  \burlalt{2610.25379}{http://arxiv.org/abs/2610.25379}.
\newblock \urlprefix\url{https://arxiv.org/abs/2610.25379}.

\bibitem[Gui76]{Guillemin}
V.~Guillemin.
\newblock The {R}adon transform on {Z}oll surfaces.
\newblock {\em Advances in Math.}, 22(1):85--119, 1976.
\newblock \doi{10.1016/0001-8708(76)90139-0}.

\bibitem[Han96]{Hangan96}
T.~Hangan.
\newblock On the {R}iemannian metrics in $\mathbb{R}^n$ which admit all
  hyperplanes as minimal hypersurfaces.
\newblock {\em Journal of Geometry and Physics}, 18(4):326--334, 1996.
\newblock \doi{10.1016/0393-0440(95)00007-0}.

\bibitem[Hat83]{Hatcher83}
Allen~E. Hatcher.
\newblock A proof of the {S}male conjecture, {${\rm Diff}(S\sp{3})\simeq {\rm
  O}(4)$}.
\newblock {\em Ann. of Math. (2)}, 117(3):553--607, 1983.
\newblock \doi{10.2307/2007035}.

\bibitem[Hel99]{HelgasonBook1999}
Sigurdur Helgason.
\newblock {\em The {R}adon transform}, volume~5 of {\em Progress in
  Mathematics}.
\newblock Birkh\"{a}user Boston, Inc., Boston, MA, second edition, 1999.
\newblock \doi{10.1007/978-1-4757-1463-0}.

\bibitem[Her12]{Herreros12}
Pilar Herreros.
\newblock Scattering boundary rigidity in the presence of a magnetic field.
\newblock {\em Comm. Anal. Geom.}, 20(3):501--528, 2012.
\newblock \doi{10.4310/CAG.2012.v20.n3.a3}.

\bibitem[HR26]{HanRahimi26}
Rui Han and Yaghoub Rahimi.
\newblock On the small denominator problem for generalized minkowski--funk
  transforms, 2026, \burlalt{2601.09547}{http://arxiv.org/abs/2601.09547}.
\newblock \urlprefix\url{https://arxiv.org/abs/2601.09547}.

\bibitem[Kiy87]{Kiyohara87}
Kazuyoshi Kiyohara.
\newblock Riemannian metrics with periodic geodesic flows on projective spaces.
\newblock {\em Japan. J. Math. (N.S.)}, 13(2):209--234, 1987.
\newblock \doi{10.4099/math1924.13.209}.

\bibitem[LM02]{LeBrun_Mason02}
Claude LeBrun and L.~J. Mason.
\newblock Zoll manifolds and complex surfaces.
\newblock {\em J. Differential Geom.}, 61(3):453--535, 2002.
\newblock \urlprefix\url{http://projecteuclid.org/euclid.jdg/1090351530}.

\bibitem[LN71]{LeungNomizu71}
Dominic~S. Leung and Katsumi Nomizu.
\newblock The axiom of spheres in {R}iemannian geometry.
\newblock {\em J. Differential Geometry}, 5:487--489, 1971.
\newblock \urlprefix\url{http://projecteuclid.org/euclid.jdg/1214430010}.

\bibitem[LS47]{LusternikSchnirelmann47}
L.~Lusternik and L.~Schnirelmann.
\newblock Topological methods in variational problems and their application to
  the differential geometry of surfaces.
\newblock {\em Uspehi Matem. Nauk (N.S.)}, 2:166--217, 1947.

\bibitem[MS18]{MazzucchelliSuhr18}
Marco Mazzucchelli and Stefan Suhr.
\newblock A characterization of {Z}oll {R}iemannian metrics on the 2-sphere.
\newblock {\em Bull. Lond. Math. Soc.}, 50(6):997--1006, 2018.
\newblock \doi{10.1112/blms.12200}.

\bibitem[MZ25]{MazurowskiZhou25}
Liam Mazurowski and Xin Zhou.
\newblock Pairs of embedded spheres with pinched prescribed mean curvature,
  2025, \burlalt{2511.08228}{http://arxiv.org/abs/2511.08228}.
\newblock \doi{10.48550/arXiv.2511.08228}.

\bibitem[Rub00]{Rubin00}
Boris Rubin.
\newblock Generalized minkowski--funk transforms and small denominators on the
  sphere.
\newblock {\em Fract. Calc. Appl. Anal.}, 3(2):177--203, 2000.

\bibitem[Son23]{ASong23}
Antoine Song.
\newblock Existence of infinitely many minimal hypersurfaces in closed
  manifolds.
\newblock {\em Ann. of Math. (2)}, 197(3):859--895, 2023.
\newblock \doi{10.4007/annals.2023.197.3.1}.

\bibitem[SS26]{SarnataroStrykerCurves}
Lorenzo Sarnataro and Douglas Stryker.
\newblock Existence of closed embedded curves of constant curvature via
  min-max.
\newblock {\em Journal f{\"u}r die reine und angewandte Mathematik},
  831:185--211, 2026.
\newblock \doi{10.1515/crelle-2025-0090}.

\bibitem[Str81]{Strichartz81}
Robert~S. Strichartz.
\newblock {$L^p$ estimates for Radon transforms in {E}uclidean and
  non-{E}uclidean spaces}.
\newblock {\em Duke Mathematical Journal}, 48(4):699 -- 727, 1981.
\newblock \doi{10.1215/S0012-7094-81-04839-0}.

\bibitem[Str83]{Strichartz83}
Robert~S. Strichartz.
\newblock Analysis of the laplacian on the complete riemannian manifold.
\newblock {\em J. Functional Analysis}, 52(1):48--79, 1983.
\newblock \doi{10.1016/0022-1236(83)90090-3}.

\bibitem[SY94]{Schoen_Yau94}
R.~Schoen and S.-T. Yau.
\newblock {\em Lectures on differential geometry}.
\newblock Conference Proceedings and Lecture Notes in Geometry and Topology, I.
  International Press, Cambridge, MA, 1994.
\newblock Lecture notes prepared by Wei Yue Ding, Kung Ching Chang [Gong Qing
  Zhang], Jia Qing Zhong and Yi Chao Xu, Translated from the Chinese by Ding
  and S. Y. Cheng, With a preface translated from the Chinese by Kaising Tso.

\bibitem[Sze75]{Szego75}
G\'{a}bor Szeg\H{o}.
\newblock {\em Orthogonal Polynomials}, volume~23 of {\em American Mathematical
  Society Colloquium Publications}.
\newblock American Mathematical Society, Providence, RI, fourth edition, 1975.

\bibitem[Tri86]{Triebel86}
Hans Triebel.
\newblock Spaces of {Besov--Hardy--Sobolev} type on complete {Riemannian}
  manifolds.
\newblock {\em Arkiv f{"o}r Matematik}, 24:299--337, 1986.
\newblock \doi{10.1007/BF02384402}.

\bibitem[Vil68]{Vilenkin68}
N.~Ja. Vilenkin.
\newblock {\em Special Functions and the Theory of Group Representations},
  volume~22 of {\em Translations of Mathematical Monographs}.
\newblock American Mathematical Society, Providence, RI, 1968.

\bibitem[Whi87]{White87Parametric}
Brian White.
\newblock The space of $m$-dimensional surfaces that are stationary for a
  parametric elliptic functional.
\newblock {\em Indiana University Mathematics Journal}, 36(3):567--602, 1987.
\newblock \doi{10.1512/iumj.1987.36.36031}.

\bibitem[Whi91]{BWhite91}
Brian White.
\newblock The space of minimal submanifolds for varying {R}iemannian metrics.
\newblock {\em Indiana Univ. Math. J.}, 40(1):161--200, 1991.
\newblock \doi{10.1512/iumj.1991.40.40008}.

\bibitem[Whi16]{White16}
Brian White.
\newblock Controlling area blow-up in minimal or bounded mean curvature
  varieties.
\newblock {\em J. Differential Geom.}, 102(3):501--535, 2016.
\newblock \doi{10.4310/jdg/1456754017}.

\bibitem[Yau82]{Yau82}
Shing~Tung Yau.
\newblock Problem section.
\newblock In {\em Seminar on {D}ifferential {G}eometry}, Ann. of Math. Stud.,
  No. 102, pages 669--706. Princeton Univ. Press, Princeton, NJ, 1982.

\bibitem[Zol03]{Zoll1903}
O.~Zoll.
\newblock {Über Flächen mit Scharen geschlossener geodätischer Linien.}
\newblock {\em Mathematische Annalen}, 57:108--133, 1903.
\newblock \urlprefix\url{http://eudml.org/doc/158089}.
\newblock ISSN: 0025-5831.

\end{thebibliography}

\end{document}